\newif\ifconfver
\confverfalse     
\confvertrue        

\newif\ifplainver  
\plainverfalse

\ifplainver
    \confverfalse   
\fi
\ifconfver
     \documentclass[10pt,twocolumn,twoside]{IEEEtran}
\else
    \ifplainver
        \documentclass[11pt]{article}
        \usepackage{fullpage}
    \else
        \documentclass[11pt,draftcls,twocolumn]{IEEEtran}
    \fi
\fi

\usepackage{graphicx,algorithmic,bm,amsmath,amsthm,amssymb,color,cite,subcaption,tcolorbox,xargs,enumitem,bbm,enumitem,nicefrac, booktabs, multirow}
\usepackage{mathrsfs}
\usepackage[colorlinks, citecolor=blue]{hyperref}
\usepackage[linesnumbered,vlined,ruled,commentsnumbered]{algorithm2e}

\usepackage{makecell}

\newtheorem{lemma}{Lemma}
\newtheorem{proposition}{Proposition}
\newtheorem{theorem}{Theorem}
\newtheorem{definition}{Definition}
\newtheorem{corollary}{Corollary}

\newtheorem{assumption}{H\!}
\newtheorem{remark}{Remark}
\newtheorem{assumSGD}{SG\!}
\newtheorem{cvxSGD}{CVX\!}

\newtheorem{assumSGDvar}{CSA\!}
\newtheorem{assumTD}{TD\!}
\newtheorem{assumSA}{SA\!}
\newtheorem{assumNA}{NA\!}
\newtheorem{assumEM}{EM\!}
\newtheorem{assumVR}{VR\!}

\usepackage{pgfplots}
\usetikzlibrary{arrows,shapes,calc,tikzmark,backgrounds,matrix,decorations.markings}
\usepgfplotslibrary{fillbetween}

\pgfplotsset{compat=1.3}

\usepackage{relsize}
\tikzset{fontscale/.style = {font=\relsize{#1}}
    }

\definecolor{lavander}{cmyk}{0,0.48,0,0}
\definecolor{violet}{cmyk}{0.79,0.88,0,0}
\definecolor{burntorange}{cmyk}{0,0.52,1,0}

\definecolor{asuorange}{rgb}{1,0.699,0.0625}
\definecolor{asured}{rgb}{0.598,0,0.199}
\definecolor{asuborder}{rgb}{0.953,0.484,0}
\definecolor{asugrey}{rgb}{0.309,0.332,0.340}
\definecolor{asublue}{rgb}{0,0.555,0.836}
\definecolor{asugold}{rgb}{1,0.777,0.008}
\usepackage{ulem}
\usepackage{enumitem}

    \makeatletter
    \def\multilimits@{\bgroup
  \Let@
  \restore@math@cr
  \default@tag
 \baselineskip\fontdimen10 \scriptfont\tw@
 \advance\baselineskip\fontdimen12 \scriptfont\tw@
 \lineskip\thr@@\fontdimen8 \scriptfont\thr@@
 \lineskiplimit\lineskip
 \vbox\bgroup\ialign\bgroup\hfil$\m@th\scriptstyle{##}$\hfil\crcr}
    \def\Sb{_\multilimits@}
    \def\endSb{\crcr\egroup\egroup\egroup}
\makeatother

\newtheoremstyle{t}         
    {\baselineskip}{2\topsep}      
    {\rm}                   
    {0pt}{\bfseries}  
    {}                      
    { }                      
    {\thmname{#1}\thmnumber{#2}.}

\theoremstyle{t}

\makeatletter
\DeclareRobustCommand*\cal{\@fontswitch\relax\mathcal}
\makeatother
\usepackage{cleveref}
\crefname{assumption}{{\textbf{H}} \hspace*{-4pt}}{{\textbf{H}} \hspace*{-4pt}}
\Crefname{assumSGD}{{\textbf{SG}}\hspace*{-3pt}}{{\textbf{SG}}\hspace*{-3pt}}
\Crefname{cvxSGD}{{\textbf{CVX}}\hspace*{-3pt}}{{\textbf{CVX}}\hspace*{-3pt}}

\Crefname{assumSGDvar}{{\textbf{CSA}}\hspace*{-3pt}}{{\textbf{CSG}}\hspace*{-3pt}}
\Crefname{assumTD}{{\textbf{TD}}\hspace*{-3pt}}{{\textbf{TD}}\hspace*{-3pt}}
\Crefname{assumSA}{{\textbf{SA}}\hspace*{-3pt}}{{\textbf{SA}}\hspace*{-3pt}}
\Crefname{assumNA}{{\textbf{NA}}\hspace*{-3pt}}{{\textbf{NA}}\hspace*{-3pt}}
\Crefname{assumEM}{{\textbf{EM}}\hspace*{-3pt}}{{\textbf{EM}}\hspace*{-3pt}}
\Crefname{assumVR}{{\textbf{VR}}\hspace*{-3pt}}{{\textbf{VR}}\hspace*{-3pt}}

\Crefname{paragraph}{Section}{Section}

\Crefname{algorithm}{{Algorithm}}{Algorithm}

\usepackage[textwidth=3cm,textsize=tiny]{todonotes}

\makeatletter
\newcommand{\settitle}{\@maketitle}
\makeatother

\newcommand{\ie}{{\em i.e.,~}}

\newcommand{\eg}{{\em e.g.,~}}

\newcommand{\wrt}{{\em w.r.t.~}}

\def\cost{\operatorname{cost}}
\def\RR{\mathbb{R}}

\def\PE{\mathbb{E}}

\def\timecum{\mathrm{t}}
\def\PP{\mathbb{P}}

\def\E{\mathbf{E}}

\def\SA{\textsf{SA}}
\def\SG{\textsf{SG}}

\def\hg{{\bm h}}
\def\bias{{\bm b}}
\def\Bias{{\bm b}}
\def\Hg{{\bm{H}}}

\def\boundbias{\tau}

\def\boundstd{\sigma}
\def\boundvar{\boundstd^2}
\def\step{\gamma}
\def\Bterm{\operatorname{B}}
\def\Vinit{\mathscr{V}}
\def\State{{\bm{X}}}

\newcommand{\1}{\mathbbm{1}} 

\newcommand{\beq}{\begin{equation}}
\newcommand{\eeq}{\end{equation}}
\newcommand{\beo}{\begin{array}{rl}}
\newcommand{\eeo}{\end{array}}

\def\mcf{\mathcal{F}}
\def\mcg{\mathcal{F}}

\def\bU{\mathbf{U}}

\def\iid{i.i.d.}

\usepackage{color}

\DeclareMathOperator{\rank}{\operatorname{rk}}
\DeclareMathOperator{\diag}{\operatorname{diag}}
\DeclareMathOperator{\Proj}{Prj}

\DeclareMathOperator{\Span}{span}

\newcommand{\argmin}{\mathop{\mathrm{arg\,min}}}

\newcommand{\argmax}{\mathop{\mathrm{arg\,max}}}

\newcommand{\grd}{{\nabla}}

\def\prm{\bm{w}}
\def\dprm{{\bm{\mathsf{w}}}}
\def\dpy{\bm{\mathsf{y}}}
\def\dpY{\bm{\mathsf{Y}}}

\def\dPrm{\bm{\mathsf{W}}}
\def\dNoise{\bm{\mathsf{U}}}
\def\dBias{\bm{\mathsf{B}}}

\def\totstep{T}

\newcommand{\eqdef}{\mathrel{\mathop:}=}
\def\lyap{{\operatorname{V}}}
\newcommand{\superlyap}{{\operatorname{W}}}
\def\Liplyap{L_{\lyap}}

\def\rholyap{\varrho}
\def\clyap{c}

\def\bsf{\operatorname{b}}
\newcommand{\Lip}[1]{L_{{#1}}}
\def\Id{\operatorname{I}}
\newcommandx\sequence[3][2=,3=]
{\ifthenelse{\equal{#3}{}}{\ensuremath{\{ #1_{#2}\}}}{\ensuremath{\{ #1_{#2}, \eqsp #2 \in #3 \}}}}
\newcommand{\coint}[1]{\left[#1\right)}
\newcommand{\ocint}[1]{\left(#1\right]}
\newcommand{\ooint}[1]{\left(#1\right)}
\newcommand{\ccint}[1]{\left[#1\right]}

\def\schat_const{s}

\def\rset{\ensuremath{\mathbb{R}}}
\def\nset{\ensuremath{\mathbb{N}}}

\newcommand{\eqsp}{\;}

\newcommandx{\as}[1][1=\PP]{\ensuremath{#1\, -\mathrm{a.s.}}}

\def\PE{\mathbb{E}}

\def\bY{\mathbf{Y}}
\def\bZ{\mathbf{Z}}
\def\bW{\mathbf{W}}

\newcommandx{\indi}[2][1=]{\1^{#1}_{#2}}

\def\dist{\operatorname{d}}

\def\iff{if and only if}

\def\valuefuncgen{\mathcal{V}}
\newcommandx{\valuefunc}[1][1=]{\valuefuncgen_{#1}}
\def\policy{\pi}
\def\stateMRP{\mathcal{S}}
\def\kerMRP{\mathcal{P}}
\def\rewardMRP{\mathcal{R}}
\newcommandx{\bellman}[1][1=]{\operatorname{B}_{#1}}

\def\statdistMRP{\varpi}
\def\feature{\boldsymbol{\phi}}
\def\Feature{\boldsymbol{\Phi}}
\def\Dtd{{\mathbf{D}}}
\def\bSigma{{\mathbf{\Sigma}}}
\def\vmintd{v_{\mathrm{min}}}
\def\bX{\mathbf{X}}
\def\bx{\mathbf{x}}
\def\bv{\mathbf{v}}

\def\rmd{\mathrm{d}}
\def\noise{\mathbf{u}}
\def\dnoise{\bm{\mathsf{u}}}
\def\dbias{\bm{\mathsf{b}}}

\def\Noise{{\bm{u}}}
\def\timeinc{s}
\newcommandx{\flow}[2][1=,2=]{\Phi^{#2}_{#1}}
\def\orbit{\operatorname{Orb}}
\def\chainrecurrent{\operatorname{CR}}
\def\equilibrium{\operatorname{EQ}}
\def\limset{\operatorname{Lim}}

\newcommand{\nofrac}[2]{#1 / #2} 

\newcommandx{\diff}[1]{\operatorname{D}_{#1}}
\def\rhs{RHS}
\def\DIhg{\hg_{\operatorname{DI}}^{\boundbias_0}}
\def\SolutionDI{\mathcal{C}^{\DIhg}}
\newcommandx\ball[3][1=]{\mathrm{B}_{#1} (#2,#3)}

\def\prmo{\bm{\theta}}
\newcommandx{\QEM}[1][1=]{\mathcal{Q}^{\mathrm{EM}}_{#1}}
\newcommandx{\crossentropy}[1][1=]{\mathcal{H}_{#1}}
\newcommandx{\chunk}[3]{{#1}_{#2}^{#3}}

\def\lbatchgd{\mathsf{b}}
\def\Zset{\mathcal{Z}}
\def\sem{\mathsf{S}}
\def\barsem{\mathsf{\bar{s}}}
\def\phiem{\boldsymbol{\phi}}
\def\psiem{\boldsymbol{\psi}}
\def\mapem{\mathsf{T}}
\def\Bem{\mathsf{B}}
\def\vminem{v_{\mathrm{min}}}
\def\vmaxem{v_{\mathrm{max}}}
\def\nbrMC{\mathsf{m}}
\def\lbatchem{\mathsf{b}_{\tiny{\mathrm{EM}}}}

\def\boundvarem{\bar{\sigma}^2}
\def\contvarem{{\bf V}}
\def\MC{\mathsf{MC}}
\def\Hgrv{\Hg^{\mathrm{vr}}}
\def\lbatchrv{\mathsf{b}_{\mathrm{vr}}}
\def\batchrv{\mathcal{B}}
\def\kin{\mathrm{k_{in}}}
\def\kout{\mathrm{k_{out}}}
\def\Liphrv{L} 
\def\Deriv{\mathsf{D}}

\def\loss{\ell}

\def\omgbiased{(1-\delta_{\compressor})}
\def\oneminusomgbiased{\delta_{\compressor}}

\def\omgunbiased{\omega_{\compressor}}
\def\omguniform{\kappa_{\compressor}}

\def\linearunbiased{\Delta_{\compressor}}

\def\linearset{{\cal B}_{\compressor}}

\newcommandx{\CPEx}[3][1=]{\PE_{#1}\left[ \left. #2 \, \right| #3 \right]}
\newcommandx{\CPE}[3][1=]{\PE_{#1}^{{#3}}\bigl[#2 \, \bigr]}
\newcommandx{\CPEs}[3][1=]{\PE_{#1}^{{#3}}[ #2 ]}
\newcommandx{\pscal}[3][1=]{\langle#2\,|\,#3\rangle_{#1}}
\newcommandx{\Pscal}[3][1=]{\left\langle#2\,|\,#3\right\rangle_{#1}}
\newcommandx{\compressor}{\boldsymbol{\mathcal{C}}}

\newcommand{\revisionupdates}[1]{{\color{black} #1}}
\usepackage{multirow}

\begin{document}
\title{Stochastic Approximation Beyond Gradient for Signal Processing and Machine Learning}
\author{Aymeric Dieuleveut, Gersende Fort, Eric Moulines, Hoi-To Wai
\thanks{AD and EM are with Ecole Polytechnique, CMAP, UMR 7641, France. GF is with Institut de Mathématiques de Toulouse, UMR5216, Université de Toulouse,  CNRS; UPS, F-31062 Toulouse Cedex 9, France. HTW is with CUHK, Hong Kong. 
Work partly supported by the {\it Fondation Simone et Cino Del Duca, Institut de France}, ANR under the program MaSDOL-19-CE23-0017-01, ANR-19-CHIA-SCAI-002, HKRGC Project 24203520, Hi!Paris \textit{FLAG} project, and been carried out under the auspices of Lagrange research Center for Mathematics and Calculus. 
}
}

{\let\newpage\relax\maketitle}

\begin{abstract}
    Stochastic Approximation ({\SA}) is a classical algorithm that has had since the early days a huge impact on signal processing, and nowadays on machine learning, due to the necessity to deal with a large amount of data observed with uncertainties. An exemplar special case of {\SA} pertains to the popular stochastic (sub)gradient  algorithm which is the working horse behind many important applications. A lesser-known fact is that the {\SA} scheme also extends to non-stochastic-gradient algorithms such as compressed stochastic gradient, stochastic expectation-maximization, and a number of reinforcement learning algorithms. The aim of this article is to overview and introduce the non-stochastic-gradient perspectives of {\SA} to the signal processing and machine learning audiences through presenting a design guideline of {\SA} algorithms backed by theories. Our central theme is to propose a general framework that unifies existing theories of {\SA}, including its non-asymptotic and asymptotic convergence results, and demonstrate their applications on popular non-stochastic-gradient algorithms. We build our analysis framework based on classes of Lyapunov functions that satisfy a variety of mild conditions. We draw connections between non-stochastic-gradient algorithms and scenarios when the Lyapunov function is smooth, convex, or strongly convex. Using the said framework, we illustrate the convergence properties of the non-stochastic-gradient algorithms using concrete examples. Extensions to the emerging variance reduction techniques for improved sample complexity will also be discussed.  
\end{abstract}


\allowdisplaybreaks[4]


\section{Introduction}
Stochastic Approximation (\SA) is a classical iterative algorithm that has a long history of over 70 years \cite{robbins1951stochastic,blum1954multidimensional}. The goal of the \SA\ scheme is to determine the roots of a nonlinear system 
when the mean field
cannot be explicitly computed but a random oracle exists.
The {\SA} scheme has had since the early days a huge impact  in signal processing and automatic control: the first applications focused on the adaptive identification of systems \cite{ljung1977analysis,ljung1983theory,benveniste2012adaptive,haykin2002adaptive}.
Recently, the spectrum of use of {\SA} schemes has widened considerably with the applications to statistical machine learning; see \cite{bottou2003stochastic,bottou2010large,principe2011kernel}.



An extensive literature on stochastic optimization is devoted to the stochastic (sub)gradient (\SG) algorithm, which is by far the most popular application of the {\SA} scheme, see \eg \cite{bottou2018optimization,cevher2014convex, lan2020first} and the references therein. The stochastic gradient algorithms are characterized by having an update recursion featuring an \emph{unbiased} mean field which is the gradient of a loss function to be minimized. However, a lesser known fact is that {\SA} scheme also includes \emph{non-stochastic-gradient (non-{\SG})} algorithms whose expected oracles are not the gradients of any function and whose oracles are possibly \emph{biased} even in the asymptotic sense. Recently, these \emph{non-\SG} algorithms have gained attention in many scenarios of Signal Processing (SP) and Machine Learning (ML).
Classical examples include pre-conditioned least mean square  for linear system identification \cite{kushner2003stochastic};
other more subtle examples include  natural gradient methods \cite{ren2021tensor,gunasekar2021mirrorless}, online blind source separation \cite{cichocki1994robust,cardoso1998blind,hyvarinen2000independent,cichocki2002adaptive}, straight-through compressed gradient estimator \cite{Courbariaux2015,shekhovtsov2021reintroducing} and randomized coordinate descent algorithm \cite{nutini2015coordinate}.

The {\SA} schemes that are non-stochastic-gradient algorithms appear quite naturally in modern statistical learning. 
An example is the determination of stationary points appearing in the stochastic versions of the Majorize-Minimization (MM) methods \cite{lange:2013,sun2016majorization} such as the {\SA} versions of the Expectation-Maximization (EM) algorithm \cite{delyon1999convergence,fort2021fast,karimi2019non}.
Another example is the algorithms used in reinforcement learning, such as TD-learning and $Q$-learning, in which the iterative mapping is deduced from the Bellman operator
\cite{sutton1988learning,sutton2018reinforcement,tsitsiklis1996analysis,karimi2019non}.

To illustrate the application of the {\SA} scheme as non-stochastic-gradient algorithms, we may take a closer look at the family of stochastic EM algorithms introduced in \cite{delyon1999convergence}.
While EM can be derived by the majorization-minimization method, a powerful perspective presented in \cite{delyon1999convergence} is to view the expectation step ('E-step') as fixed-point iterations in \emph{sufficient statistics} space (the $s$-space). This perspective allows us to build highly efficient stochastic EM algorithms for streaming data and big data \cite{cappe2009line,chen2018stochastic,karimi2019non,fort2021fast,fort2020stochastic}. The challenge in their analysis is that this fixed point map in $s$-space is not the stochastic gradient map of \emph{any} function. Therefore, the stochastic EM algorithms are in fact \emph{not}  {\SG} algorithms. The convergence analysis of these algorithms to cope with the maximum likelihood estimation will be analyzed using the $\SA$ scheme, which includes the {non-\SG} algorithms.

There are many excellent overview articles or books on {\SA} scheme. We note that classical books such as \cite{kushner2003stochastic,benveniste2012adaptive, borkar2009stochastic} 
focused on asymptotic convergence for \emph{unbiased} {\SA} under a set of restrictive stepsize conditions. 
Along the line of applications on adaptive filtering, \cite{sayed2011adaptive} presents finite-time analysis on the family of adaptive filtering algorithms such as least mean squares, recursive least squares, etc. The recent articles \cite{cevher2014convex, bottou2018optimization, lan2020first} are devoted to the \emph{gradient-based} {\SA} schemes. While they provide a modern treatment on the convergence of {\SA} schemes, the discussions are limited to stochastic gradient algorithms. 
The recent book \cite{meyn2022book} includes results that are applicable to non-{\SG} algorithms, but is otherwise focused on $\SA$ schemes matched to find the roots of the gradient of a strongly convex objective function. 

Most of the existing overviews on the convergence for $\SA$ schemes are not comprehensive when it comes to discussing the results for non-stochastic-gradient algorithms, or they are limited to asymptotic convergence for algorithms with decreasing step sizes sometimes accompanied by limit distributions (in favorable cases).
A possible reason behind this is the lack of a \emph{proper} Lyapunov function to set the convergence analysis framework, and the potential bias may destabilize the {\SA} recursion.
To this end, there are no convergence results  (or only in certain cases)  for the beyond-gradient-{\SA} scheme in the literature at the level of generality of {\SG} algorithms such as \cite{lan2020first}.
Also missing is a principled guide to designing or improving algorithms for SP and ML that do not originate from Stochastic Gradient.


\begin{figure*}
    \centering
    \includegraphics[width=0.7\linewidth,trim=0 14cm 1cm 0,clip]{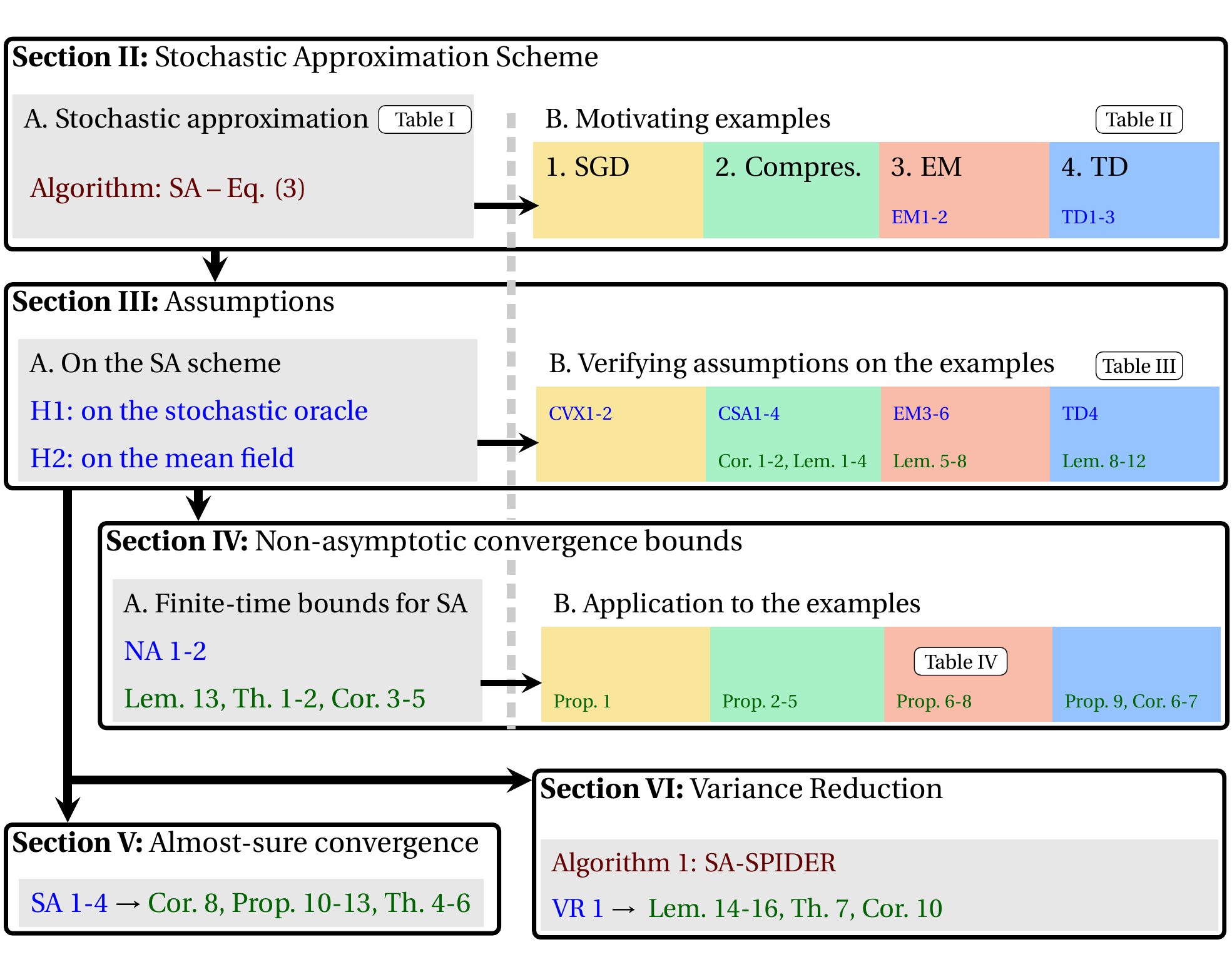}
    \caption{\revisionupdates{Guide for navigating through the content in this overview article. 
\textbf{Outline: }{\bf \Cref{sec:sa}} introduces the central~\SA~algorithm and its  derivation on four motivating examples, and \Cref{sec:assumptions_and_verif} the main two assumptions on the oracle and the field for  the analysis to hold, with corresponding derivation on examples. 
Theoretical  results are given in {\bf \Cref{sec:nonasymp}},  {\bf \Cref{sec:asymptotic-convergence}}, and {\bf \Cref{sec:vr}}. 
\textbf{Dependencies and reading guide:} 
Gray background indicates generic results on the \SA~scheme, and colored backgrounds indicate their application to the four examples. 
 In  
\Cref{sec:sa,sec:assumptions_and_verif,sec:nonasymp}, subsection A~deals with the generic scheme  and subsection B.1 to B.4~with the four examples. 
Black font indicates the sections,  blue font indicates assumptions that are made for the analysis; green font corresponds to the results (lemmas, propositions, corollaries and theorems), and red font to algorithms.  Arrows indicate dependencies:
readers may choose to read only the generic results on the \SA~scheme, together with some of the examples. 
    }
    }\vspace{-.2cm}
    \label{fig:over}
\end{figure*}

This article fills a gap in the literature by proposing a general framework that summarizes 
recent advances in the theories of the $\SA$ scheme, with emphasis on the recent applications to SP and ML. Our aim is to provide an up-to-date overview of the classical {\SA} scheme for researchers working in the fields of SP and/or ML. We shall cover results from the basic insights behind the {\SA} scheme, to standard analysis on convergence in expectation, to the advanced analysis with almost-sure convergence. Our results will handle the challenging settings of generic non-{\SG} (possibly biased) algorithms.

In the \textbf{first part} of this paper, the {\SA} scheme is introduced as a 
root-finding algorithm using a stochastic oracle, designed through Euler discretization of an ordinary differential equation flow (see \Cref{sec:SAtuto}).
Then, in \Cref{sec:ex}, examples of {\SA} schemes will be presented with a focus on non-{\SG} algorithms such as compressed {\SG} algorithms, stochastic expectation maximization, and policy evaluation via temporal difference learning. \Cref{sec:assume} discusses the general assumptions required for the convergence of the {\SA} scheme, where we formalize a set of conditions for 
a \emph{proper Lyapunov} function design with respect to the {\SA} scheme that can be potentially biased. \Cref{sec:verif:assume} shows how these conditions are satisfied in the applications listed. Readers are recommended to read this part first to get themselves familiar with the basics of {\SA} scheme. 

We next present the general convergence theories for {\SA} in two different favors. 
The \textbf{second part} is devoted to non-asymptotic convergence bounds of the {\SA} scheme that focuses on the \emph{expected convergence} towards the root(s) of a nonlinear system in a finite number of iterations. \Cref{subsec:finite-time-bounds} describes a unified result on finite-time bounds and sample complexity. One of the features of our study is to investigate the effects of the presence of bias in the stochastic oracle, a situation that often arises in applications, e.g., when the stochastic oracle uses compression or quantization, or uses Monte Carlo methods - importance sampling or Markov chain Monte Carlo that are inherently biased - see \cite{fort:etal:2016,karimi2019non,fort:moulines:2022}. In \Cref{subsec:finite-time-bounds-examples}, we illustrate the application of these results using the examples introduced in \Cref{sec:ex} and we discuss the  obtained  results with the state of the art.

The \textbf{third part} is devoted to the almost-sure convergence of the sequence of iterates generated by the {\SA} scheme towards the root(s) of a nonlinear system when the number of iterations goes to infinity. \Cref{subsec:ODE-method} gives a brief overview of the theory of Stochastic Approximation with decreasing step size.
In \Cref{subsec:limit-sets} we consider the case where the approximation bias vanishes asymptotically. We use the \emph{Ordinary Differential Equation} (ODE) method, which relates the almost sure limit sets of the stochastic approximation process to the limit sets of the flow of an ODE. 
We also establish almost-sure convergence results under the same assumptions as for the non-asymptotic approaches.
In \Cref{subsec:asymptotic-convergence-bias}, we extend these results to the case where the bias in the approximation does not vanish asymptotically. This is a case that has been much less studied in the literature. While we use \cite{tadic2017asymptotic,ramaswamy2017analysis}, which deals with stochastic-gradient-type {\SA} schemes, the results presented in this section are original.

Finally, in the \textbf{fourth part} we review a recent advance in {\SA}: we discuss in particular  the \emph{Stochastic Path-Integrated Differential EstimatoR} method (SPIDER) originally introduced for stochastic gradient algorithms \cite{fang2018spider,wang:etal:2019} and then extended to EM \cite{fort2020stochastic,dieuleveut2021federated}.  We provide here an algorithmic description and a non-asymptotic convergence analysis for a general {\SA} scheme (see also \cite{fort:moulines:2022}).

We refer the readers to Fig.~\ref{fig:over} for a guide to navigate through this overview article.

\vspace{.1cm}
\noindent \textbf{Notations.}
$\nset$ is the set of the non-negative integers, $\rset$ is the set of the real numbers, $\rset_+$ is the set of the non-negative real numbers, and $\bar \rset_+ \eqdef \rset_+ \cup \{+\infty\} $  is its completed version. 
When necessary, we use the convention $0^{-1}=\infty$ and $\infty \times 0 = 0$.
The set of the minimizers and maximizers of a function $F$ is denoted by $\argmin$  (resp. $\argmax$) when it is a singleton, and $\mathrm{Argmin}$ (resp. $\mathrm{Argmax}$) otherwise.  
For two real numbers $x$ and $y$, $x\wedge y$ (resp. $x\vee y$) denotes the minimum (resp. maximum) of $x$ and $y$. $\lfloor x \rfloor$  is the floor of $x$.

For a vector or matrix $\Dtd$,  $\Dtd^\top$ is the transpose of $\Dtd$.  Vectors are column-vectors. 
For two vectors $\prm, \prm'$ in $\rset^d$, $\pscal{\prm}{\prm'} \eqdef \prm^\top \prm'$ is the dot product of $\prm$ and $\prm'$. We set $\| \prm \| \eqdef \sqrt{\pscal{\prm}{\prm}}$. For a positive-definite matrix $\Dtd$, we denote by $\| \prm\|_{\Dtd} \eqdef \sqrt{\prm^\top \Dtd \prm}$ the norm associated to the scalar product induced by $\Dtd$.

For a differentiable function $\hg$,  its gradient is denoted by $\nabla \hg$. When $\Hg$ is a function of two variables $(\prm, {\bf x}) \mapsto \Hg(\prm, {\bf x})$, we will write $\Deriv_{ij}\Hg(\prm, {\bf x})$ for the $i$-th derivative  \wrt\ the variable $\prm$ and the $j$-th derivative \wrt\ the variable ${\bf x}$ of the function $\Hg$, evaluated at $(\prm, {\bf x})$.

All the random variables are defined on a probability space $(\Omega, \mathcal{A}, \PP)$. $\PE$ is the expectation associated to the probability $\PP$. When $\mcf$ is a filtration, $\CPE{\State}{\mcf}$ is the expectation of the random variable $\State$ conditionally to $\mcf$. When $\mcf$ is equal to $\sigma({\bf U})$, the $\sigma$-field generated by the random variable ${\bf U}$, we will write $\CPE{\State}{{\bf U}}$.
The limit set of a sequence is denoted by $\limset(\sequence{\dprm}[k][\nset])$: $\dprm_\star \in \limset(\sequence{\dprm}[k][\nset])$ if  $\lim _{k \rightarrow \infty} \dprm_{n_k}=\dprm_\star$ for some subsequence $\sequence{n}[k][\nset]$ such that  $\lim _{k \rightarrow \infty} n_k=+\infty$.

\section{Stochastic Approximation Scheme}\label{sec:sa}
\subsection{\revisionupdates{Stochastic Approximation}}
\label{sec:SAtuto}

Stochastic Approximation ({\SA}) is a class of stochastic algorithms aiming at a solution to the root-finding problem:
\beq \label{eq:rootfind}
\text{find}~~\prm^\star \in \RR^d~~\text{such that}~~\hg( \prm^\star ) = {\bm 0} \eqsp,
\eeq
where $\hg: \rset^d \to \rset^d$, $d \in \nset$, is known as the mean field function. Problem \eqref{eq:rootfind} is motivated by many tasks in SP and ML, such as convex or non-convex optimization \cite{lan2020first},
statistical estimation \cite{cappe2009line},
policy evaluation \cite{sutton2018reinforcement}, etc., some of them will be described later in \Cref{sec:ex}.

Solving \eqref{eq:rootfind} is challenging when \emph{only} stochastic estimates of the map $\hg(\cdot)$ are available. We follow the ordinary differential equation (ODE) approach introduced in \cite{ljung1977analysis,ljung1983theory}. This approach consists in identifying the limit points of the {\SA} sequence with a sub-class of invariant sets of the ODE flow (see \Cref{sec:asymptotic-convergence} for precise definitions)
\beq \label{eq:ode}
{\rmd \prm} / {\rmd  t} = \hg( \prm (t) ) \eqsp.
\eeq
In particular, the equilibrium points of the ODE (\ie the points $\prm_\star$ satisfying $\hg(\prm_\star)=0$) are the solutions of \eqref{eq:rootfind}.

To obtain a discrete-time algorithm, the common trick is to discretize the continuous time algorithm \eqref{eq:ode} through the Euler's scheme. \revisionupdates{Set $\prm_k \equiv \prm( k \gamma )$, with a sufficiently small discretization parameter $\gamma > 0$, the time-derivative of $\prm(t)$ can be approximated by the finite difference
\begin{equation} \label{eq:finite-diff}
{\rmd \prm} / {\rmd t} \approx (\prm_{k+1} - \prm_k ) / \gamma \eqsp.
\end{equation}
Substituting into \eqref{eq:ode} yields the algorithm $\prm_{k+1} = \prm_k + \gamma \hg( \prm_k )$.}
Further replacing $\hg(\prm)$ by a stochastic oracle leads to the \emph{stochastic approximation} recursion.\vspace{.1cm}

\begin{tcolorbox}[boxsep=1pt,left=4pt,right=4pt,top=3pt,bottom=3pt]
\underline{{\SA} Recursion:} let $\prm_0 \in \rset^d$ be the initialization,
\beq \label{eq:sa}
\prm_{k+1} = \prm_k + \step_{k+1} \Hg(\prm_k, \State_{k+1} ) ,k \in \nset\eqsp,
\eeq
where $\step_{k+1} > 0$ is the step size, $\State_{k+1}$ is a ${\sf X}$-valued random variable, and $\Hg: \rset^d \times {\sf X} \to \rset^d$ is a vector field that describes stochastic oracles of $\hg(\cdot)$.
\end{tcolorbox}

\revisionupdates{Before delving further into the {\SA} scheme \eqref{eq:sa}, we remark that alternatives to \eqref{eq:ode} can be used to derive other (stochastic) algorithms for solving \eqref{eq:rootfind}. For instance, the second order ODE
\begin{equation}
{\rmd^2 \prm} / {\rmd t^2} + (3/t) \, {\rmd \prm} / {\rmd t} = \hg( \prm (t) )
\end{equation}
with the Taylor approximation formulas 
\begin{equation} \notag
\frac{\prm_{k+1} - \prm_k}{\gamma} \approx \frac{\rmd \prm}{\rmd t} + 
\frac{\gamma}{2} \frac{\rmd^2 \prm}{\rmd t^2}, ~  
\frac{\prm_{k} - \prm_{k-1} }{\gamma} \approx \frac{\rmd \prm}{\rmd t} - 
\frac{\gamma}{2} \frac{\rmd^2 \prm}{\rmd t^2}
\end{equation}
lead to the Nesterov's accelerated method \cite{su2016differential}, and admit the same fixed points as \eqref{eq:ode}. Notice that an active area of research is to understand different momentum methods from an ODE perspective \cite{shi2021understanding}. However, in the stochastic setting, there are limited results even in the stochastic gradient case, see \cite{aybat2020robust}. }

\revisionupdates{Our focus is on the behavior of the {\SA} scheme \eqref{eq:sa} derived from \eqref{eq:ode}.}
To give some insights, it is useful to write the stochastic oracle as $\Hg(\prm_{k},\State_{k+1})= \hg(\prm_k)+ \noise_{k+1}$ where $\noise_{k+1}$ is a perturbation vector that distorts the mean field update direction $\hg( \prm_k )$. In the simplest setting, the conditional expectation given the past history of the algorithm evaluates to $\CPE{\Hg(\prm_k, \State_{k+1})}{\mcf_k} = \hg(\prm_k)$ for any $\prm_k \in \rset^d$; the past history $\mcf_k$ is the filtration  $\sigma(\prm_0, \State_\ell, \ell \leq k)$  up to the $k$-th iteration. \revisionupdates{For this case, the conditional expectation
\begin{equation}
    \CPE{ \prm_{k+1} }{\mcf_k} = \prm_k + \step_{k+1} \hg( \prm_k )
\end{equation} 
coincides with the deterministic algorithm obtained by discretizing \eqref{eq:ode} with \eqref{eq:finite-diff}. As such, intuitively \eqref{eq:sa} may have similar behavior as the ODE flow in \eqref{eq:ode}.}
Furthermore, $\sequence{\noise}[k+1][\nset]$ is a martingale difference sequence with respect to (\emph{w.r.t.}) the filtration $\sequence{\mcf}[k][\nset]$. 
Then, \eqref{eq:sa} is a Cauchy-Euler approximation for solving \eqref{eq:ode} with the  stepsize sequence $\sequence{\step}[k][\nset]$.

Eq.~\eqref{eq:sa} gives a prototype \emph{stochastic algorithm} for tackling many ML and SP problems. A general algorithm design procedure is to find a desired $\hg(\cdot)$ and embed the problem at hand into~\eqref{eq:rootfind}. Subsequently, we design the stochastic field $\Hg(\cdot)$ to approximate $\hg(\cdot)$ and apply the {\SA} recursion \eqref{eq:sa}. Note that the design of the stochastic field shall also respect practical constraints such as limited computation complexity, hardware limitation on arithmetic, and the availability of stochastic samples, to list a few.

\revisionupdates{The stepsize sequence $(\gamma_k)_{k\geq1}$ plays a critical role on the behavior of the {\SA} scheme, and is thus one of the most significant hyper-parameters. One motivation of the theoretical analysis is to propose theoretically grounded rules to tune this sequence, and possibly to ensure that a given choice results in a (worst case) optimal behavior. In the {\SA} case, from a high-level standpoint, 
there exists a tradeoff between the average process $\dprm_{k+1}= \dprm_k + \step_{k+1} \hg(\dprm_k)$, which necessitates sufficiently large steps to converge towards a critical point (e.g., $\sum_{k\geq 1} \gamma_k =\infty$), and the bias and randomness of the oracle $\Hg(\prm_k, \State_{k+1})$, which typically requires limiting the magnitude of the learning rate (e.g., $\sum_{k\geq 1} \gamma_k^2 <\infty$). }

\revisionupdates{One example is the stochastic optimization problem $\min_{ \prm  \in \rset^d} F(\prm)$ where $F(\prm) \eqdef  \PE\left[ \ell (\prm, \State) \right]$. We can set the mean field $\hg( \prm ) \eqdef -\grd F( \prm )$ as the gradient. Under regularity conditions, the stochastic oracle $\Hg(\prm_k, \State_{k+1} )$ satisfying $\CPE{\Hg(\prm_k, \State_{k+1})}{\mcf_k} = \hg(\prm_k)$ is a \emph{stochastic gradient} defined by $\Hg(\prm_k, \State_{k+1}) \eqdef -\Deriv_{10} \ell ( \prm_k, \State_{k+1})$ at the $(k+1)$th iteration, and $\State_{k+1}$ is an i.i.d.~random sample following the same law of $\State$ in $\PE\left[F(\prm, \State) \right]$. In this case, the {\SA} recursion \eqref{eq:sa} yields the popular stochastic gradient ({\SG}) algorithm \cite{lan2020first}.
However, the {\SA} recursion \eqref{eq:sa} is not limited to the {\SG} algorithms: it can cover more general scenarios, e.g., the function $\hg(\prm_k)$ is not necessarily the gradient for any objective function in an optimization problem, or the expectation $\CPE{ \Hg(\prm_k, \State_{k+1}) }{\mcf_k}$ is not a gradient map. }

\revisionupdates{For readability,  \Cref{tab:all_not_genSA} aggregates all notations that are used for the generic analysis of the \SA~scheme, specifically in \Cref{sec:SAtuto,sec:assume,subsec:finite-time-bounds,sec:asymptotic-convergence,sec:vr}.}

\renewcommand{\arraystretch}{1}
\begin{table}
    \caption{Summary of notations used in the general analysis of \SA, in \Cref{sec:SAtuto,sec:assume,subsec:finite-time-bounds}}
    \label{tab:all_not_genSA}
\resizebox{\linewidth}{!}{\begin{tabular}{p{0.33\linewidth}p{0.7\linewidth}p{0.21\linewidth}}\toprule
Notation & Object & Def.~in   \\
\midrule
\multicolumn{3}{r}{In \Cref{sec:SAtuto}}\\
\midrule
$h$ & Mean  field &  \cref{eq:rootfind} \\
$d$ & Problem dimensionality &  \cref{eq:rootfind} \\
$(\prm_{k})  $& Seq. of iterates &  \cref{eq:sa} \\
$(\step_{k})$ & Seq. of step-size &   \cref{eq:sa} \\
$(\Hg(\prm_k, \State_{k+1} ))$ & Stochastic oracle on the field &   \cref{eq:sa} \\
$\noise_{k+1} $& Noise perturbation  &  \\ 
$\State_{k+1} $ & $X$-valued random variable &   \cref{eq:sa} \\
$\mcf_k=\sigma(\prm_0, (\State_\ell)_{\ell \leq k})$ & filtration adap. to $(\prm_{k}) $ &    \cref{eq:sa} \\
\midrule
\multicolumn{3}{r}{In \Cref{sec:assume}}\\
\midrule
$\superlyap: \rset^d \to \rset_+$ &    non-negative Borel function &   \\
$\clyap_{\hg,0},  \clyap_{\hg,1} \in \rset_+$   & Constants controlling the mean field $\hg$   &    \Cref{assum:field}\\
$\boundbias_{0,k},\boundbias_{1,k} \in \rset_+$ & Seq. of constants controlling the expected oracle field      &    \Cref{assum:field}\\
$\boundvar_0,\boundvar_1 \in \rset_+$ &  Constants controlling the variance of the  oracle field   &    \Cref{assum:field}\\
 $\lyap$  &   smooth Lyapunov function   &    \Cref{assum:lyapunov}\\
 $\Liplyap$  &   smoothness of  $\lyap$   &    \Cref{assum:lyapunov}\\
 $\rholyap$  &   link between $\superlyap$ and $\lyap$   &    \Cref{assum:lyapunov}\\
$ \clyap_{\lyap}$ &  upper bound on $\nabla \lyap$ w.r.t.~$\superlyap$    &    \cref{eq:definition:clyap-lyap} \\
$\Lambda_{\lyap}$&  set of points where $\pscal{\nabla \lyap(\prm)}{\hg(\prm)}=0$    &    \cref{eq:definition-Lambda-V} \\
$\equilibrium(\hg) $&   set of equilibrium points of the vector field $\hg$  &    \cref{eq:equilibrium-set}\\
$\varepsilon>0$&   target precision  &   \\
$R$ & (random) stopping time &   \\
$T$ & total number of iterations &  \\
$R_{\totstep}$ & stopping rule strategy after $T$ steps &  \\ 
\midrule
&& \Cref{subsec:finite-time-bounds} \\
\midrule
$(\bsf_\ell, \eta_\ell)_{\ell\in \{0,1\}} , \step_{\max}, \omega_k $ & Constants function of $ (\boundvar_\ell, \boundbias_\ell,  \clyap_{\hg,\ell})_{\ell\in \{0,1\}} $, $\Liplyap, \clyap_{\lyap}$ &  \cref{eq:definiton-bsf-0,eq:definition-eta,eq:step-max,eq:definition-u-k} \\
$\bar{\prm}_T $ & weighted average of the parameters & \Cref{rem:convex-superlyap}\\
$\Vinit$ & Initial condition term & eq. \eqref{eq:definition-biais} \\
$\Bterm$ & Non vanishing bias  & eq. \eqref{eq:definition-biais}\\
$\lambda_{k}, b_{k} $ & Time dependent constants in \Cref{theo:main-result-quantitative-fast} & eqs. \eqref{eq:definition-lambda-k}, \eqref{eq:definition-b-k}\\
$\Lambda_{j:k} $ & Shortcut for $  \prod_{i=j}^k \lambda_i$ & \\
$\beta \in \ocint{0,1}$ & Step-size decrease rate & \eqref{eq:diminishing-step} \\
\bottomrule
\end{tabular}
}
\end{table}
\subsection{Motivating Examples}\label{sec:ex}

\subsubsection{Stochastic Gradient Descent}
\label{subsec:SGDintro}
We begin our exposition on {\SA} schemes with stochastic gradient ({\SG}) algorithms, the simplest yet most popular setting in SP and ML. 
Stochastic gradient algorithms are the workhorse of modern machine learning and data-driven optimization  \cite{bottou2003stochastic,bottou2007tradeoffs,bottou2018optimization}. Much of the success is due to their broad applicability -- the stochastic gradient algorithm generally works for any problem for which there is an unbiased gradient estimator.
For convex problems, a long line of work \cite{agarwal2009information,nemirovski1994parallel,nemirovski2009robust,juditsky2011first,juditsky2011first-ii,hazan2014beyond,cevher2014convex,qiu2016survey} sheds lights on convergence properties of {\SG}, and they are by now well-understood; while non-convex problems are discussed in \cite{bottou2018optimization, ghadimi2013stochastic}. 


To describe the general setting, the mean field $\hg( \prm )$ of stochastic gradient algorithm as an {\SA} scheme corresponds to the gradient of an objective function that we aim at minimizing.  We consider a differentiable objective function $F: \rset^d \to \rset$. A necessary condition on a model $\prm$ to be a minimizer of $F$ is to satisfy $\nabla F(\prm)= 0$. We consider $\hg(\prm) \eqdef -\nabla F(\prm)$ and look for points $\prm$ such that $\hg(\prm)=0$.

We highlight two fundamental situations in the case of discriminative learning.
\begin{enumerate}[label=\alph*), leftmargin=*]
    \item \emph{Expected Risk Minimization (ERM) for streaming data.} The function $F$ is the expected loss $\ell$ on an observation $\State$ (with distribution $\rho$) of a model $\prm$: $F(\prm) \eqdef \PE[\loss( \prm, \State )]$.
 At  iteration $k+1$, a new observation $\State_{k+1} \sim \rho$, independent from the past, is revealed. The random field is 
 \begin{equation}
 \Hg( \prm_k, \State_{k+1} ) \eqdef - \Deriv_{10} \ell ( \prm_k, \State_{k+1}). \label{eq:randomfield_Streaming_SGD}
 \end{equation}
    \item \emph{ERM for Batch Data}. The function $F$ is the empirical loss $\ell$ over a set of observations $(Z_1, \dots, Z_n)$  for a model $\prm$: $F(\prm) \eqdef n^{-1} \sum_{i=1}^n \loss( \prm, Z_i)$.
 At  iteration $k+1$, a random index $\State_{k+1} \in \{1,\dots, n\}$, independent from the past and with uniform distribution on $\{1, \cdots, n\}$, 
 is sampled by the learner. The random field is defined as \begin{equation}
 \Hg( \prm_k, \State_{k+1}) \eqdef  -\Deriv_{10} \ell ( \prm_k, Z_{\State_{k+1}}). \label{eq:randomfield_ERM_SGD}
 \end{equation}
\end{enumerate}

\begin{remark}
Extension to mini-batch {\SG}~\cite{dekel2012optimal}: at iteration $k+1$, the learner may receive (resp. sample) a number $\lbatchgd$ (called mini-batch size) of observations (resp. indices).  For the streaming case, the random field $\Hg$ is then  $\Hg( \prm_{k}, \State_{k+1} ) \eqdef  \lbatchgd^{-1}\sum_{i=1}^{\lbatchgd} \Deriv_{10} \ell (\prm_k, Z_{\lbatchgd k+i} )$ where $\sequence{Z}[k][\nset]$ is a sequence of \iid\  observations and $\State_{k+1} \eqdef (Z_{\lbatchgd k+1}, \ldots, Z_{\lbatchgd (k+1)})$.  For the batch case, it is  $\Hg( \prm_{k}, \State_{k+1} ) \eqdef \lbatchgd^{-1}   \sum_{i\in \State_{k+1}} \Deriv_{10} \ell ( \prm_k, Z_{i} ) $, with $\State_{k+1}$ is a subset  of size $\lbatchgd$  sampled at random with or without replacement in $\{1,\dots, n\}$. \revisionupdates{The choice of $\lbatchgd$ leads to tradeoffs between per-iteration computation cost and overall convergence rate, interested readers are referred to \cite[Sec.~4]{bottou2018optimization} for details.}
\end{remark}
In the case of batch data, it is possible to use a non-uniform distribution to sample indices - see, for example, \cite{needell2014stochastic}. Moreover, the objective function has a finite-sum structure, and the problem is often rewritten as follows for simplicity:
\begin{equation} \label{eq:finitesum}
\textstyle \min_{ \prm \in \rset^d } F( \prm ), \qquad F(\prm) \eqdef \frac{1}{n} \sum_{i=1}^n f_i( \prm ) \eqsp,
\end{equation}
where for all $i \in \{1, \dots, n\}$,  $f_i: \rset^d \rightarrow \rset$, and $ f_i( \prm )\eqdef \ell (\prm, Z_i)$.
The random field in~\eqref{eq:randomfield_ERM_SGD} can then be written as
\begin{equation}\label{eq:SGD-finite-sum}
\textstyle \Hg(\prm,\State)= - \sum_{i=1}^n X^i \, \nabla f_i(\prm),
\end{equation}
where $\State\eqdef (X^1,\dots, X^n) \in \{0,1\}^n$ is an $n$-dimensional binary sampling vector with $\sum_{i=1}^n X^i = 1$. If the law of the sampling vector $\State$ is such that $\PE[X^i]= 1/n$ for any $i$, then $\PE[\Hg(\prm,\State)]= -\nabla F(\prm)$.
We observe that {\SA} scheme using the above $\Hg(\cdot)$ yields the classical stochastic gradient algorithm for the finite-sum problem~\eqref{eq:finitesum}:
\begin{equation} \label{eq:definition-SGD}
\textstyle \prm_{k+1} = \prm_k -  {\gamma_{k+1}} \sum_{i=1}^n X_{k+1}^i \nabla f_i(\prm_k) \eqsp.
\end{equation}
We note that the computational complexity of \eqref{eq:definition-SGD} is independent of $n$, since $\sum_{i=1}^n X_{k+1}^i = 1$ and the learner must compute $\nabla f_i(\prm_k)$ for the only index $i$ satisfying $X_{k+1}^i = 1$.

\subsubsection{Compressed Stochastic Approximation}
\label{subsec:ex-variantsSGD}
We study compressed {\SA} methods where a compression operator is used in the update scheme.
The goal is either to reduce transmission, storage, or computational costs.
 We focus here on instantiating methods for the {\SG} case.
As a first example, we consider the Gauss-Southwell coordinate descent estimator. For high-dimensional problems ($d \gg 1$), coordinate descent methods reduce computational complexity by restricting the update to a subset of the coordinates. We consider the same optimization 
\revisionupdates{problem as}
\eqref{eq:finitesum}. Let $j_{k+1} \in \{1, \ldots, d\}$ is the chosen coordinate in the $k$-th iteration, we have
\begin{equation}\label{eq:gauss_south}
\prm_{k+1}=\prm_k-
\gamma_{k+1} \nabla_{j_{k+1}} F (\prm_k ) \, {\bm e}_{j_{k+1}} \eqsp,
\end{equation}
where $\{\bm e_1, \ldots, \bm e_d\}$ is the canonical basis of $\rset^d$ and $\nabla_j F$ is the $j$-th coordinate of the gradient.
The \emph{Gauss-Southwell selection} rule \cite{nutini2015coordinate} uses:
\begin{equation}\label{eq:gauss_south_choice_i}
j_{k+1} \eqdef \underset{j\in \{1, \dots, d\}}{\operatorname{Argmax}}\left|\nabla_j F\left(\prm_k\right)\right| \eqsp.
\end{equation}
This corresponds to a greedy selection procedure, since at each iteration we select the coordinate with the largest directional derivative. The Gauss-Southwell rule, on the other hand, corresponds to a deterministic algorithm: 
$\Hg(\prm_{k}, \sim) \eqdef \nabla_{j_{k+1}} F\left(\prm_k\right)  {\bm e}_{j_{k+1}} $ with $j_{k+1}$ in \eqref{eq:gauss_south_choice_i}, a straightforward extension to {\SG} is possible by replacing the $\nabla_j F( \prm_k )$ in \eqref{eq:gauss_south_choice_i} with $\sum_{i=1}^n X_{k+1}^i \nabla_{j} f_i(\prm_k)$ and construct the stochastic oracle as $\Hg(\prm_{k}, \State_{k+1} ) \eqdef  \sum_{i=1}^n X_{k+1}^i \nabla_{j_{k+1}} f_i \left( \prm_k \right)  {\bm e}_{j_{k+1}} $.

The coordinate descent algorithm \eqref{eq:gauss_south} is a special case of the general compressed {\SG} methods that aim to reduce the storage and/or transmission cost of {\SG}. To formally discuss the general methods, we introduce the concept of \emph{compression operators}. A compression operator on $\rset^d$ is a mapping $\compressor: \rset^d \times \sf U \to \rset^d$, where $\sf U$ is a general state space equipped with a sigma field and distribution $\mu_U$.
The operator $\compressor$ is called (random) compression if, for any $\bx \in \rset^d$, the cost of storing/transmitting $\compressor (\bx, \bU)$, with $\bU\sim \mu_U$, is almost-surely (or on average) less than the cost of storing/transmitting $\bx$ itself.

\begin{remark}\label{rem:expl-compr}
Two prevailing strategies have been used and combined to create such compression schemes~\cite{seide20141,alistarh_qsgd_2017}:
\begin{enumerate}[leftmargin=*]
    \item (random) projection: $\bx \in \rset^d$ is projected onto a smaller dimensional subspace. E.g., the \revisionupdates{$\rm{Rand}_h$} operator projects $\bx$ onto a random space generated by \revisionupdates{${\rm h}$} canonical vectors; and \revisionupdates{$\rm{Top}_h$} operator projects $\bx$ onto a deterministic space generated by the \revisionupdates{${\rm h}$} canonical vectors corresponding to the largest values of $\bx$. For example, \eqref{eq:gauss_south_choice_i} corresponds to using $\rm{Top}_1$ compressor;
    \item (random) quantization, that (randomly) maps each coordinate of $\bx$ onto a scalar codebook   and transmits the index of the corresponding codeword \cite{gersho2012vector}. E.g., assuming a uniform quantizer converts a floating point value $x \in \rset$ to the closest quantized value as:
\begin{equation}
\label{eq:deterministic-rounding}
Q_d(x, \sim) \eqdef \operatorname{sign}(x) \, \Delta \, \left\lfloor\frac{|x|}{\Delta}+\frac{1}{2}\right\rfloor,
\end{equation}
where $\Delta$ denotes the quantization resolution.
On the other hand, the quantization function for stochastic rounding is defined as:
\begin{equation}
\label{eq:stochastic-rounding}
Q_s(x,U) \eqdef \Delta \times \begin{cases}\left\lfloor\frac{x}{\Delta}\right\rfloor+1 & \text { if } U \leq \frac{x}{\Delta}-\left\lfloor\frac{x}{\Delta}\right\rfloor, \\ \left\lfloor\frac{x}{\Delta}\right\rfloor & \text { otherwise,}
\end{cases}
\end{equation}
where $U$ is uniformly distributed on $\coint{0,1}$. 
These scalar quantizers can be extended to the vector case by applying the quantization operation on each coordinate.
\end{enumerate}
\end{remark}
Let $\Hg( \prm, \State )$ denotes the {\SA} random field. We introduce three classes of compressed {\SA} methods as follows. The first class corresponds to compressing the random field:
\begin{equation} \label{eq:cSGD}
\prm_{k+1} = \prm_k + \step_{k+1} \, \compressor (\Hg(\prm_k, \State_{k+1} ) , \bU_{k+1}  ) \eqsp.
\end{equation}
It includes \eqref{eq:gauss_south} as a special case.
In the literature, these methods have attracted attention with the increasing interest in distributed optimization, \eg inexact gradient descent methods~\cite{de2020worst,gannot2022frequency}, low-precision coordinate descent methods~\cite{de2015taming}, compressed {\SG} methods~\cite{magnusson2017convergence, seide20141,karimireddy2019error,stich2020error,alistarh_qsgd_2017,bernstein2018signsgd, reisizadeh2019exact}, quantized algorithms for wireless sensor networks \cite{rabbat2005quantized, ribeiro2006bandwidth, msechu2011sensor, yi2014quantized} and the references therein.

The second class of compressed {\SA} refers to recursion where $\Hg$ is observed at a point (slightly) different from $\prm_k$ called the \textit{perturbed iterate}:
\begin{equation}
\label{eq:binary-connect}
\prm_{k+1} = \prm_k + \step_{k+1} \, \Hg( \compressor(\prm_k, \bU_{k+1}) ,\State_{k+1} ) \,.
\end{equation}
In the {\SG} for deep learning, the above setup is known as a Straight-Through Estimator (STE) introduced by \cite{Courbariaux2015}, which quantizes the model $\prm$ before computing the gradient oracle; then the {\SG} update is performed using a full precision buffer. In the convex optimization literature, recursion has been studied as perturbed iterations by \cite{mania2017perturbed}, and also in the randomized- smoothing approach, where the observation point is intentionally perturbed (e.g., by a Gaussian noise) to achieve better regularity~\cite{duchi2012randomized,scaman2018optimal}.
The same approach includes the study of {\SG} with asynchrony (i.e., the field can be measured on an 'old' iterated model)~\cite{agarwal2011distributed,chaturapruek2015asynchronous} or in distributed systems where the gradient is observed on a local model held only by the local workers~\cite{stich2020error,philippenko2021preserved}.

The third class of compressed {\SA} method is the recursion: 
\begin{equation}
\label{eq:de-sa}
\prm_{k+1} = \compressor( \prm_k + \step_{k+1} \, \Hg( \prm_k, ,\State_{k+1} ) , \bU_{k+1} ) \,.
\end{equation}
Note that \eqref{eq:de-sa} is a special case of \eqref{eq:sa} with the random field $\widetilde{\Hg}( \prm_k, \bU_{k+1}, \State_{k+1} )$ given by
\begin{equation} \label{eq:de-sa-H}
\frac{1}{\gamma_{k+1}} \left( \compressor( \prm_k + \step_{k+1} \, \Hg( \prm_k, ,\State_{k+1} ) , \bU_{k+1} ) - \prm_k \right) .
\end{equation}
In the {\SG} case, the above setup is known as the low-precision {\SG}, introduced in \cite{courbariaux2014training, gupta2015deep, de2018high, li2017training}, which quantizes the model $\prm$ after computing the gradient oracle. 

The use of low-precision arithmetic plays an essential role in SP and ML, where frugal algorithms are often mandatory. To train models with a low-precision representation of the parameters, we can apply the quantization function $\compressor(\cdot, \bU)$ to convert  entries of the parameter vector $\prm$ into a quantized/rounded version $\hat{\prm}=\compressor(\prm, \bU)$.
The first application of the STE rule was in the BinaryConnect algorithm of \cite{Courbariaux2015} for training neural networks with Boolean weights; in this case, the weights are binary $\{-\Delta, \Delta\}$. STE has been applied to many different settings and improved; see \cite{le2021adaste,liu2022nonuniform,tjandra2019end}.

\subsubsection{Stochastic EM algorithms}
\label{sec:EM:introduction}
The Expectation-Maximisation algorithm (EM) proposed by the popular work \cite{dempster:1977} was used to solve the optimization problem $\argmin_{\prmo\in \rset^d} F(\prmo)$ when $F$ is defined by an (possibly intractable) integral $F(\prmo) \eqdef - \log \int_{\widetilde \Zset} p(z; \prmo) \tilde \mu(\rmd z)$, where $p(z; \prmo)$ is a positive function and  $\tilde \mu$ is a sigma-finite measure on a measurable set $\widetilde \Zset$; see \cite{moon1996expectation,mclachlan2007algorithm} and references therein. There are numerous applications of EM, including inference of mixture distributions \cite{redner1984mixture}, robust inference in the presence of heavy tailed noise \cite{swami2000non}, Hidden Markov Models \cite{rabiner1986introduction,cappe2005springer}, factor analysis \cite{ghahramani1995factorial}, graphical models with missing data. In the subsequent discussions, we consider the case where both the function $p(z; \prmo)= \prod_{i=1}^{n} p_i(z_i;\prmo)$ and the dominating measure $\tilde \mu= \mu^{\otimes n}$  have a product form that yields 
\begin{equation}
\label{eq:EM:objective:prmo-2} 
\textstyle F(\prmo) \eqdef - \frac{1}{n} \sum_{i=1}^n \log g_i(\prmo),
\end{equation}
where
\begin{equation}\label{eq:EM:likelihoodY}
\textstyle g_i(\prmo) \eqdef   \int_\Zset p_i(z_i; \prmo) \, \mu(\rmd z_i),
\end{equation}
and $p_i(z_i;\prmo)$ is positive for all $i \in \nset$ and $(z_i, \prmo) \in \Zset \times \rset^d$.
Such an optimization problem is motivated by minimizing the Kullback-Leibler divergence computed along the examples indexed by $i \in \nset$. As expressed by \eqref{eq:EM:objective:prmo-2}-\eqref{eq:EM:likelihoodY}, in EM the divergence/loss  $F$ is not explicit and is given by an integral over a \emph{latent variable} $z=(z_1, \cdots, z_n)$.

A popular application of EM 
is the computation of Maximum Likelihood estimator. In this case,
$g_i(\prmo)$ is the log-likelihood function of an observation $Y_i$ in a latent variable model (see, e.g., \cite{dempster:1977,everitt:1984,fessler1994space,moon1996expectation,mclachlan2007algorithm,gupta2011theory}). The positive quantity $p_i(z_i; \prmo)$
is the joint probability of the observation $Y_i$ and the latent variable $z_i$ for a given value of the parameter $\prmo$; it is a shorthand notation for $p_{Y_i}(z_i; \prmo)$.
Subsequently, $F(\prmo)$ corresponds to the case where the observations are independent and  the statistical analysis is based on a given set of $n$ examples that are not necessarily identically distributed.

EM is a \emph{Majorize-Minimization} algorithm (see, e.g., \cite[chapter~8]{lange:2013}) that handles the minimization of \eqref{eq:EM:objective:prmo-2} by iterating between an \emph{Expectation step} (E-step) and a \emph{Minimization step} (M-step). Given the current iterate $\prmo_k$, the E-step defines a surrogate function $\QEM[\prmo_k] (\prmo)$ such that:
$F(\prmo_k) = \QEM[\prmo_k](\prmo_k)$ and $F(\prmo) \leq \QEM[\prmo_k](\prmo)$ for all $\prmo \in \rset^d$.
The M-step updates the parameter by selecting a minimizer 
\begin{equation} \label{eq:EM:iter:theta}
\prmo_{k+1} \eqdef \mathrm{argmin}_{\prmo \in \rset^d} \QEM[\prmo_k](\prmo) \,,
\end{equation}
which is assumed to be unique for simplicity.
Under regularity conditions, since $\QEM[\prmo_k](\cdot)$ majorizes $F(\cdot)$, $\prmo_{k+1}$ is a stationary point of the difference $\prmo \mapsto \QEM[\prmo_k](\prmo)- F(\prmo)$.
This yields 
\[
\nabla F(\prmo_{k}) = \nabla \QEM[\prmo_k](\prmo_{k})\eqsp.
\]

As the surrogate $\QEM[\prmo_k](\cdot)$ is an upper bound on the objective, an improvement on the surrogate translates to an improvement on the objective.
This descent property lends EM algorithm with remarkable numerical stability. Since each iteration is defined as a minimization of a function, EM algorithm is invariant under changes of parametrization, which is a significant advantage over first-order (gradient) method.
The EM surrogate used by \cite{dempster:1977} is given up to an additive constant (which depends on $\prmo'$) as
\begin{align}\label{eq:EM:QEM}
\QEM[\prmo'](\prmo) & \eqdef - \frac{1}{n} \sum_{i=1}^n \int_\Zset \log p_i(z_i; \prmo) \, \pi_i(z_i; \prmo') \,  \mu(\rmd z_i),
\end{align}
where $z \mapsto \pi_i(z; \prmo)$ is the probability density function (p.d.f.) on $\Zset$ defined by
\begin{equation} \label{eq:EM:posterior}
\pi_i(z; \prmo) \eqdef p_i(z; \prmo)/g_i(\prmo).
\end{equation}
In many applications, $\pi_i(z_i; \prmo)$ is the posterior distribution of the latent variable $z_i$ given the observation $\# i$ when the value of the parameter is $\prmo$.
A useful decomposition of this surrogate, which is the key discovery of  \cite{dempster:1977}, is given by
\begin{equation}
\label{eq:surrogate-decomposition}
\QEM[\prmo'](\prmo) = -\frac{1}{n} \sum_{i=1}^n \log g_i(\prmo) + \crossentropy[\prmo'](\prmo)
\end{equation}
where $\crossentropy[\prmo'](\prmo)$ is the cross-entropy of the distribution $\pi(\chunk{z}{1}{n};\prmo) \eqdef \prod_{i=1}^n \pi_i(z_i;\prmo)$ relative to the distribution $\pi(\cdot; \prmo')$:
\begin{align*}
\crossentropy[\prmo'](\prmo)
&
\eqdef - \frac{1}{n} \sum_{i=1}^n \int_\Zset   \log( \pi_i(z_i;\prmo) ) \, \pi_i(z_i;\prmo') \, \mu(\rmd z_i) \eqsp.
\end{align*}
The cross-entropy $\prmo \mapsto \crossentropy[\prmo'](\prmo)$ is minimized  at $\prmo'$ such that $\crossentropy[\prmo'](\prmo')$ is the entropy of $\pi(\cdot; \prmo')$. Under appropriate regularity conditions, this in particular implies that, for all $\prmo \in \rset^d$,
\begin{equation}
\label{eq:characterization-minimum-CE}
\nabla \crossentropy[\prmo](\prmo)= 0 \eqsp.
\end{equation}
Moreover, the surrogate decomposition  \eqref{eq:surrogate-decomposition}  and the inequality $\crossentropy[\prmo_k](\prmo_{k+1}) \geq \crossentropy[\prmo_k](\prmo_{k})$ imply
\begin{multline*}
\textstyle -\frac{1}{n} \sum_{i=1}^n \log g_i(\prmo_{k+1}) + \frac{1}{n}  \sum_{i=1}^n \log g_i(\prmo_{k}) \\
\leq  \QEM[\prmo_k](\prmo_{k+1}) - \QEM[\prmo_k](\prmo_k)  \leq 0 \eqsp,
\end{multline*}
where we used the definition of $\prmo_{k+1}$ in the \rhs. Hence, any update of the EM algorithm leads to an increase of the function $\prmo \mapsto \frac{1}{n} \sum_{i=1}^n \log g_i(\prmo)$.

The limiting point of the EM algorithms are the fixed point of the EM mapping, \ie the parameters $\prmo_\star$ which satisfy
\[ \textstyle
\prmo_\star = \argmin_{\prmo \in \rset^d} \QEM[\prmo_\star](\prmo).
\]
Under appropriate regularity conditions (see \eg \cite{wu1983convergence}), we have from \eqref{eq:EM:iter:theta} and \eqref{eq:surrogate-decomposition} that
\begin{equation}
\label{eq:differential-EM}
\textstyle \nabla \crossentropy[\prmo_k](\prmo_{k+1}) = \frac{1}{n} \sum_{i=1}^n \nabla  \log g_i(\prmo_{k+1}).
\end{equation}
If $\prmo_\star$ is a fixed point of \eqref{eq:EM:iter:theta}, together with \eqref{eq:characterization-minimum-CE}, we have
\begin{equation}
\label{eq:differential-EM}
\textstyle 0= \nabla \crossentropy[\prmo_\star](\prmo_\star) = \frac{1}{n} \sum_{i=1}^n \nabla  \log g_i(\prmo_\star) \eqsp,
\end{equation}
showing that the fixed points of the EM coincide with the roots of the gradient of the objective function (see \eqref{eq:EM:objective:prmo-2}).

We restrict our attention to the case where
$p_i$ belongs to the curved exponential family (see, e.g., \cite{brown:1986}).
Exponential family models are important special cases as the  E-step amounts only to computing a conditional expectation. In particular,
\begin{assumEM}
  \label{assumEM:expo} There exist measurable functions
  $\phiem: \rset^d \to \rset^d$,   $\psiem: \rset^d \to \rset$,  and for all $ i \in \{1, \cdots, n\}$, $\sem_i: \Zset \to \rset^d$
such that
\[
\log p_i(z_i;\prmo) \eqdef  \pscal{\sem_i(z_i)}{\phiem(\prmo)} - \psiem(\prmo).
\]
\end{assumEM}
In the terminology used for exponential families, the functions $\sem_i$'s are called the \emph{sufficient statistics}.
In the applications, $\sem_i$ depends on $i$ through the observation $Y_i$. Moreover,
\begin{assumEM}
  \label{assumEM:Mstep} There exists a measurable function $\mapem: \rset^d \to \rset^d$ such that for any $s \in \rset^d$
\[
\mapem(s) \eqdef \mathrm{argmin}_{\prmo \in \rset^d}  \{ \psiem(\prmo) - \pscal{s}{\phiem(\prmo)} \}.
\]
\end{assumEM}
In most cases, the optimization problem defined for $\mapem$ is strongly convex, and in some cases, it can be solved in closed form. 
Under \Cref{assumEM:expo}, the $\QEM[\prmo']$ function writes
\[
\QEM[\prmo'](\prmo) \eqdef  \psiem(\prmo)  - \pscal{\barsem(\prmo')}{\phiem(\prmo)},
\]
where  $\barsem(\prmo')$ is the mean value of the expectations of the $\sem_i$ functions under the p.d.f. $\pi_i$'s:
\begin{align}\label{eq:EM:def:bars}
&\barsem(\prmo') \eqdef \frac{1}{n} \sum_{i=1}^n \barsem_i(\prmo'), \\
&\text{where} \quad  \barsem_i(\prmo'):=  \int_\Zset \sem_i(z_i)  \, \pi_i(z_i; \prmo') \mu(\rmd z_i).
\end{align}
At each iteration of the EM algorithm, the computation of the surrogate function  $\prmo \to \QEM[\prmo_k](\prmo)$
boils down to the computation of the expectation $\barsem(\prmo_k)$.
\Cref{assumEM:expo,assumEM:Mstep} imply that a step of the EM algorithm is expressed as
\[
\prmo_{k+1} = \mapem \circ \barsem (\prmo_k),
\]
thus showing that the fixed points of the EM mapping are the roots of
$\prmo \mapsto \mapem \circ \barsem(\prmo) -\prmo$.

Note that $\barsem(\prmo_k)$ might be seen as a double expectation: the \emph{inner} expectation amounts to evaluating $\barsem_i(\prmo')$,
and the \emph{outer} integral is an average over the $n$ functions $\barsem_i$. In large-scale learning, the outer integral is intractable or has prohibitive computational cost; it may also be the case that the inner integrals are not explicit, e.g., when $\pi_i$ is known except for a normalizing constant, and its expression or the geometry of $\Zset$ is complicated. The \emph{Stochastic EM} algorithms were developed to avoid a
full scan of the $n$ functions at each iteration and to allow a stochastic approximation of the inner integrals by Monte Carlo sampling. The stochastic EM algorithms described in \cite{celeux:diebolt:1985,Wei:tanner:1990,delyon1999convergence,fort:moulines:2003}
address the intractability of inner expectation; the algorithms in \cite{neal:hinton:1998,Ng:mclachlan:2003,cappe2009line,chen2018stochastic,karimi:etal:2019,fort2020stochastic,fort2021fast}
address the intractability of the  outer expectation; \cite{fort:moulines:2022} addresses both intractabilities. 

Many if not all stochastic EM  algorithms are instances of {\SA}. The key ingredient is the following result (see, for example, \cite[section 7]{delyon1999convergence}):
if $\prmo_\infty$ is a root of $\prmo \mapsto \mapem \circ \barsem(\prmo) -\prmo$ on $\rset^d$, then
$\prm_\infty \eqdef \barsem(\prmo_\infty)$ is a root of
$\prm \mapsto \barsem \circ \mapem(\prm) - \prm$ on
$\rset^d$; \textit{(ii)} if $\prm_\infty$ is a root of
$\prm \mapsto \barsem \circ \mapem(\prm) - \prm$ on $\rset^d$, then
$\prmo_\infty \eqdef \mapem(\prm_\infty)$ is a root of
$\prmo \mapsto \mapem \circ \barsem(\prmo) -\prmo$ on $\rset^d$. Consequently, EM can be executed
 in the $s$-space by running a {\SA} algorithm to solve the root-finding problem
\begin{equation}\label{eq:EM:statisticspace}
\prm \in \rset^d, \qquad \barsem \circ \mapem(\prm) - \prm = 0.
\end{equation}
Stochastic EM algorithms define a sequence $\sequence{\prm}[k][\nset]$ by the
iteration $\prm_{k+1} = \prm_k + \step_{k+1} \Hg(\prm_k, \State_{k+1})$
where $\Hg(\prm_k, \State_{k+1})$ is a random oracle of  $\hg(\prm_k) \eqdef \barsem \circ \mapem(\prm_k)
- \prm_k$, the mean field $\hg$
evaluated at the current iterate $\prm_k$;  and $\sequence{\step}[k][\nset]$ is a deterministic stepsize
sequence. The close links between these versions of EM in the $s$-space and mirror descent algorithms are highlighted in the recent work \cite{kunstner2021homeomorphic}.

To develop a stochastic EM algorithm from \eqref{eq:sa}, we note the randomness $\State_{k+1}$
defines a stochastic approximation of the inner and/or outer
expectations in the map $\barsem$, see \eqref{eq:EM:def:bars}. Given a {\SA} sequence $\sequence{\prm}[k][\nset]$ converging to a solution $\prm_\star$ of the fixed point problem \eqref{eq:EM:statisticspace}, the sequence  $\prmo_k \eqdef \mapem(\prm_k)$ converges to $\prmo_\star$ which is, thanks to \eqref{eq:differential-EM}, a stationary point of the function $F$, \ie\ $\nabla F(\prmo_\star)=0$. Of course, formulating precisely this technical result requires assumptions on the regularity of the model.
Among the many stochastic versions of EM in the literature  (see references above), let us make the stochastic mean field $\Hg$ explicit for two of them.

\paragraph{Mini-batch EM} \label{sec:EM:intro-MinibatchEM} This algorithm is an adaptation to the finite-sum context of the \emph{Online EM} \cite{cappe2009line}, which is designed to process a data stream. Mini-batch EM avoids computing the exact outer expectation in \eqref{eq:EM:def:bars} at each iteration of EM; it replaces the sum over $n$ terms with a sum over a mini-batch chosen at random. The stochastic oracle is given by
\begin{equation}\label{eq:onlineEM:oracle}
\Hg(\prm_k, \State_{k+1}) \eqdef \frac{1}{\lbatchem} \sum\nolimits_{i \in \State_{k+1}} \barsem_i\left( \mapem(\prm_k) \right) - \prm_k,
\end{equation}
where $\State_{k+1}$  is a set  of size
$\lbatchem$, collecting indices picked at random with or without replacement  in $\{1, \cdots, n\}$. See e.g. \cite{maire:etal:2012,nguyen:hal-02415068,oudoumanessah:hal-03824951} for an application.

\paragraph{Stochastic Approximation EM (SAEM)} \label{sec:EM:intro-SAEM} This algorithm, proposed by \cite{delyon1999convergence}, deals with the case the inner expectations in \eqref{eq:EM:def:bars} are intractable and must be approximated by Monte Carlo, while the outer expectations in \eqref{eq:EM:def:bars} can be computed with reasonable computational effort; see \eg \cite{tourneret2003bayesian,sahmoudi:etal:2005,allassonniere:etal:2006,septier:etal:2008, richard:etal:2009,yildirim2009hybrid,lindstein:2013,zhang:singh:etal:2013,svensson:etal:2014,boisbunon:zerubia:2014,braham:etal:2016,liu:etal:2019,liu:etal:2019:eusipco,zhou2020student,zhou2021parameter,yildirim:etal:2022,sauty:durrleman:2022} for applications. The stochastic oracle is given by
\begin{equation}
\label{eq:oracle-SAEM}
\Hg(\prm_k, \State_{k+1}) \eqdef \frac{1}{n} \sum_{i=1}^n  \sum_{j=1}^{\nbrMC} \rho_{i}^j(\prm_k)  \,  \sem_i(Z_{i,k+1}^{j}) - \prm_k
\end{equation}
where $\State_{k+1} \eqdef (Z_{i,k+1}^{j}, 1 \leq i \leq n, 1 \leq
j \leq \nbrMC)$ collect \revisionupdates{the $\nbrMC$ outputs} of $n$ distinct  Monte Carlo
samplers designed to approximate the distributions
$\pi_i(z_i; \mapem(\prm_k))$, for $i = 1,\ldots, n$.

If i.i.d. sampling is possible from the p.d.f. $\pi_i(z_i; \mapem(\prm_k))$, then $\rho_{i}^j(\prm_k) \eqdef 1/\nbrMC$. If sampling according to the conditional distribution $\pi_i(z_i; \mapem(\prm_k)))$ is not possible, more sophisticated sampling techniques must be used. For illustration, we consider below the importance sampling procedure, a method of using independent samples from a proposal distribution $\tilde \pi_{i}(z_i; \mapem(\prm_k))$ to approximate the expectation with respect to the target distribution $\pi_i(z_i; \mapem(\prm_k))$.  The importance sampling estimator approximates  the target distribution by a random probability measure using weighted samples that are generated from the proposal; (see e.g. \cite[Chapter V]{asmussen:glynn:2007} and \cite{agapiou:etal:2017}).
More precisely, the self-normalized Importance Sampling estimator works as follows. We first sample independently $Z_{i,k+1}^1, \ldots, Z_{i,k+1}^{\nbrMC}$ from the proposal $\tilde \pi_{i}(z_i; \mapem(\prm_k))$ and define the normalized importance weights\vspace{-.1cm}
\begin{equation}
\label{eq:self-normalized}
\rho_{i}^j(\prm_k) \eqdef \frac{p_i(Z_{i,k+1}^j; \mapem(\prm_k))}{\tilde \pi_{i}(Z_{i,k+1}^j; \mapem(\prm_k))} \, \left( \sum_{\ell=1}^{\nbrMC} \frac{p_i(Z_{i,k+1}^\ell; \mapem(\prm_k))}{\tilde \pi_{i}(Z_{i,k+1}^\ell; \mapem(\prm_k))}\right)^{-1}.
\end{equation}
Instead of IS, another option is to use Markov Chain Monte Carlo (MCMC) samplers targeting the p.d.f. $\pi_i(z_i; \mapem(\prm_k))$; see for example \cite{delyon1999convergence,tourneret2003bayesian,fort:moulines:2003,karimi:etal:2019,fort2021fast}.
\subsubsection{TD Learning}
\label{sec:definition-TD-learning}
Many tasks in reinforcement learning such as policy evaluation, Q-learning \cite{sutton2018reinforcement}, etc., can be formulated as root finding problems whose effective solutions are often given as non-stochastic-gradient {\SA} recursions. Below, we select the temporal difference (TD) learning algorithm \cite{sutton1988learning} for policy evaluation to illustrate another aspect of general principle of stochastic algorithm design.

We follow the derivations of \cite{bhandari2018finite} in this example.
Consider the problem of evaluating the value function of applying a policy $\policy$ in a Markov Decision Process.
The policy $\policy$ specifies the conditional probability of choosing an action given a certain state, It induces a Markov Reward Process (MRP) given by the tuple $( \stateMRP, \kerMRP, \rewardMRP, \lambda )$:  $\stateMRP= \{s_1,\dots,s_n\}$ is the state-space (assumed for simplicity to be finite);  $\kerMRP$ is the $n \times n$ state transition matrix of the probability of transition from a given state to another;  $\rewardMRP$ is the reward function such that $\rewardMRP(s, s')$ associates a reward with each state transition; $\lambda \in \ooint{0,1}$ is the discount factor.
With a slight abuse of notation, we define the expected instantaneous reward from state $s $, for $s \in \stateMRP$,
\[ \textstyle
\rewardMRP(s) \eqdef \sum_{s' \in \stateMRP} \kerMRP(s,s') \rewardMRP(s,s');
\]
As a standing assumption, we concentrate on ergodic MRPs:
\begin{assumTD}
\label{assum:TD:stationary-policy}
The non-negative matrix $\kerMRP$ is irreducible and has a unique stationary distribution $\statdistMRP$, \ie\ $\statdistMRP \kerMRP= \statdistMRP$.
\end{assumTD}
\noindent Note that $\statdistMRP(s) > 0$ for any $s \in \stateMRP$ under \Cref{assum:TD:stationary-policy}.

The value function $\valuefunc$ of the above MRP is the expected cumulative discounted reward for a given state $s \in \stateMRP$, \ie
\begin{align}\label{eq:valuefunc_def}
\valuefunc(s)
& \eqdef \CPEx{\sum_{k=0}^{\infty} \lambda^k \rewardMRP(S_k)}{S_0=s} = \sum_{k=0}^\infty \lambda^k \kerMRP^k \rewardMRP(s) ,
\end{align}
where the expectation is over the distribution of the Markov chain $\sequence{S}[k][\nset]$, started at $S_0=s$, with Markov kernel $\kerMRP$. Note that 
$\kerMRP^k \rewardMRP(s) = \sum_{s' \in \stateMRP} \kerMRP^k ( s, s' ) \rewardMRP( s' )$.
This value function obeys the Bellman equation $\bellman  \valuefunc= \valuefunc$ where the Bellman operator $\bellman$ is defined as 
\begin{equation}
\label{eq:definition-bellman}
\textstyle \left[\bellman V\right](s) \eqdef \rewardMRP(s)+\lambda \sum_{s^{\prime} \in \stateMRP} \kerMRP(s,s') V(s'),~\forall~s \in \stateMRP \eqsp,
\end{equation}
for any function $V: \stateMRP \to \rset$.
\revisionupdates{Assume bounded reward function,}
the value function $\valuefunc(s)$ is well defined and is the only solution to \eqref{eq:definition-bellman} \cite{bertsekas2012dynamic,sutton2018reinforcement}. In most applications, the cardinality $n$ of the state space $\stateMRP$ is large. It is often advocated in such case to resort to parametric approximations, for example by using a linear function approximation or deep neural networks \cite{mnih2015human,fujimoto2018addressing}.
For simplicity, we focus on linear function approximation, where $\valuefunc(s) \approx \valuefunc[\prm](s) \eqdef \feature(s)^{\top} \prm$;
$\feature(s) \in \rset^d$ is called the \emph{feature vector} for the state $s \in \stateMRP$, and $\prm \in \rset^d$ is a parameter vector to be estimated. \revisionupdates{Without loss of generality, we assume that
\begin{assumTD}\label{assum:TD:normed:calR}
For all $s,s' \in \stateMRP$, $|\rewardMRP(s,s')| \leq 1$, $\max_{s \in \stateMRP} \| \feature(s) \| \leq 1$.
\end{assumTD}}

Since the state space is finite, the value function $\valuefunc[\prm]$ can be represented as a vector in $\rset^n$, whose $i$-th coordinate is $\valuefunc[\prm](s_i)$. This vector can be written compactly as
\begin{equation}
\label{eq:definition_feature-matrix}
\valuefunc[\prm]= \Feature \prm, \quad \Feature \eqdef [\feature(s_1),\dots,\feature(s_n)]^{\top} \,,
\end{equation}
such that $\Feature$ is the $n \times d$ feature matrix.
The linear function approximation restricts the set of admissible value function $\valuefunc[\prm]$ to $\Span(\Feature)$, the subset of $\rset^n$ spanned by  the columns of $\Feature$. As a result, the Bellman equation $\bellman  \valuefunc= \valuefunc$ may no longer be satisfied by $\valuefunc[\prm]$ for any $\prm \in \rset^d$.

Among the many TD learning algorithms (see e.g. \cite[Section B]{tsitsiklis1997analysis}), we consider in this paper the so-called TD(0)-learning algorithm. TD(0) is an {\SA} scheme: under the on-policy setting where the data is generated from the MRP induced by $\policy$, the algorithm starts with an initial  estimate $\prm_0$ and at iteration $(k+1)$, it gets a new observation $\State_{k+1} =(S_{k+1}, S_{k+1}')$,  and computes the next iterate by $\prm_{k+1} = \prm_{k} + \step_{k+1} \Hg(\prm_{k},\State_{k+1})$
where
\begin{equation}\label{eq:iteration-TD0}
\Hg(\prm,(s,s')) \eqdef \left(\rewardMRP(s,s')+\lambda \feature(s')^{\top} \prm-\feature(s)^{\top} \prm \right) \feature(s) \eqsp.
\end{equation}
We  make a simplifying assumption about $\State_{k+1}$:
\begin{assumTD}
\label{assum:TD:sampling}
The sequence $\sequence{S}[k][\nset]$ is sampled independently from the stationary distribution $\statdistMRP$ and, for each $k$,  $S'_k$ is sampled from $\kerMRP(S_k,\cdot)$.
\end{assumTD}
\noindent This assumption is  classical in reinforcement learning and is suitable for algorithms that use a replay buffer \cite{di2022analysis}.

Under mild conditions, the expected value of $S'_{k+1}$ conditionally to $\{S_{k+1} = s\}$ of $\rewardMRP(S_{k+1},S'_{k+1}) + \lambda \feature(S'_{k+1})^{\top} \prm-\feature(S_{k+1})^{\top} \prm$ evaluate{} to $[\bellman \valuefunc[\prm] ](s) - \valuefunc[\prm](s)$ which gives the \emph{temporal difference error} for the Bellman equation with the estimate $\prm$.
For any function $g: \stateMRP \times \stateMRP \to \rset$, denote
\begin{equation}
\label{eq:definition-expectation}
\textstyle \PE_{\statdistMRP}[g(S_0,S'_0)] \eqdef  \sum_{s, s' \in \stateMRP} \statdistMRP(s)\kerMRP(s,s') g(s,s') \eqsp.
\end{equation}
The mean field function of the TD(0) algorithm is  given by
\begin{align}
\hg(\prm) & \eqdef \PE_{\statdistMRP}[\feature(S_0) \rewardMRP(S_0,S'_0)]  \nonumber \\
& \qquad + \PE_{\statdistMRP}[\feature(S_0)\{\lambda \feature(S'_0) - \feature(S_0)\}^\top] \, \prm  \nonumber \\
& = \Feature^{\top} \Dtd_{\statdistMRP}(\bellman \Feature \prm - \Feature \prm) \label{eq:TD:mean-field}
\end{align}
where we have set the diagonal matrix
\begin{equation}
\label{eq:definition-diagonal}
\Dtd_{\statdistMRP} \eqdef  \diag(\statdistMRP(s_1),\dots,\statdistMRP(s_n)).
\end{equation}
If $\prm_\star$ is a root of $\hg(\prm)= 0$, then for all $\prm'  \in \rset^d$,
\begin{equation}
\label{eq:condition-root-mean-field}
\pscal[\Dtd_{\statdistMRP}]{\Feature \prm'}{\bellman \Feature \prm_\star - \Feature \prm_\star} =0
\end{equation}
showing that the Bellman error $\bellman \Feature \prm_\star - \Feature \prm_\star$ is orthogonal to the linear subspace spanned by the column of the feature matrix $\Feature$ in the scalar product  $\pscal[\Dtd_\statdistMRP]{\cdot}{\cdot}$.
In light of \eqref{eq:condition-root-mean-field},
we may now characterize the root $\prm_\star$ to $\hg( \prm ) = {\bm 0}$ using the projected Bellman equation
\begin{equation}
\label{eq:projected-bellman}
\valuefunc[\prm]= \Proj_{\statdistMRP} \bellman \valuefunc[\prm] \eqsp;
\end{equation}
 $\Proj_{\statdistMRP}$ is the projection operator onto $\Span(\Feature)$ \wrt ${\| \cdot \|}_{\Dtd_\statdistMRP}$. In \Cref{subsec:TD-learning-checking}, we will show that the equation $\hg(\prm) = 0$ has a unique root, $\prm_\star$, which is the fixed point to \eqref{eq:projected-bellman}.

\vspace{.1cm}
\noindent
\revisionupdates{\textbf{Summary.} In this section, we introduced the general \SA~scheme under consideration, and showed that it can be instantiated in several major applications of ML and SP. We summarize in \Cref{tab:all_algos} the algorithms that we introduce.}

\begin{table}[]
    \centering
    \caption{Summary of algorithms introduced in \Cref{sec:ex}\label{tab:all_algos}, that are variants of the \SA~scheme \eqref{eq:sa}.}
   \resizebox{\linewidth}{!}{ \revisionupdates{\begin{tabular}{lll}
       \toprule
       Setting & Name & Reference \\
       \midrule
       SGD &  SGD for ERM & \cref{eq:SGD-finite-sum}\\
       \midrule
       \multirow{4}{*}{Compressed SA} & Gauss-Southwell (GD with $\rm{Top}_1$)& \cref{eq:gauss_south,eq:gauss_south_choice_i} \\
& \SA~with compressed field & \cref{eq:cSGD}  \\
& STE, \SA~with quantized iterates & \cref{eq:binary-connect}  \\
& Low precision \SA & \cref{eq:de-sa-H}  \\
\midrule
   \multirow{3}{*}{EM}     &  Mini-batch EM & \cref{eq:onlineEM:oracle} \\
    &  SAEM with independent MC  & \cref{eq:oracle-SAEM} \\
    &  SAEM with importance sampling MC  & \cref{eq:oracle-SAEM,eq:self-normalized} \\
    \midrule
       TD  &  TD(0)   & \cref{eq:iteration-TD0}\\
       \bottomrule
    \end{tabular}}}
\end{table}


\section{Assumptions}
\label{sec:assumptions_and_verif}

\revisionupdates{In  this section, we first introduce in \Cref{sec:assume}  the main assumptions on the mean field $\hg$ and the random oracle $\Hg$ under which theoretical results will be derived in \Cref{sec:nonasymp,sec:asymptotic-convergence,subsec:variance-reduction}. We then establish conditions under which those assumptions are satisfied on the four previously detailed examples, in respectively \Cref{subsec:SGD-checking,sec:EM:verifH,subsec:csgd-checking,subsec:TD-learning-checking}.}

\subsection{\revisionupdates{Assumptions on the \SA~scheme.}}\label{sec:assume}
Assume that the root-finding problem \eqref{eq:rootfind} is to be solved by the {\SA} algorithm \eqref{eq:sa} which acquires the mean-field $\hg(\prm)$ via subsequent calls to a \emph{stochastic  oracle}. At iteration $(k+1)$, we denote $\prm_k$ as the current value of the model and the stochastic oracle outputs  $\Hg(\prm_k, \State_{k+1})$, where $\sequence{\State}[k][\nset]$ are random variables taking values in  ${\sf X} \subseteq \rset^\ell$. We denote by  $\mcf_k \eqdef \sigma(\prm_0, \State_\ell, 1 \leq \ell \leq k)$ the sigma-algebra generated by the random variables $\{\State_\ell\}_{\ell=1}^k$ and the initial model $\prm_0$.

We shall consider a set of assumptions for the Borel functions $\Hg: \rset^d \times {\sf X}  \to \rset^d$ and $\hg : \rset^d \to \rset^d$ in relation to the {\SA} recursion. This gives rise to the first ingredient of our analysis framework that considers a non-negative Borel function $\superlyap: \rset^d \to \rset_+$ which  controls the growth to infinity of the variance of the stochastic oracle and its bias. We  use the same function in \Cref{assum:lyapunov}, where it measures the coercivity of the mean field.
\begin{assumption}
\label{assum:field}
\begin{enumerate}[label=\emph{\alph*}),nosep,leftmargin=*]
\item \label{item:field:boundsecondmoment}For all $k \geq 0$, $\PE[\| \Hg(\prm_k, \State_{k+1})\|^2] < \infty$.
\item \label{item:clyap-hg} There exist $\clyap_{\hg,0},  \clyap_{\hg,1} \in \rset_+$ such that for all $\prm \in \rset^d$
\[
\| \hg(\prm) \|^2 \leq \clyap_{\hg,0} + \clyap_{\hg,1}\superlyap(\prm).
\]
\item \label{item:field:boundbias}
For any $k \geq 0$, there exist $\boundbias_{0,k},\boundbias_{1,k} \in \rset_+$ such that,  a.s.,
\begin{equation}
\label{eq:bound-bias}
\| \CPE{\Hg(\prm_{k},\State_{k+1})}{\mcf_{k}} - \hg(\prm_{k}) \|^2
\leq \boundbias_{0,k} + \boundbias_{1,k} \superlyap (\prm_k).
\end{equation}
\item \label{item:field:conditional-variance}
There exist  $\boundvar_0,\boundvar_1 \in \rset_+$  such that for any $k \geq 0$, a.s.,
\begin{multline}
\label{eq:bound-variance}
\CPE{\| \Hg(\prm_{k},\State_{k+1}) - \CPE{\Hg(\prm_{k},\State_{k+1})}{\mcf_{k}}\|^2}{\mcf_{k}} \\
\leq \boundvar_0 + \boundvar_1 \superlyap(\prm_{k}).
\end{multline}
\end{enumerate}
\end{assumption}
Under these assumptions, the random oracle $\Hg(\prm_k, \State_{k+1})$ may be a  biased estimator of the mean-field $\hg(\prm_k)$, with $\CPE{\Hg(\prm_{k},\State_{k+1})}{\mcf_{k}} - \hg(\prm_{k})$ being a possibly time-varying conditional bias. A special case often considered in the literature pertains to unbiased {\SA} where $\CPE{\Hg(\prm_{k},\State_{k+1})}{\mcf_{k}} = \hg(\prm_{k})$. In this case, $\boundbias_{0,k}= \boundbias_{1,k} \eqdef  0$ for all $k \in \nset$ in \Cref{assum:field}-\ref{item:field:boundbias}.
\begin{definition}[Unbiased stochastic oracle (USO)]
\label{def:USO}
The stochastic oracle in {\SA} scheme  is said to be \emph{unbiased} if $\boundbias_{0,k}= \boundbias_{1,k}= 0$ for any $k \in \nset$.
\end{definition}
\noindent In USO, the sequence $\sequence{\Noise}[k][\nset]$ where $\Noise_{k+1}\eqdef \Hg(\prm_{k},\State_{k+1}) - \hg(\prm_{k})$, is a martingale increment sequence: $\CPE{\Noise_{k+1}}{\mcf_{k}}=0$, $\PP$-a.s.

An example of USO is the case where the mean field $\hg$ has a finite sum structure: $\hg = n^{-1} \sum_{i=1}^n \hg_i$.
If $\State$ is a uniform random variable on $\{1, \ldots, n\}$, then $\Hg(\prm,\State) \eqdef \hg_\State(\prm)$ is an unbiased stochastic oracle. A generalization of this example is the case where $\hg$ is defined as an 
expectation $\hg(\prm) \eqdef \int S(z,\prm) \, \pi(\rmd z; \prm)$ with respect to a distribution $\pi$ that may depend on the current value of the parameter $\prm$; then $\Hg(\prm,\State) \eqdef N^{-1} \sum_{i=1}^N S(Z_i,\prm)$, where $\State \eqdef (Z_1, \ldots, Z_N)$ ares i.i.d. samples from $\pi(\cdot; \prm)$, induces an unbiased stochastic oracle. If instead a self-normalized importance sampling is used (see \eg\ \eqref{eq:self-normalized}), then
the oracle is  biased. In this case, the bias is inversely proportional to the number of Monte Carlo samples used in the importance sampling estimate; see \eg\ \cite{agapiou:etal:2017}.

It is worth noting that in the \emph{Robbins-Monro} setting - in reference to \cite{robbins1951stochastic} - the sequence $\sequence{\State}[k][\nset]$ is assumed to be \iid~
In comparison, \Cref{assum:field}-\ref{item:field:boundbias} and \Cref{assum:field}-\ref{item:field:conditional-variance} used here are slightly weaker since we do not assume \revisionupdates{that} $\sequence{\State}[k][\nset]$ are independent nor that they have the same distribution.  However, this excludes more subtle dependency structures and time-dependent distributions: in some situations $\sequence{\State}[k][\nset]$ is a Markov chain (possibly) controlled by $\sequence{\prm}[k][\nset]$; considering such dependency structures requires the use of sophisticated probabilistic methods, that go beyond the scope of this survey; see \cite{metivier1987theoremes,andrieu2005stability,benveniste2012adaptive,andrieu2014markovian,fort:etal:2016,karimi2019non,wang2020finite,ramaswamy2018stability} and the references therein.

\Cref{assum:field}-\ref{item:field:conditional-variance} implies the conditional variance of the random oracle $\Hg(\prm_k, \State_{k+1})$ is either bounded ($\boundvar_1=0$) or does not grow faster than  $\superlyap(\prm_k)$.
The case $\boundvar_1=0$ is called the \emph{bounded variance} case; considered in the analysis of {\SG} by \cite{ghadimi2013stochastic}.

Our analysis for {\SA} follows that of the Lyapunov function approach. In this setup, we introduce the second ingredient of our analysis framework which considers a \emph{smooth} Lyapunov function $\lyap$ for the flow of the nonlinear ODE $\rmd \dprm / \rmd t= \hg(\dprm)$.
Formally, we have the following set of assumptions:
\begin{assumption}
\label{assum:lyapunov}
There exists a function $\lyap: \rset^d \to \rset$ such that,
\begin{enumerate}[label=\emph{\alph*}),nosep]
\item \label{item:vstar}  $\lyap_\star \eqdef \inf_{\prm \in \rset^d} \lyap(\prm) > -\infty$.
\item \label{item:lyap} $\lyap$ is continuously differentiable and $\Liplyap$-smooth, \ie\ there exists $\Liplyap \in \rset_+$ such that for all $\prm,\prm' \in \rset^d$, $\| \nabla \lyap(\prm) - \nabla \lyap(\prm') \| \leq \Liplyap \|\prm - \prm'\|$.
\item \label{item:rholyap} There exists $\rholyap > 0$ such that, for all $\prm \in \rset^d$, $\pscal{\nabla V(\prm)}{\hg(\prm)} \leq -  \rholyap \superlyap(\prm)$.
\end{enumerate}
\end{assumption}
It may seem a bit complicated to have  many different constants that may seem redundant: For example, we could have fixed $\superlyap(\prm) \eqdef  - \pscal{ \nabla \lyap(\prm)}{\hg(\prm)}$ and choose $\rholyap=1$, then change the constants $\clyap_{\hg,1}$, $\boundbias_{1,k}$, $\boundvar_1$, and $\clyap_{\lyap}$ accordingly.
The motivation behind \Cref{assum:lyapunov}-\ref{item:rholyap} is to add an additional degree of flexibility that will later facilitate analysis and allow covering settings associated with different choices for $\superlyap$.

We will often use the constant $\clyap_{\lyap} > 0$ which satisfies:
\begin{equation}
\label{eq:definition:clyap-lyap}
\| \nabla \lyap( \prm) \| \leq \clyap_{\lyap} \sqrt{\superlyap(\prm)}, \quad \prm \in \rset^d \eqsp.
\end{equation}
Note that $\clyap_{\lyap}$ might be  $+\infty$ (see e.g. \Cref{subsec:SGD-checking} for some examples).
On the set
\begin{equation}
\label{eq:definition-Lambda-V}
\Lambda_{\lyap} \eqdef  \{\prm \in \rset^d, \pscal{\nabla \lyap(\prm)}{\hg(\prm)}=0\} \eqsp,
\end{equation}
the function $\superlyap$ is   identically
zero by \Cref{assum:lyapunov}-\ref{item:rholyap}:  for all $\prm \in \Lambda_{\lyap}$, $\superlyap(\prm)=0$.  However the condition $\superlyap(\prm) =0$  is not sufficient yet to be able to derive useful information for solving the root finding problem $\hg( \prm ) = {\bm 0}$. The set of equilibrium points of the vector field $\hg$ is the set
\begin{equation}
\label{eq:equilibrium-set}
\equilibrium(\hg) \eqdef \{ \prm \in \rset^d, \hg(\prm )=0\} = \{ \hg = 0 \}\eqsp.
\end{equation}
Obviously, one has $\equilibrium(\hg) \subset \Lambda_{\lyap}$, but the converse may not hold.
As we are trying to approach the equilibrium points, it is sensible to assume that $\lyap$ is a \emph{strict Lyapunov} function,
\ie
\begin{equation}
\label{eq:definition-strict-lyapunov}
\equilibrium(\hg) = \Lambda_{\lyap} \eqsp.
\end{equation}
If this is the case, the function $\superlyap$ on $\equilibrium(\hg)$ is zero.  In the case of {\SG}, the Lyapunov function $\lyap$ is the objective function to be minimized $F$. The mean field $\hg$ is the negated gradient of the objective $\hg= -\nabla F$. In this case, $\pscal{\nabla \lyap}{\hg}= -\| \nabla F \|^2$ and thus  $\lyap$ is automatically a strict Lyapunov function. For the function $\superlyap$, we typically choose $\superlyap(\prm) \eqdef  \|\nabla F(\prm) \|^2 = \| \hg(\prm)\|^2$, so that $\{\superlyap =0 \} = \equilibrium(\hg)$. See \Cref{sec:verif:assume}.

To obtain meaningful results, the function $\superlyap$ should be lower bounded outside any open neighborhood of $\equilibrium(\hg)$: for any $\delta > 0$, it is typically required that
\begin{equation}
\label{eq:condition-W}
\min_{\dist(\prm, \{\hg=0\}) \geq \delta} \superlyap(\prm) =: \epsilon_\superlyap(\delta) > 0 \eqsp,
\end{equation}
where $\dist(\prm,\mathcal{A})$ is the Euclidean distance of $\prm$ to the set $\mathcal{A}$. We have chosen not to include this condition in \Cref{assum:lyapunov} because it is not involved in the proofs of the results presented below, but only in their interpretation.
Under \eqref{eq:condition-W}, if  for any $\varepsilon>0$, there exists a stopping time $R$, possibly random, such that $\PE\left[\superlyap(\prm_{R})\right] \leq \varepsilon$, then  for any $\delta > 0$, we have
\begin{align*}
\PP( \dist(\prm_{R},\{\hg= 0\}) \geq \delta) & \leq
\PP( \superlyap(\prm_{R}) \geq \epsilon_\superlyap(\delta)
) \\
& \leq \PE[\superlyap(\prm_{R})] / \epsilon_\superlyap(\delta) \leq \varepsilon / \epsilon_\superlyap(\delta).
\end{align*}
This upper bound can be set arbitrarily small with a convenient choice of $\varepsilon$. We will prove in \Cref{sec:nonasymp} that under \Cref{assum:field} and \Cref{assum:lyapunov}, given a total number of iterations $\totstep$, there exists a  stopping rule strategy $R_{\totstep}$ such that $\PE\left[\superlyap(\prm_{R_\totstep})\right]$ goes to zero when $\totstep \to \infty$.
To understand why, we analyze the deterministic sequence $\dprm_{k+1}= \dprm_k + \step_{k+1} \hg(\dprm_k)$. Note that, by \Cref{assum:lyapunov}-\ref{item:clyap-hg},
\begin{multline*}
\lyap(\dprm_{k+1}) \leq \lyap(\dprm_k) + \step_{k+1} \pscal{\nabla \lyap(\dprm_k)}{\hg(\dprm_k)}
\\ + (\Liplyap/2) \step_{k+1}^2 \|\hg(\dprm_{k+1}) \|^2 \,.
\end{multline*}
Using that $\pscal{\nabla \lyap(\dprm)}{\hg(\dprm)} \leq -\rholyap \superlyap(\dprm)$ and assuming for simplicity $ \|\hg(\dprm_{k+1}) \|^2 \leq \clyap_{\hg,0}$, we immediately get that
\begin{multline*}
\rholyap \step_{k+1} \superlyap(\dprm_k) \leq \lyap(\dprm_k) - \lyap(\dprm_{k+1}) +  (\Liplyap/2) \step_{k+1}^2 \clyap_{\hg,0} \eqsp.
\end{multline*}
For any $\totstep > 0$, setting $\step_{k}= \frac{1}{\sqrt{\totstep}}$ for $k \in \{1,\dots,\totstep\}$, we get
\[
\frac{1}{\totstep} \sum_{k=0}^{\totstep-1} \superlyap(\dprm_k) \leq \frac{2\{\lyap(\dprm_0) - \lyap_{\star} \}  + \Liplyap \clyap_{\hg,0}}{2 \rholyap \sqrt{\totstep}} \eqsp.
\]
If $R_{\totstep}$ is a uniform random variable on \revisionupdates{$\{0,\dots,\totstep-1\}$}, then
\[
\PE[ \superlyap(\dprm_{R_{\totstep}})] \leq \frac{2\{\lyap(\dprm_0) - \lyap_{\star} \}  + \Liplyap \clyap_{\hg,0}}{2  \rholyap \sqrt{\totstep}}\eqsp.
\]
The \rhs\ goes to zero if $\totstep \to \infty$. This discussion is essentially an anticipation of the results to be obtained in \Cref{sec:nonasymp}, which introduce the bias and variance of the random oracle.

\subsection{Verifying Assumptions for the Examples}\label{sec:verif:assume}

\subsubsection{Stochastic Gradient Descent}
\label{subsec:SGD-checking}

We recall from \Cref{subsec:SGDintro} the stochastic gradient algorithm \eqref{eq:definition-SGD}. Here,
the mean field $\hg(\prm) = -\nabla F(\prm)$ is the negated gradient of the objective function. For simplicity, we consider the algorithm for the batch case~\eqref{eq:randomfield_ERM_SGD}, but results and assumption could easily be extended to include the streaming data framework~\eqref{eq:randomfield_Streaming_SGD}.
We concentrate on the case where the sampling of the index $\State_{k+1}$ in \eqref{eq:randomfield_ERM_SGD} is uniform over $\{1, \dots, n\}$. The oracle $\Hg(\prm, \State)$ of the mean field $\hg(\prm)$ is then unbiased, satisfying \Cref{assum:field}-\ref{item:field:boundbias} with $\boundbias_{0,k}= \boundbias_{1,k}=0$. Next, we describe conditions under which \Cref{assum:field}, \ref{assum:lyapunov} are valid and derive  the required constants.\vspace{.1cm}



\renewcommand{\arraystretch}{2}

\begin{table*}[t]
    \centering
\resizebox{.95\linewidth}{!}{
\begin{tabular}{ l l l l l l l l l l}
    \toprule
    Algorithm & $\lyap(\prm)$ & $\superlyap(\prm)$ & $\lyap_{\star}$ & $(\clyap_{\hg,0}, \clyap_{\hg,1})$ & $(\boundbias_0, \boundbias_1)$ & $(\boundvar_0, \boundvar_1)$ & $\Liplyap$ & $\rholyap$ & $\clyap_{\lyap}$ \\
    \midrule
    {\SG}-\Cref{assum:SGD} & $F(\prm)$ & $\| \nabla F(\prm) \|^2$ & $F_\star$ & $(0,1)$ & $(0,0)$ & $(M^2/n,0)$ & $\Lip{\nabla F}$ & $1$ & $1$ \\
    \hline
    {\SG}-\Cref{assum:SGD}-\Cref{assum:SGD:convex} & $\frac{1}{2} \| \prm - \prm_\star \|^2$ & $\pscal{\nabla F(\prm)}{\prm - \prm_\star}$ & $0$ & $(0,\Lip{\nabla F})$ & $(0,0)$ & $(M^2/n, 0)$ & $1$ & $1$ & $\infty$\\
    \hline
     {\SG}-\Cref{assum:SGD}-\Cref{assum:SGD:strongly-convex} & $\frac{1}{2} \| \prm - \prm_\star \|^2$ & $\pscal{\nabla F(\prm)}{\prm - \prm_\star}$ & $0$ & $(0,\Lip{\nabla F})$ & $(0,0)$ & $(M^2/n, 0)$ & $1$ & $1$ & $1/\sqrt{\mu}$\\
    \hline
    {\SG}-\Cref{assum:SGD}-\Cref{assum:SGD:strongly-convex} & $\frac{1}{2} \| \prm - \prm_\star \|^2$ & $\frac{1}{2} \| \prm - \prm_\star \|^2$ & $0$ & $(0, 2 \Lip{\nabla F}^2)$ & $(0,0)$ &  $(M^2/n, 0)$ & $1$ & $2 \mu$ & $\sqrt{2}$ \\
    \hline
    Mini-batch EM & $F \circ \mapem(\prm)$ & $\| h( \prm ) \|^2$ & $F_\star$ & $(0,1)$ & $(0,0)$ & $(\nicefrac{\boundvarem_0}{n \nbrMC},\nicefrac{\boundvarem_1}{n \nbrMC})$ & $\Liplyap$ & $\vminem$ & $\vmaxem$ \\
    \hline
    SAEM (\iid)& $F \circ \mapem(\prm)$ & $\| h( \prm ) \|^2$ & $F_\star$ & $(0,1)$ & $(0,0)$ & $(\nicefrac{\boundvarem_0}{\lbatchem},\nicefrac{\boundvarem_1}{\lbatchem})$ & ${\Liplyap}$ & $\vminem$ & $\vmaxem$ \\
     \hline
    SAEM (IS)& $F \circ \mapem(\prm)$ & $\| h( \prm ) \|^2$ & $F_\star$ & $(0,1)$ &
    \eqref{eq:EM:SAEM:boundbias}  & \eqref{eq:EM:SAEM:boundvar}  & ${\Liplyap}$ & $\vminem$ & $\vmaxem$ \\
    \hline
    TD(0) & $\frac{1}{2} \| \prm - \prm_\star \|^2$ & $\|\Feature \prm - \Feature \prm_\star \|_{\Dtd_{\statdistMRP}}^2$ & $0$ & $(0, (1+\lambda)^2 )$ & $(0,0)$ & \eqref{eq:TD:boundvar_0} & $1$ & $1-\lambda$ & $1/ \sqrt{\vmintd}$ \\
    \bottomrule
    \end{tabular}}
    \caption{Required constants for \Cref{assum:field} and \Cref{assum:lyapunov}, for each example detailed in \Cref{sec:verif:assume}. 
    }\vspace{-.2cm}
\label{tab:constant} 
\end{table*}
\renewcommand{\arraystretch}{1}



\vspace{.1cm}
\noindent \textbf{Smooth (possibly non-convex) objective functions.}
We instantiate the assumptions by setting the Lyapunov function as $\lyap(\prm) \eqdef F(\prm)$ and choosing $\superlyap (\prm)\eqdef \| \nabla F(\prm)\|^2$.
With these choices, \Cref{assum:field}-\ref{item:clyap-hg} is satisfied with $(\clyap_{\hg,0}, \clyap_{\hg,1}) \eqdef (0,1)$, \Cref{assum:lyapunov}-\ref{item:rholyap} is satisfied with $\rholyap \eqdef 1$ since $\pscal{\nabla \lyap(\prm)}{\hg(\prm)}= - \|\nabla F(\prm)\|^2$. Finally, we have $\clyap_{\lyap}\eqdef 1$ in \eqref{eq:definition:clyap-lyap}. To establish the other conditions, namely \Cref{assum:field}-\ref{item:field:conditional-variance} and \Cref{assum:lyapunov}-\ref{item:vstar}, \ref{item:lyap}, we add the following assumption.
\begin{assumSGD}
\label{assum:SGD}
\begin{enumerate}[label=\emph{\alph*}), nosep]
\item \label{item:smoothF} For any $i \in \{1,\dots,n\}$, the function $f_i$ is differentiable and its gradient is $\Lip{\nabla f_i}$-Lipschitz, \ie\  for all $\prm, \prm' \in \rset^d$, $\| \nabla f_i(\prm) - \nabla f_i(\prm') \| \leq \Lip{\nabla f_i} \| \prm -\prm'\|$.
\item \label{item:SGD:boundedF} For any $i \in \{1,\dots,n\}$,  $M \eqdef \max_{i \in \{1,\dots, n\}} \sup_{\prm \in \rset^d} \| \nabla f_i (\prm) - \nabla F(\prm)\| < \infty$.
\end{enumerate}
\end{assumSGD}
Thus \Cref{assum:field}-\ref{item:field:conditional-variance} holds with $\boundvar_0:= M^2/n$ and $\boundvar_1 \eqdef 0$.
Under the above condition, \Cref{assum:lyapunov}-\ref{item:vstar} holds with $\lyap_\star  \eqdef \inf_{\prm \in \rset^d} F(\prm)$.
\Cref{assum:lyapunov}-\ref{item:lyap} holds with $\Liplyap\eqdef n^{-1} \sum_{i=1}^n \Lip{\nabla f_i}$.\vspace{.1cm}



\vspace{.1cm}
\noindent \textbf{Smooth convex objective functions.} When $F$ is convex, we can choose other functions $\lyap, \superlyap$ and obtain different constants in \Cref{assum:field}-\ref{assum:lyapunov}.
We consider
\begin{cvxSGD}
\label{assum:SGD:convex}
The function $F$ is convex with an optimal solution  \ie\ $\mathrm{Argmin}_{\prm \in \rset^d} F(\prm) \neq \emptyset$. 
\end{cvxSGD}
Note that a differentiable function is convex on $\rset^d$ \iff\ for all $\prm,\prm' \in \rset^d \times \rset^d$,
\begin{equation}
\label{eq:condition-convexity}
F(\prm') \geq F(\prm) + \pscal{\nabla F(\prm)}{\prm'-\prm} \eqsp,
\end{equation}
and equivalently, its gradient is monotone
\begin{equation} \label{eq:convex:monotoneGdt}
\pscal{\nabla F(\prm) - \nabla F(\prm')}{\prm - \prm'}
\geq 0 \eqsp;
\end{equation}
see \eg Proposition 17.45 and Proposition 17.7 in \cite{bauschke:combettes:2011}. 
Second, under \Cref{assum:SGD} and \Cref{assum:SGD:convex}, the mapping $\nabla F$ is \emph{co-coercive}  \ie for all $\prm, \prm' \in \rset^d$,
\begin{multline} \label{eq:cocoercive:convex}
\pscal{\nabla F(\prm) - \nabla F(\prm')}{\prm - \prm'} \\
\geq ({1} / {\Lip{\nabla F}} ) \| \nabla F(\prm) - \nabla F(\prm')\|^2 \eqsp;
\end{multline}
see \eg\ Proposition 17.45 and Corollary 18.14 in \cite{bauschke:combettes:2011}. 
 When $F$ is strictly convex, then the inequalities in \eqref{eq:convex:monotoneGdt} and \eqref{eq:cocoercive:convex} are strict (see \eg\ Proposition 17.10. in  \cite{bauschke:combettes:2011}).


Let $\prm_\star \in \rset^d$ be a root of $\nabla F$:  $\nabla F(\prm_\star)=0$. We get from \eqref{eq:cocoercive:convex} that $\|\nabla F(\prm) \|^2 \leq \Lip{\nabla F} \pscal{\nabla F(\prm)}{\prm - \prm_\star}$ for all $\prm \in \rset^d$.
This naturally suggests the definitions
\[
\lyap(\prm) \eqdef  (1/2) \| \prm - \prm_\star \|^2, \quad  \superlyap(\prm) \eqdef \pscal{\nabla F(\prm)}{\prm - \prm_\star}.
\]
This yields
\[
\|\hg(\prm)\|^2= \| \nabla F(\prm) \|^2 \leq \Lip{\nabla F} \superlyap(\prm) \eqsp,
\]
so that  \Cref{assum:field}-\ref{item:clyap-hg} is satisfied with $(\clyap_{\hg,0}, \clyap_{\hg,1})=(0,\Lip{\nabla F})$.  
\Cref{assum:lyapunov}-\ref{item:lyap} trivially holds with $\Liplyap \eqdef 1$ and from
\[
\pscal{\nabla \lyap(\prm)}{\hg(\prm)}= -\pscal{\prm-\prm_\star}{\nabla F(\prm)}= - \superlyap(\prm),
\]
\Cref{assum:lyapunov}-\ref{item:rholyap} is satisfied with $\rholyap=1$.
Moreover, we note that the condition \eqref{eq:definition:clyap-lyap} reads $\| \prm- \prm_\star\| \leq \clyap_{\lyap} \sqrt{\pscal{\nabla F(\prm)}{\prm-\prm_\star}}$; then $\clyap_{\lyap} = \infty$ when $F$ is convex but  $\clyap_{\lyap}$ is finite when $F$ is strongly convex (see below). \vspace{.1cm}

\vspace{.1cm}
\noindent \textbf{Smooth strongly convex objective functions.} We now strengthen the condition on $F$ to strongly convex functions.
\begin{cvxSGD}
\label{assum:SGD:strongly-convex}
The function $F$ is strongly convex with modulus $\mu >0$.
\end{cvxSGD}
In other words, there exists $\mu>0$ such that for any $\prm, \prm' \in \rset^d$ and $\lambda \in \ooint{0,1}$, it holds
\begin{multline*}
\lambda F(\prm) + (1-\lambda) F(\prm') \geq F(\lambda \prm + (1-\lambda) \prm') \\
+  \lambda (1-\lambda) ({\mu}/{2}) \|\prm - \prm'\|^2.
\end{multline*}

Note that a strongly convex function possesses an unique minimizer (see \eg Corollary 11.17  in \cite{bauschke:combettes:2011}), here denoted by $\prm_\star$.  Under \Cref{assum:SGD}, \Cref{assum:SGD:strongly-convex} is equivalent to $\nabla F$ \revisionupdates{being} $\mu$-strongly monotone  \ie\  for all $\prm, \prm' \in \rset^d$,
\begin{equation}
\label{eq:strongly-monotone}
\pscal{\nabla F(\prm') - \nabla F(\prm)}{\prm'-\prm} \geq \mu \| \prm' - \prm \|^2 \eqsp,
\end{equation}
(see \eg\ Definition 2.23 and Exercise 17.5 in \cite{bauschke:combettes:2011}), which is stronger that \eqref{eq:convex:monotoneGdt}.
As in the convex case above, we may again instantiate the assumptions with 
\[
\lyap(\prm) \eqdef  (1/2) \| \prm - \prm_\star \|^2, \quad \superlyap(\prm)\eqdef \pscal{\nabla F(\prm)}{\prm - \prm_\star},
\] 
Since \Cref{assum:SGD:strongly-convex} implies \Cref{assum:SGD:convex}, \eqref{eq:cocoercive:convex} holds and $(\clyap_{h,0}, \clyap_{h,1})  \eqdef (0,\Lip{\nabla F})$. Under \Cref{assum:SGD}, $(\boundvar_0, \boundvar_1) \eqdef (M^2 /n, 0)$. We also have $\Liplyap \eqdef 1$ and $\rholyap \eqdef 1$. Finally, $\clyap_\lyap \eqdef \sqrt{2/\mu}$ in \eqref{eq:definition:clyap-lyap}.

Moreover, an alternative choice for $\superlyap$ is possible, \ie
\[
\superlyap(\prm)= \lyap(\prm) \eqdef  (1/2) \| \prm - \prm_\star \|^2.
\]
Observe that, from \eqref{eq:strongly-monotone},
\[
\pscal{\nabla \lyap(\prm)}{\hg(\prm)}= -\pscal{\prm-\prm_\star}{\nabla F(\prm)} \leq - 2 \mu \superlyap(\prm) \eqsp.
\]
From \Cref{assum:SGD}-\ref{item:smoothF} we get
$\|\nabla F(\prm)\|^2 \leq 2 \Lip{\nabla F}^2 \superlyap(\prm)$. It yields that \Cref{assum:field}-\ref{item:clyap-hg} is satisfied with $(\clyap_{\hg,0}, \clyap_{\hg,1}) \eqdef (0, 2 \Lip{\nabla F}^2)$. 
In addition, under \Cref{assum:SGD}, $(\boundvar_0, \boundvar_1) \eqdef ( M^2/n, 0)$. We also have $\Liplyap \eqdef 1$ and from  \eqref{eq:strongly-monotone}, \Cref{assum:lyapunov}-\ref{item:rholyap} is satisfied with $\rholyap \eqdef 2 \mu$. Finally,  $\clyap_{\lyap} \eqdef \sqrt{2}$ in \eqref{eq:definition:clyap-lyap}.
\subsubsection{Compressed {\SA}} \label{subsec:csgd-checking}
We show that the assumptions are valid for the methods introduced in \Cref{subsec:ex-variantsSGD}. Proofs for this section can be found in \Cref{subsec:proof:csgd-checking}.
We first focus on the compressed {\ SA } in \eqref{eq:cSGD}.
Consider the following assumption about the compressor $\compressor$:
\begin{assumSGDvar}
\label{assum:SGD:contractive}
There exists $\oneminusomgbiased \in \ccint{0,1} $ such that the compression operator $\compressor$ satisfies for any $\bx\in \rset ^d$:
\begin{equation}\label{eq:comp_contractive}
    \E\left[\|\compressor (\bx, \bU) - \bx\|^2\right] \le \omgbiased \|\bx\|^2.
\end{equation}
\end{assumSGDvar}
The Gauss-Southwell estimator (see \eqref{eq:gauss_south}), obtained by using the projection operator $\rm{Top}_1$, satisfies \Cref{assum:SGD:contractive} with $\oneminusomgbiased \eqdef \nicefrac{1}{d}$: note indeed that $ \sum_{i=1}^d x_i^2 \leq d (\max_i x_i^2)$ which implies that $
\|\compressor(\bx,\bU) -\bx \|^2 = \|\bx\|^2 - \max_i x_i^2 \leq (1-1/d) \|\bx\|^2$.
The projection operator $\rm{Top}_h$ satisfies \Cref{assum:SGD:contractive} with $\omgbiased \eqdef 1 -\nicefrac{h}{d}$. Other compression operators satisfy this assumption for various constants $\oneminusomgbiased$, see \eg \cite{stich2018sparsified, karimireddy2019error}, \cite[Table~1]{safaryan2022uncertainty}.

We now discuss \Cref{assum:field} for the compressed oracles. Since compression can be a random operator, we need to adjust the definition of the filtration $\sequence{\mcf}[k][\nset]$ as follows: $\mcg_0 \eqdef \sigma(\prm_0)$ and for all $k \geq 0$
\[
\mcg_{k+1} \eqdef \sigma\left( \prm_0, \State_1, \bU_1, \ldots, \State_{k+1}, \bU_{k+1}\right);
\]
$\bU_{k+1}$ denotes the random variable $\bU$ sampled at iteration $ k+1$ when the compression operator $\compressor$ is applied. With these definitions, $\prm_k \in \mcg_k$ for all $k \geq 0$.
We assume that the oracle $\Hg$ satisfies \Cref{assum:field}. We show that the oracle $\compressor(\Hg(\prm,\State), \bU)$
also satisfies \Cref{assum:field} with different constants.
\begin{lemma} \label{lem:csgd-var1}
Assume that $\Hg$ satisfies \Cref{assum:field}, with constants $(\clyap_{\hg,0}, \clyap_{\hg,1})$, $(\boundbias_{0,k}, \boundbias_{1,k})$ and $(\boundvar_0, \boundvar_1)$ and that $\sup_k (\boundbias_{0,k} + \boundbias_{1,k}) < \infty$. If $\compressor$ satisfies~\Cref{assum:SGD:contractive}, then, the compressed {\SG} in \eqref{eq:cSGD} satisfies \Cref{assum:field} with constants in \Cref{assum:field}~\ref{item:field:boundbias} and \Cref{assum:field}~\ref{item:field:conditional-variance} given for $\ell \in \{0,1\}$ and $k \geq 0$, by
\begin{align}\label{eq:boundbias_res_contractive}
 \boundbias_{\ell,k; \compressor} &\eqdef  ( (1+\zeta_1) +(1+\zeta_2)  (1+\zeta_1^{-1})\omgbiased) \boundbias_{\ell,k}  \nonumber \\ 
 & \quad +(1+\zeta_2^{-1})(1+\zeta_1^{-1}) \omgbiased \clyap_{h,\ell} \nonumber \\ 
 &   \quad+ (1+\zeta_1^{-1}) \omgbiased \boundvar_\ell. \\
 \boundvar_{\ell;\compressor} &\eqdef  \omgbiased \Big( (1+\zeta_2)\sup_k \boundbias_{\ell,k}+  (1+\zeta_2^{-1}) \clyap_{h,\ell}  \Big) \nonumber \\ 
  & \quad +  \omgbiased \boundvar_\ell.
\label{eq:increase_var_contractive} 
\end{align}
for any $\zeta_1, \zeta_2 \in \bar \rset_+$. 
    $\clyap_{h,0}$ and $\clyap_{h,1}$ are unchanged. 
\end{lemma}
An application of \Cref{lem:csgd-var1} is the Gauss-Southwell update.
\begin{corollary}\label{cor:GSascomprGD}
Consider the Gauss-Southwell update~\eqref{eq:gauss_south_choice_i}, with the corresponding gradient field $\Hg(\prm, \sim) = \rm{Top_1}(\nabla F(\prm))$ satisfies \Cref{assum:field} with $(\clyap_{\hg,0;\compressor}, \clyap_{\hg,1;\compressor}) =(0,1)$, $(\boundbias_{0,k;\compressor}, \boundbias_{1,k;\compressor})=(0,1-\nicefrac{1}{d})$ and $(\boundvar_{0;\compressor}, \boundvar_{1;\compressor})=(0,0)$, for  $\superlyap = \|\nabla F (\cdot)\|^2$).
\end{corollary}
Indeed, the $\rm{Top}_1$ compressor satisfies~\Cref{assum:SGD:contractive} with ${\omgbiased}=(1-\nicefrac{1}{d})$ and  the deterministic field $\hg(\prm) \equiv \nabla F(\prm ) $ satisfies \Cref{assum:field}, for $\superlyap = \|\nabla F (\cdot)\|^2$,  with $(\clyap_{\hg,0}, \clyap_{\hg,1}) =(0,1)$, $(\boundbias_{0,k}, \boundbias_{1,k})=(0,0)$ and $(\boundvar_0, \boundvar_1)=(0,0)$. Using \Cref{lem:csgd-var1}, with $(\zeta_1, \zeta_2)=(\infty, \infty)$,  gives the result. 

We note that the compressed random field can be biased even if the original $\Hg$ is not biased: $\boundbias_{\ell,k} =0$ does not imply $\boundbias_{\ell,k; \compressor} = 0$ if there is compression (\ie\ $\omgbiased \neq 0$).
As a workaround, a stronger assumption is \emph{unbiased} compression.
\begin{assumSGDvar}
\label{assum:SGD:URBV}
There exists $\omgunbiased \geq 0$ such that the compression operator $\compressor$ satisfies for any $\bx\in \rset ^d$:
\begin{equation}\label{eq:URVB}
\PE[\compressor (\bx, \bU)] = \bx \eqsp, \quad \PE\left[\|\compressor(\bx, \bU) - \bx\|^2\right] \le \omgunbiased \|\bx\|^2.
\end{equation}
\end{assumSGDvar}
Note that among the operators satisfying this assumption, (scaled) $\rm{Rand}_h$ with $1+\omgunbiased=\nicefrac{d}{h}$, (scaled) $p$-sparsification with $1+\omgunbiased=\nicefrac{1}{p}$, stochastic rounding quantization \eqref{eq:stochastic-rounding} with $\omgunbiased$ as a function of $\Delta$ \cite{alistarh_qsgd_2017,horvoth2022natural}.
\begin{lemma}\label{lem:csgd-var2}
Assume that $\Hg$ satisfies \Cref{assum:field}, with constants $(\clyap_{\hg,0}, \clyap_{\hg,1})$, $(\boundbias_{0,k}, \boundbias_{1,k})$ and $(\boundvar_0, \boundvar_1)$ and that $\sup_k (\boundbias_{0,k} + \boundbias_{1,k}) < \infty$. If $\compressor$ satisfies~\Cref{assum:SGD:URBV}, then, the compressed {\SA} in \eqref{eq:cSGD} satisfies \Cref{assum:field} with constants in \Cref{assum:field}~\ref{item:field:boundbias} and \Cref{assum:field}~\ref{item:field:conditional-variance} given for $\ell \in \{0,1\}$ and $k \geq 0$, by $\boundbias_{\ell,k;\compressor} \eqdef \boundbias_{\ell,k}, \clyap_{\ell,k;\compressor} \eqdef \clyap_{\ell,k}$, and
\begin{equation}\label{eq:increase_var_URBV}
        \boundvar_{\ell; \compressor} \eqdef (1+\omgunbiased) \boundvar_{\ell} + 2 \omgunbiased (\clyap_{\hg,\ell} + \sup_k \boundbias_{\ell,k} ).
    \end{equation}
\end{lemma}
Without compression ($\omgunbiased=0$), the constants remain unchanged; using compression introduces additional variance.
Next, we consider the following result for the compressed {\SA } with perturbed iterate \eqref{eq:binary-connect}.
We introduce a third assumption about compression operators that covers the case of a uniformly bounded quantization error in space.
\begin{assumSGDvar}
\label{assum:SGD:UnifBoundedVar}
There exists $\omguniform \geq 0$ such that the compression operator $\compressor$ satisfies for any $\bx\in \rset ^d$:
\begin{equation}\label{eq:unif_bounded_C}
  \PE\left[\|\compressor(\bx, \bU) - \bx\|^2\right] \le \omguniform.
\end{equation}
\end{assumSGDvar}
Such an assumption is satisfied for operators satisfying \Cref{assum:SGD:contractive} or \Cref{assum:SGD:URBV} on a bounded domain $\mathcal{X}\subset \rset^d $, or by using an adaptive number of bits to compress the signal, depending on its scale. This assumption was used, for example, in \cite{Zheng19,tang2019doublesqueeze}. For deterministic rounding, as defined in \eqref{eq:deterministic-rounding}, with a quantization step $\Delta$, it holds with $\omguniform= d\Delta^2$. Such an assumption can be adapted to handle \emph{asynchrony}, for a.s. bounded fields and bounded delays (the compression scheme can be defined \emph{only} on the points to which it is applied, and depends on the previous sequence of iterates)~\cite{agarwal2011distributed}.
\begin{lemma} \label{lem:csgd-var3}
Assume that $\Hg$ satisfies \Cref{assum:field} with constants $(\clyap_{\hg,0}, \clyap_{\hg,1})$, $(\boundbias_{0,k}, \boundbias_{1,k})$ and $(\boundvar_0, \boundvar_1)$ such that $\sup_k (\boundbias_{0,k} ) < \infty$,  for any $k\geq 0$, $ \boundbias_{1,k} = 0$,  $\boundvar_1 =0 $. Assume that $\compressor$ satisfies~\Cref{assum:SGD:UnifBoundedVar},  that  $\hg$ is $\Lip{\hg}$ Lipschitz,  and that $\CPE{ \Hg( \cdot , \State_{k+1}) }{\mcg_{k}} $ is $\Lip{\PE{\Hg}}$-Lipschitz. Then the compressed {\SA} in \eqref{eq:binary-connect}
satisfies \Cref{assum:field}-\ref{item:field:boundbias} and \Cref{assum:field}-\ref{item:field:conditional-variance}, for all $k\geq 0$, for any $\zeta\in \bar\rset_{+}$ with $\boundbias_{1,k,\compressor} \eqdef 0 $, $\boundvar_{1;\compressor} \eqdef 0$, and 
\begin{align*}
\boundbias_{0,k,\compressor} &\eqdef(1+\zeta) \boundbias_{0,k}  + (1+\frac{1}{\zeta})\Lip{\hg}^2 \omguniform,  \ \ 
\boundvar_{0; \compressor}  \eqdef   \boundvar_{0}+  \, \Lip{\PE{\Hg}}^2\omguniform.
\end{align*}
\end{lemma}
In other words, the {\SA} scheme with perturbed iterates results in an additional bias and variance term.
We get the following corollary for compressed {\SG} algorithm STE. 
\begin{corollary}\label{cor:STEfield}
Let $\Hg$ be the oracle given by  \eqref{eq:SGD-finite-sum} and consider the STE compression {\SG} algorithm~\eqref{eq:binary-connect}. Assume \Cref{assum:SGD}-\ref{item:smoothF} and $\compressor$ satisfies~\Cref{assum:SGD:UnifBoundedVar}. Then the resulting field satisfies \Cref{assum:field}-\ref{item:field:boundbias} and \Cref{assum:field}-\ref{item:field:conditional-variance}, for all $k\geq 0$, with $\boundbias_{1,k,\compressor} \eqdef 0 $, $\boundvar_{1;\compressor} \eqdef 0$, and 
\begin{align*}
\boundbias_{0,k,\compressor} &\eqdef \Lip{{\nabla F}}^2 \omguniform,  \qquad
\boundvar_{0; \compressor} \eqdef   \boundvar_{0}+  \, \Lip{\nabla F}^2\omguniform.
\end{align*}
\end{corollary}
Indeed, we have $\hg=\nabla F$ and an unbiased stochastic oracle, thus $\CPE{ \Hg( \cdot , \State_{k+1}) }{\mcg_{k}} =\nabla F$.

Lastly, we study the low-precision {\SA} introduced in \eqref{eq:de-sa} where we work with the compressor satisfying:
\begin{assumSGDvar}
\label{assum:SGD:UniformQuant}
There exists $\linearunbiased \geq 0$ such that for any $\bx \in \rset^d$, the compression operator $\compressor$ satisfies: (i) $\PE[ \compressor( \bx , \bU ) ] = \bx$, (ii) denote $\linearset^d$ as the image of $\compressor(\cdot)$ and it holds for any $\bar{\bx} \in \linearset^d$,  
$\compressor( \bar{\bx} , \bU ) = \bar{\bx}$ almost surely,
and for any $\bv \in \rset^d$, $\bx \in \rset^d$, 
\begin{equation}\label{eq:linearQBound}
  \PE\left[\|\compressor( \bar{\bx} + \bv , \bU) - \bx\|_2^2\right] \le \| \bar{\bx} + \bv - \bx \|_2^2 + \linearunbiased \| \bv \|_1.
\end{equation}
\end{assumSGDvar}
The above assumption is valid for the stochastic rounding quantizer in \eqref{eq:stochastic-rounding} with $\linearunbiased = \Delta$, 
as shown in \cite[Lemma 3]{Zheng19} with $\delta = \linearunbiased, b = \infty$ therein. Notice that it is also known as a \emph{linear quantizer} in \cite{Zheng19}.
\begin{lemma}\label{lem:csgd-var4}
Assume $\Hg$ satisfies \Cref{assum:field}, with constants $(\clyap_{\hg,0}, \clyap_{\hg,1})$, $(\boundbias_{0,k}, \boundbias_{1,k})$ and $(\boundvar_0, \boundvar_1)$ and that $\sup_k (\boundbias_{0,k} + \boundbias_{1,k}) < \infty$. Consider \eqref{eq:de-sa}, \eqref{eq:de-sa-H} with constant stepsize $ \step_{k} = \bar{\step}$ for all $k$.
Assume that $\compressor$ satisfies \Cref{assum:SGD:UniformQuant}. Then, the random field in \eqref{eq:de-sa-H} satisfies \Cref{assum:field} with constants in \Cref{assum:field}~\ref{item:field:boundbias} and \Cref{assum:field}~\ref{item:field:conditional-variance} given for $\ell \in \{0,1\}$ and $k \geq 0$, by $\boundbias_{\ell,k;\compressor} \eqdef \boundbias_{\ell,k}$ and 
\begin{equation}\label{eq:increase_var_linearQuant}
\begin{aligned}
\boundvar_{0; \compressor} & \textstyle \eqdef \boundvar_0 + \frac{ \linearunbiased \sqrt{d} }{ 2 \bar{\step} } \left( 3 + \sup_{k \geq 0} \boundbias_{0,k} + \clyap_{\hg,0} + \boundvar_0  \right), \\
\boundvar_{1; \compressor} & \textstyle \eqdef \boundvar_1 + \frac{ \linearunbiased \sqrt{d} }{ 2 \bar{\step} } \left( \sup_{k \geq 0} \boundbias_{1,k} + \clyap_{\hg,1} + \boundvar_1 \right),
\end{aligned}
\end{equation}
while $\clyap_{h,0}$ and $\clyap_{h,1}$ are unchanged.
\end{lemma}
We note that the random field  in \eqref{eq:de-sa-H} inherits the same bias properties from $\Hg$, while its variance is increased to ${\cal O}( 1 + \frac{\linearunbiased \sqrt{d} }{ \bar{\step} } )$. Note that the variance is now inversely proportional to the step size. As we  show below in \Cref{subsub:exrate:CSGD}, such a compressed {\SA} Algorithm converges only to an ${\cal O}( \linearunbiased \sqrt{d} )$ approximate stationary solution.

\subsubsection{Stochastic EM algorithms}
\label{sec:EM:verifH} We
use the definitions and notations introduced
in \Cref{sec:EM:introduction}.
The Stochastic EM  algorithms solve the  root-finding problem  for the function
\begin{equation}
\hg(\prm)  \eqdef \barsem \circ \mapem(\prm) - \prm; \label{eq:EM:meanfield}
\end{equation}
see \eqref{eq:EM:statisticspace}.
Our choice for the Lyapunov function is:
\begin{equation} \label{eq:EM:lyapunovfct}
\lyap(\prm) \eqdef F \circ \mapem(\prm).
\end{equation}
In addition to \Cref{assumEM:expo,assumEM:Mstep},  consider the following assumption
\begin{assumEM} \label{assumEM:supplementary}
\begin{enumerate}[label=\alph*),nosep]
\item 
There exists $F_\star>-\infty$ such that $F(\prmo) \geq F_\star$  for any $\prmo \in \rset^d$.
\item  
The function $\lyap$ is continuously differentiable on $\rset^d$ and there exists  $\Liplyap \in \rset_+$ such that for any $\prm, \prm' \in \rset^d$,
\[
\left\| \nabla \lyap(\prm) - \nabla  \lyap(\prm') \right\| \leq \Liplyap \| \prm - \prm' \|.
\]
\item 
For any $\prm \in \rset^d$, there exists a $d \times d$ positive definite matrix $\Bem(\prm)$ such that
$\nabla \lyap(\prm) = - \Bem(\prm) \hg(\prm)$. In addition, there exist positive constants $\vminem \leq \vmaxem$ such that  for any $\prm, \prm' \in \rset^d$,
\[ \vminem \| \prm' \|^2 \leq \| \prm'\|^2_{\Bem(\prm)} \leq \vmaxem \| \prm'\|^2 \, .\]
\end{enumerate}
\end{assumEM}
Lemma~2 in \cite{delyon1999convergence} provides sufficient conditions
 on the regularity of the functions $\phiem, \psiem, \sem_i$ and
 $\mapem$, implying that
 \[
\Bem(\prm) = \left( \nabla \mapem(\prm) \right)^{\top} \, \Deriv_{02} \mathsf{L}(\prm, \mapem(\prm)) \, \left( \nabla \mapem(\prm) \right),
\]
where $\mathsf{L}(\prm, \prmo) \eqdef \psiem(\prmo) - \pscal{\prm}{\phiem(\prmo)}$. Since under \Cref{assumEM:Mstep}, $\mapem(\prm)$ is the unique minimizer of $\prmo \mapsto \mathsf{L}(\prm, \prmo)$, $\Bem(\prm)$ is positive semi-definite. $\nabla \mapem(\prm)$ is a $d  \times d$ matrix and positive definiteness results from conditions on the rank of $\nabla \mapem(\prm)$. 

\revisionupdates{\begin{lemma} 
Under \Cref{assumEM:expo,assumEM:Mstep,assumEM:supplementary},  for any $\prm \in \rset^d$,
$\lyap(\prm)  \geq F_\star$, $\pscal{\nabla \lyap(\prm)}{\hg(\prm)} \leq - \vminem \| \hg(\prm) \|^2$, and
$\| \nabla \lyap(\prm) \|^2  \leq \vmaxem^2 \| \hg(\prm) \|^2$. Thus  \Cref{assum:lyapunov}   is  satisfied with $\lyap_\star \eqdef F_\star$ and
\begin{equation}
\label{eq:constants-minibatch-EM}
\superlyap(\prm) \eqdef \| \hg(\prm) \|^2, \quad  \rholyap \eqdef \vminem, \quad \clyap_{\lyap} \eqdef \vmaxem;
\end{equation}
 \Cref{assum:field}-\ref{item:clyap-hg} is satisfied with $\clyap_{\hg,0} \eqdef 0$ and $\clyap_{\hg,1} \eqdef 1$.
\end{lemma}}    
Let us now check \Cref{assum:field}-\ref{item:field:boundbias}  and \Cref{assum:field} -\ref{item:field:conditional-variance} for some specific examples of stochastic field $\Hg$.

\paragraph{Mini-batch EM}
\label{sec:EM:verifH:onlineEM} 
We recall from \Cref{sec:EM:intro-MinibatchEM} the form of mini-batch EM, \revisionupdates{see \cref{eq:onlineEM:oracle}. 
We consider the  following condition, and obtain the subsequent lemma.}
\begin{assumEM} \label{assumEM:var:onlineEM}
There exist $ \boundvarem_{0}, \boundvarem_{1} \in \rset_+$ such that  for any $\prm \in \rset^d$
\[ 
\sup_{i \in \{1, \ldots, n\}} \| \barsem_i(\mapem(\prm))\|^2 \leq  \boundvarem_0 + \boundvarem_1 \superlyap(\prm).
\]
\end{assumEM}

\revisionupdates{
\begin{lemma}
    Under \Cref{assumEM:expo,assumEM:Mstep,assumEM:supplementary,assumEM:var:onlineEM}, for $\Hg$ given by \eqref{eq:onlineEM:oracle}, \Cref{assum:field}-\ref{item:field:boundbias}  and \Cref{assum:field} -\ref{item:field:conditional-variance} are verified with   $\boundbias_{\ell,k} = 0$ for all $k\in \nset$, and  $\boundvar_\ell \eqdef \nicefrac{\boundvarem_\ell}{\lbatchem}$, for $\ell\in \{1,2\}$.
\end{lemma}}

\revisionupdates{We refer e.g. to \cite[Lemma~7.1.]{fort:moulines:2022}: first, 
$\CPE{\Hg(\prm_k, \State_{k+1})}{\mcf_k} = \hg(\prm_k)$, thus  \cref{eq:onlineEM:oracle} provides an unbiased stochastic oracle. Moreover, the conditional variance of $\Hg(\prm_k, \State_{k+1})$ is upper bounded by 
\[ \textstyle
\frac{1}{\lbatchem} \frac{1}{n} \sum_{i=1}^n \| \barsem_i(\mapem(\prm_k)) - \frac{1}{n}\sum_{j=1}^n \barsem_j(\mapem(\prm_k)) \|^2.
\]
Second,  using $ \sum_{i=1}^n \|a_i - n^{-1} \sum_{j=1}^n a_j \|^2 \leq \sum_{i=1}^n \|a_i\|^2$,  with  \Cref{assumEM:var:onlineEM},  \Cref{assum:field}-\ref{item:field:conditional-variance} is satisfied with $\boundvar_\ell \eqdef \nicefrac{\boundvarem_\ell}{\lbatchem}$.}

\paragraph{SAEM} \label{sec:EM:verifH:SAEM}
\revisionupdates{We now focus on SAEM, defined in \eqref{eq:oracle-SAEM}.} First, consider the case when conditionally to the past, the random variables $\{Z_{i,k+1}^j, 1  \leq j \leq \nbrMC, 1 \leq i \leq n \}$ are independent, and for all $i \in \{1, \ldots, n\}$ and $j$, the distribution of $Z_{i,k+1}^j$ is $\pi_{i}(z_i;  \mapem(\prm_k))$.
Then
\[
 \CPE{ \rho_{i}^j(\prm_k) \, \sem_i(Z_{i,k+1}^j)}{\mcf_k} = \frac{1}{\nbrMC }\barsem_i(\prm_k)
 \]
and the SAEM algorithm is an
unbiased SA. Moreover, the conditional variance of $\Hg(\prm_k, \State_{k+1})$ is equal to
\[
\frac{1}{n^2 \nbrMC} \sum_{i=1}^n \CPE{\| \sem_i(Z_{i,k+1}^1) - \frac{1}{n} \sum_{j=1}^n \barsem_j(\prm_k) \|^2 }{\mcf_k}.
\]
Following the same lines as in the discussions about the Mini-batch EM, we consider the following condition.
\begin{assumEM} \label{assumEM:var:SAEM:exact}
There exist constants $ \boundvarem_{0}, \boundvarem_{1} \in \rset_+$ such that for any $k \in \nset$, almost-surely,
\[
\sup_{i \in \{1, \ldots, n\}} \CPE{\| \sem_i(Z_{i,k+1}^1)\|^2}{\mcf_k} \leq  \boundvarem_0 + \boundvarem_1 \superlyap(\prm_k).
\]
\end{assumEM}
\noindent
\revisionupdates{We have the following lemma.}
\revisionupdates{\begin{lemma}
  Assume that conditionally to the past, the random variables $\{Z_{i,k+1}^j, 1  \leq j \leq \nbrMC, 1 \leq i \leq n \}$ are independent, and for all $i \in \{1, \ldots, n\}$ and $j$, $Z_{i,k+1}^j \sim \pi_i(z_i; \mapem(\prm_k))$. Assume also \Cref{assumEM:expo,assumEM:Mstep,assumEM:supplementary,assumEM:var:SAEM:exact}. Then the oracle given by \eqref{eq:oracle-SAEM} satisfies 
\Cref{assum:field}-\ref{item:field:boundbias} and \Cref{assum:field}-\ref{item:field:conditional-variance} with   $\boundbias_{\ell,k} = 0$, and $\boundvar_\ell \eqdef \boundvarem_\ell/(n \nbrMC)$, for $\ell \in \{1,2\} $ and any $k\in \nset$.
\end{lemma}}

\revisionupdates{Finally,} let us now consider the  self-normalized Importance Sampling case; conditionally to the past, the random variables $\{Z_{i,k+1}^j, 1  \leq j \leq \nbrMC, 1 \leq i \leq n \}$ are independent, and for all $i \in \{1, \ldots, n\}$ and $j$, the distribution of $Z_{i,k+1}^j$ is $\tilde \pi_{i}(z_i; \mapem(\prm_k))$. In that case, SA is not unbiased due to the use of the self-normalized importance weights
\[
\CPE{\Hg(\prm_k, \State_{k+1})}{\mcf_k} \neq \hg(\prm_k).
\]
The expression of the bias and variance is complicated for general functions $\sem_i$ and for simplicity, we assume here that the functions $\sem_i$ are bounded (see \cite[Theorem~2.3]{agapiou:etal:2017} for an in depth study of self-normalized importance sampling). Define the second moment of the importance ratio  with respect to the proposal $\tilde{\pi}_i(\cdot; \mapem(\prm))$
\[
\chi_i(\prm) \eqdef  \int_\Zset \left( \frac{\pi_i(z_i; \mapem(\prm))}{\tilde \pi_i(z_i; \mapem(\prm))}\right)^2  \tilde \pi_i(z_i; \mapem(\prm)) \mu(\rmd z_i).
\]
Consider the following assumptions.
\begin{assumEM} \label{assumEM:var:SAEM:IS} $s_\star \eqdef \max_{1 \leq i \leq n} \sup_{z \in \Zset} \| \sem_i(z) \|$ is finite and there exist constants $ \clyap_{\chi,0} , \clyap_{\chi,1} \in \rset_+$ such that for any $\prm \in \rset^d$
\[ \textstyle
\left(\frac{1}{n}\sum_{i=1}^n \chi_i(\prm) \right)^2\leq \clyap_{\chi,0} + \clyap_{\chi,1} \, \superlyap(\prm).
\]
\end{assumEM}
From \cite[Theorem 2.1]{agapiou:etal:2017}, it holds
\begin{align*}
& \| \CPE{\Hg(\prm_k, \State_{k+1})}{\mcf_k} - \hg(\prm_k) \|   \leq s_\star \frac{12}{\nbrMC} \frac{1}{n}\sum_{i=1}^n \chi_i(\prm), \\
& \CPE{\|\Hg(\prm_k, \State_{k+1}) - \CPE{\Hg(\prm_k, \State_{k+1})}{\mcf_k} \|^2}{\mcf_k}  \\
& \qquad \leq s_\star^2 \frac{4}{\nbrMC} \frac{1}{n^2}\sum_{i=1}^n \chi_i(\prm),
\end{align*}
\revisionupdates{from which we deduce the following Lemma.
\begin{lemma} 
 Under \Cref{assumEM:expo,assumEM:Mstep,assumEM:supplementary,assumEM:var:SAEM:IS}, for $\Hg$ given by \eqref{eq:oracle-SAEM},
\Cref{assum:field}-\ref{item:field:boundbias} and \Cref{assum:field}-\ref{item:field:conditional-variance}  are satisfied with
\begin{align}
&\boundbias_{\ell,k} \eqdef  s_\star^2 \frac{144}{\nbrMC^2} \clyap_{\chi,\ell},  \qquad \ell=0,1  \label{eq:EM:SAEM:boundbias}   \\
&\boundvar_{0} \eqdef s_\star^2 \frac{4}{n \nbrMC} \sqrt{\clyap_{\chi,0} + \clyap_{\chi,1}}, \quad \boundvar_{1} \eqdef s_\star^2 \frac{4}{n \nbrMC} \sqrt{ \clyap_{\chi,1}}. \label{eq:EM:SAEM:boundvar}
\end{align}
\end{lemma}}

Finally, the case when the samples $\{Z_{i,k+1}^j, 1 \leq j \leq \nbrMC \}$ are the path of  an ergodic Markov chain with
invariant distribution $\pi_i(z_i; \mapem(\prm_k))$, is more complex. We have again
\[
\CPE{\Hg(\prm_k, \State_{k+1})}{\mcf_k} \neq \hg(\prm_k).
\]
An expression of the bias and of the conditional variance of $\Hg$ can be found in \cite[Proposition 5]{atchade:fort:moulines:2015} (see also \cite[Section~6]{fort:moulines:2003}) in terms of the iterates of the Markov kernels of the Markov Chain Monte Carlo samplers; these controls rely on Markov Chain theory results (see e.g. \cite{douc2018markov}) whose exposition is out of the scope of this paper.


\subsubsection{TD Learning}
\label{subsec:TD-learning-checking}
We use the definitions and notations of \Cref{sec:definition-TD-learning}. We follow \cite{bhandari2018finite} to set up the following assumptions and the missing proofs are relegated to \Cref{subsec:proof:TD-learning-checking}.

\revisionupdates{To design a suitable function $\superlyap$ for our assumptions about the {\SA} scheme, we first observe} the following lemma from \cite[Lemma 4]{tsitsiklis1997analysis}, which shows
that the projected Bellman operator $\Proj_{\statdistMRP} \bellman$ is a contraction. Note that under \Cref{assum:TD:stationary-policy}, $\Dtd_{\statdistMRP}$ is positive-definite.
\begin{lemma}[\cite{tsitsiklis1997analysis}]
\label{lem:contraction-projected-bellman}
Assume \Cref{assum:TD:stationary-policy,assum:TD:sampling}.   Then, $\Proj_{\statdistMRP} \bellman$ is a contraction with respect to $\|\cdot\|_{\Dtd_\statdistMRP}$ with modulus $\lambda \in (0,1)$, \ie\ for all $\prm, \prm' \in \rset^d$, we get
\begin{equation}
\label{eq:contraction-projected-bellman}
\left\|\Proj_{\statdistMRP} \bellman \valuefunc[\prm]-\Proj_{\statdistMRP} \bellman \valuefunc[\prm'] \right\|_{\Dtd_\statdistMRP} \leq \lambda \left\|\valuefunc[\prm]-\valuefunc[\prm']\right\|_{\Dtd_\statdistMRP}.
\end{equation}
\end{lemma}
Therefore, under \Cref{assum:TD:stationary-policy,assum:TD:sampling},  there exists an unique value function $\valuefunc[\star]$ in  $\Span(\Feature)$ which solves the fixed point of the projected Bellman equation \eqref{eq:projected-bellman}: the TD(0) algorithm can be interpreted as a simple {\SA} scheme for solving the projected Bellman fixed point equation.
Finally, we set
\begin{equation}
\label{eq:lyapunov-W-TD}
\superlyap(\prm) \eqdef \left\|\valuefunc[\prm]-\valuefunc[\star] \right\|_{\Dtd_{\statdistMRP}}^2 \eqsp,
\end{equation}
which measures the difference between the approximate value function with parameter $\prm$ and the function $\valuefunc[\star]$.
\revisionupdates{With the above setup, we are ready to establish \Cref{assum:field} under \Cref{assum:TD:stationary-policy,assum:TD:sampling}.}
It has been checked (see \Cref{sec:definition-TD-learning}) that under \Cref{assum:TD:sampling}, $\CPE{\Hg(\prm_k, \State_{k+1})}{\mcf_k} = \hg(\prm_k)$; hence, \Cref{assum:field}-\ref{item:field:boundbias} is satisfied with $\boundbias_{0,k}=\boundbias_{1,k} \eqdef 0$. \revisionupdates{Moreover, it can be checked that \Cref{assum:field}-\ref{item:field:boundsecondmoment} holds, and \Cref{assum:field}-\ref{item:clyap-hg} holds with $(\clyap_{\hg,0}, \clyap_{\hg,1}) = (0, (1+\lambda)^2)$.} The other conditions require some technical work, which is summarized below:
\begin{lemma}
\label{lem:TD:bound-mean-field}
Assume \Cref{assum:TD:stationary-policy,assum:TD:normed:calR,assum:TD:sampling}.  Then, for all $\prm \in \rset^d$  we get $\| \hg(\prm) \|^2 \leq (1+\lambda)^2  \, \superlyap(\prm)$
and for any $k \geq 0$, it holds almost-surely,
\begin{align*}
&\CPE{\| \Hg(\prm_{k},\State_{k+1}) - \hg(\prm_{k})\|^2}{\mcf_{k}}
\leq \boundvar_0 + \boundvar_1  \superlyap(\prm_k),
\end{align*}
\begin{align}
\label{eq:TD:boundvar_0}
\boundvar_0 & \eqdef 6  \left( 1 +     \{ \lambda^2 +1 \} \| \valuefunc[\star]\|^2_{\Dtd_{\statdistMRP}} \right),~~\boundvar_1 \eqdef 2 (1+\lambda)^2. 
\end{align}
\end{lemma}
We now check \Cref{assum:lyapunov}. Let $\prm_\star \in \rset^d$ be such that $\valuefunc[\star] = \Feature \prm_\star$; such a point exists since $\valuefunc[\star] \in \Span(\Feature)$. Define
\begin{equation}
\label{eq:lyapunov-V-TD}
\lyap(\prm)  \eqdef  (1/2) \| \prm- \prm_\star\|^2 \eqsp.
\end{equation}
Then \Cref{assum:lyapunov}-\ref{item:vstar}, \ref{item:lyap} are verified with $\lyap_\star \eqdef 0$ and $\Liplyap \eqdef 1$.
The following Lemma, borrowed from \cite[Lemma~3]{bhandari2018finite}, shows that \Cref{assum:lyapunov}-\ref{item:rholyap} is satisfied with 
$\rholyap \eqdef 1 - \lambda$.
\begin{lemma}
\label{lem:TD:lyap}
Assume \Cref{assum:TD:stationary-policy,assum:TD:sampling}.
For any $\prm \in \rset^d$,
$$
\pscal{\nabla \lyap(\prm)}{\hg(\prm)} \leq -(1-\lambda)  \superlyap(\prm)\eqsp,
$$
where $\hg$ and $\lyap$ are defined in \eqref{eq:TD:mean-field} and \eqref{eq:lyapunov-V-TD} respectively.
\end{lemma}
Finally,  let us check \eqref{eq:definition:clyap-lyap}. Since
$\superlyap(\prm) = (\prm -\prm_\star)^{\top} \Feature^\top \Dtd_{\statdistMRP} \Feature(\prm - \prm_\star)$,
we get (see \Cref{subsec:proof:TD-learning-checking} for a proof)
\begin{lemma}
\label{lem:TD:upper-lower-spectrum-covariance}
Assume \Cref{assum:TD:stationary-policy}.
Then the minimal eigenvalue $\vmintd$ of $\Feature^\top \Dtd_{\statdistMRP} \Feature$ is non-negative and for all $\prm \in \rset^d$, $\sqrt{\vmintd} \| \prm \| \leq \| \prm \|_{\bSigma_\statdistMRP} \leq  \|\prm\|$.
\end{lemma}
\noindent Hence we  set $\clyap_\lyap \eqdef 1/ \sqrt{\vmintd}$ in \eqref{eq:definition:clyap-lyap}  which can be $+\infty$.

Let us revisit this discussion under stronger conditions on the feature matrix $\Feature$. We assume in the sequel that the number of parameters $d$ needed to approximate the value function $\valuefunc[\prm] \eqdef\Feature \prm$, is smaller than or equal to the number of states $n$ and that any redundant or irrelevant feature has been removed from the feature matrix $\Feature$ defined in \eqref{eq:definition_feature-matrix}. Additionally, we assume that the features are normalized. This is formalized in the following assumption:
\begin{assumTD}
\label{assum:TD:full-rank}
The feature matrix is full rank with $\rank(\Feature)=d$. In addition, for all $s \in \stateMRP$, $\| \feature(s) \| \leq 1$.
\end{assumTD}
Under \Cref{assum:TD:stationary-policy,assum:TD:sampling,assum:TD:full-rank}, there exists an unique vector $\prm_\star \in \rset^d$ such that $\valuefunc[\star] = \Feature \prm_\star$, and  $\prm_\star$ is the unique root to the mean field $\hg$. Under \Cref{assum:TD:full-rank}, the \emph{feature covariance matrix}
\begin{equation}
\label{eq:covariance-feature}
\bSigma_{\statdistMRP} \eqdef \PE_{\statdistMRP}[ \feature(S_0) \feature^\top(S_0)]= \Feature^{\top} \Dtd_\statdistMRP \Feature \eqsp,
\end{equation}
is positive-definite: \revisionupdates{we have an equivalent expression for $\superlyap(\prm)$},
\[
\superlyap(\prm)  =\left\|\prm-\prm_\star\right\|_{\bSigma_{\statdistMRP}}^2,
\]
and $\clyap_\lyap = 1/ \sqrt{\vmintd} > 0$.  Finally, 
$(\boundvar_0, \boundvar_1)$ are given by \eqref{eq:TD:boundvar_0}.



\section{Non-asymptotic convergence bounds} \label{sec:nonasymp}
As our first attempt on the theoretical analysis for {\SA}, we overview the non-asymptotic convergence bounds that estimates some properties of the iterates after running the {\SA} scheme for a certain number of iterations. Throughout this section, we focus on the case where the number of iterations of the algorithm is bounded by $\totstep < \infty$, the optimization horizon. 
Note that the bounds to be presented are  in the form of \emph{expected convergence}, where upper bounds on the expected values of the function $\superlyap$ are obtained after $\totstep$ iterations. 

\subsection{Finite-time bounds and sample complexity of SA}
\label{subsec:finite-time-bounds}
We introduce two simplifying assumptions to streamline our forthcoming discussions.
\begin{assumNA}
\label{assumNA:uniform-bound-bias}
There exist constants $\boundbias_0, \boundbias_1  \in \rset_+$ such that for  any $k \in \nset$, $\boundbias_{0,k}= \boundbias_0$, $\boundbias_{1,k}= \boundbias_1$.
\end{assumNA}
\begin{assumNA}
\label{assumNA:bias-is-small} It holds
$\clyap_{\lyap} (\sqrt{\boundbias_0} /2 + \sqrt{\boundbias_1}) < \rholyap$, where $\clyap_{\lyap}$ is in \eqref{eq:definition:clyap-lyap}.
\end{assumNA}
Recall that $\clyap_{\lyap}$ can be infinite, and in  this case,  it is necessary to have $\boundbias_0= \boundbias_1=0$ for verifying  \Cref{assumNA:bias-is-small}.
In words, the above assume that the  bias in the oracle $\Hg( \prm_k, \State_{k+1} )$ is uniformly bounded and small \wrt  the strength of the drift $\rholyap$.
Define
\begin{align}
\label{eq:definiton-bsf-0}
\bsf_0 &\eqdef \clyap_{\lyap} \sqrt{\boundbias_0}/2, \quad \bsf_1 \eqdef \clyap_{\lyap} (\sqrt{\boundbias_0} /2 + \sqrt{\boundbias_1}).
\end{align}
For USO (see \Cref{def:USO}), $\bsf_0 = \bsf_1 = 0$.
Set  for $\ell \in \{0,1\}$
\begin{align}
\label{eq:definition-eta}
\eta_\ell& \eqdef \boundvar_\ell + \boundbias_\ell +  \clyap_{\hg,\ell} + \sqrt{\clyap_{\hg,\ell}} \left( \sqrt{\boundbias_0 }+ \sqrt{\boundbias_1 }\right) \\
\nonumber
& \qquad + \sqrt{\boundbias_\ell } \left( \sqrt{\clyap_{\hg,0}} + \sqrt{\clyap_{\hg,1}}\right) \, , \\
\label{eq:step-max}
\step_{\max} & \eqdef  2 \{\rholyap- \bsf_1\}/(\Liplyap \eta_1) \, ,   \\
\label{eq:definition-u-k}
\omega_k & \eqdef 2\{ \rholyap - \bsf_1\} - \step_{k} \Liplyap \eta_1  \eqsp.
\end{align}
If $\eta_1 = 0$, then by convention, $\step_{\max} \eqdef +\infty$. By construction and under \Cref{assumNA:bias-is-small}, if $\step_k \in \ooint{0, \step_{\max}}$, then $\omega_k > 0$.

The essential argument of our analysis is a Robbins-Siegmund type inequality \cite{robbins1971convergence}, which has played an essential role in the theory of {\SA} since the first works on this subject. Roughly speaking, we control the (expected) changes of the Lyapunov function value $\lyap( \prm_{k+1} ) - \lyap( \prm_k )$ in relation to  $\superlyap ( \prm_k )$, the stepsizes, bias, etc.
\begin{lemma}[Robbins-Siegmund type inequality]
\label{lem:main-result-quantitative}
Assume \Cref{assum:field,assum:lyapunov} and \Cref{assumNA:uniform-bound-bias}.
Then, for any $k \geq 0$, we have almost-surely
\begin{align}
\nonumber
\CPE{\lyap(\prm_{k+1})}{\mcf_k} & \leq  \lyap(\prm_k)
- (1/2) \step_{k+1} \omega_{k+1}
\, \superlyap(\prm_k) \\
\label{eq:condition-lyap}
&+  \step_{k+1} \bsf_0 + \step_{k+1}^2 \Liplyap  \eta_0/2 \,.
\end{align}
\end{lemma}
\begin{proof}
Let $k \geq 0$. 
By \Cref{assum:lyapunov}-\ref{item:lyap},  we have
\begin{align} 
\lyap(\prm_{k+1}) \leq \lyap(\prm_k) + \pscal{\nabla \lyap(\prm_k)}{\prm_{k+1}-\prm_k} \nonumber
\\ + (\Liplyap/2) \|\prm_{k+1} - \prm_k \|^2 \,. \label{eq:robinsiegmund-Lem9}
\end{align}
Define $\bias_k \eqdef \CPE{\Hg(\prm_k, \State_{k+1})}{\mcf_k} - \hg(\prm_k)$. Computing the conditional expectation of both sides of \eqref{eq:robinsiegmund-Lem9} yields:
\begin{align*}
&\CPE{\lyap(\prm_{k+1})}{\mcf_k} \leq \lyap(\prm_k) \\
&\quad + \step_{k+1} \pscal{\nabla \lyap(\prm_k)}{\hg(\prm_k)} + \step_{k+1} \pscal{\nabla \lyap(\prm_k)}{\bias_k}\\
&\quad + \step^2_{k+1} (\Liplyap /2) \,  \CPEs{\| \Hg(\prm_k,\State_{k+1}) \|^2}{\mcf_k} \,.
\end{align*}
We now show that
\begin{equation}
\label{eq:bound-bias-0}
|\pscal{\nabla \lyap(\prm_k)}{\bias_k}|  \leq \bsf_0 + \bsf_1 \superlyap(\prm_k),
\end{equation}
Note first that, using \Cref{assum:field}-\ref{item:field:boundbias} and \eqref{eq:definition:clyap-lyap}, we get
\begin{align*}
|\pscal{\nabla \lyap(\prm_k)}{\bias_k}|
&\leq \clyap_{\lyap}\sqrt{\superlyap(\prm_k)} \{ \boundbias_{0}+ \boundbias_{1} \superlyap(\prm_k) \}^{1/2} , \\
&\leq \clyap_{\lyap} \{ \sqrt{\boundbias_{0}} \sqrt{\superlyap(\prm_k)} + \sqrt{\boundbias_{1}} \superlyap(\prm_k) \} \eqsp,
\end{align*}
and \eqref{eq:bound-bias-0} follows from $\sqrt{a} \leq (1/2) (1+a)$, $a \geq 0$. Lastly,
\begin{align*}
&\CPEs{\| \Hg(\prm_k,\State_{k+1}) \|^2}{\mcf_k}
= \| \hg(\prm_k) + \bias_k \|^2 \\
&+\CPE{\| \Hg(\prm_{k},\State_{k+1}) - \CPEs{\Hg(\prm_{k},\State_{k+1})}{\mcf_{k}}\|^2}{\mcf_{k}} \\
&\leq \eta_0 + \eta_1 \superlyap(\prm_k) \, ,
\end{align*}
where $\eta_0, \eta_1$ are defined in \eqref{eq:definition-eta}.
Combining the results above and \Cref{assum:lyapunov}-\ref{item:rholyap}, we get \eqref{eq:condition-lyap}.
\end{proof}

In \Cref{app:RS-forGD-general}, we provide a slightly different version for \Cref{lem:main-result-quantitative} for the particular case of GD. An important consequence of \Cref{lem:main-result-quantitative} is that it allows us to deduce a non-asymptotic bound on $\{ \superlyap( \prm_k ) \}_{k=0}^{\totstep-1}$ as follows.
\begin{tcolorbox}[boxsep=2pt,left=4pt,right=4pt,top=3pt,bottom=3pt]
\begin{theorem}
\label{theo:main-result-quantitative}
Assume \Cref{assum:field}, \ref{assum:lyapunov} and \Cref{assumNA:uniform-bound-bias}, \ref{assumNA:bias-is-small}.
Assume in addition that the step sizes $\sequence{\step}[k][\nset]$ are chosen such that $\step_k \in \ooint{0, \step_{\max}}$.
Then, for any $\totstep \geq 1$, we get
\begin{align}
\label{eq:main-result-non-asymptotic}
\nonumber
&\sum_{k=0}^{\totstep-1} \frac{ \step_{k+1} \omega_{k+1}}{\sum_{\ell=0}^{\totstep-1}  \step_{\ell+1} \omega_{\ell+1}} \PE[ \superlyap(\prm_k) ] \\
&\quad \leq  \frac{2(\PE[\lyap(\prm_0)] - \lyap_\star) +   \Liplyap \eta_0 \sum_{k=0}^{\totstep-1} \step^2_{k+1}  }{ \sum_{\ell=0}^{\totstep-1}  \step_{\ell+1} \omega_{\ell+1}}\\
\nonumber
&\quad + \frac{2 \bsf_0  \sum_{k=0}^{\totstep-1} \step_{k+1}}{\sum_{\ell=0}^{\totstep-1} \step_{\ell+1} \omega_{\ell+1}} \eqsp.
\end{align}
\end{theorem}
\end{tcolorbox}
\begin{proof}
Taking the expectations of both sides of \eqref{eq:condition-lyap} gives
\begin{multline*}
(1/2) \step_{k+1} \omega_{k+1} \PE[\superlyap(\prm_k)]
\le \PE[\lyap(\prm_k)] - \PE[\lyap(\prm_{k+1})] \\ + \step_{k+1} \bsf_0   + \step_{k+1}^2  \Liplyap \eta_0/2.
\end{multline*}
We obtain \eqref{eq:main-result-non-asymptotic} by summing these inequalities from $k=0$ to $k=T-1$ and by using \Cref{assum:lyapunov}-\ref{item:vstar}; note that under \Cref{assumNA:bias-is-small}, $\step_{\max}>0$ and $\step_{\ell+1} \omega_{\ell+1} >0$.
\end{proof}
Suppose that $\bsf_0 = 0$, under an appropriate stepsize policy (e.g., $\step_k = \step_{\max}/\sqrt{\totstep}$), it can be shown that the {\rhs} of \eqref{eq:main-result-non-asymptotic} goes to zero when $\totstep \to \infty$.
As discussed in \Cref{sec:assume},
for strict Lyapunov functions,
$\superlyap( \prm ) = 0$ implies that $\prm$ is an equilibrium point to the vector field $\hg(\prm) = {\bm 0}$.
We further recall that under \eqref{eq:condition-W}, controlling $\PE\left[\superlyap(\prm_R)\right]$ at some random stopping time $R$ leads to a high probability bound for the distance from $\prm_R$ to $\{ \hg = 0 \}$.

The following discussions demonstrate how to apply \eqref{eq:main-result-non-asymptotic} to obtain performance estimates  for the {\SA} scheme.



\vspace{.1cm}
\noindent \textbf{Random Stopping.}
We first discuss ways to implement the LHS in \eqref{eq:main-result-non-asymptotic}. In particular, we note that it
can be viewed as the expected value $\PE[\superlyap(\prm_{R_\totstep})]$, where $R_\totstep$
is a random variable taking values in $\{0, \cdots,\totstep-1\}$, independent of $\sequence{\prm}[k][\nset]$, and with probability mass function
\begin{equation} \label{eq:def:Rtotstep}
\PP ( R_\totstep = k)= \frac{\step_{k+1} \omega_{k+1}}{\sum_{\ell=0}^{\totstep-1}  \step_{\ell+1} \omega_{\ell+1}} \eqsp.
\end{equation}
We may regard $\prm_{R_\totstep}$ as the output of the {\SA} scheme terminated at the random iteration number $R_{\totstep}$.
Equivalently, one can look at such a randomization scheme from a slightly different perspective.
Here, one can also run the {\SA} algorithm for $\totstep$ iterations, but randomly choose a point $\prm_{R_\totstep}$ from its trajectory as the output of the algorithm. Clearly, in the latter method, the algorithm only needs to be run for the first $R_\totstep$ iterations. Note, however, that the primary goal of introducing the random iteration number $R_\totstep$ is to derive complexity results, not to save computational effort in the last $\totstep-R_\totstep$ iterations of the algorithm.
\begin{remark}
\label{rem:convex-superlyap}
When the function $\superlyap$ is convex, the LHS  of \eqref{eq:main-result-non-asymptotic} is lower bounded by
$\PE[\superlyap(\bar{\prm}_T)]$, where
\[
\bar{\prm}_T \eqdef  \sum_{k=0}^{\totstep-1} \frac{ \step_{k+1} \, \omega_{k+1}}{\sum_{\ell=0}^{\totstep-1}  \step_{\ell+1} \, \omega_{\ell+1}}  \prm_k \eqsp.
\]
In this case, one may adopt the convex combination $\bar{\prm}_T$ to achieve the LHS of  \eqref{eq:main-result-non-asymptotic} as an alternative to random stopping. \revisionupdates{Such averaging techniques are referred to as Polyak-Ruppert averaging, and were introduced in \cite{polyak1992acceleration,ruppert1988efficient}.}
\end{remark}

\vspace{.1cm}
\noindent \textbf{Constant step size.}
Let us analyze a special case when a constant stepsize policy is used, \ie $\step_{k+1}=\step$ for each $k \in \{0, \dots, \totstep -1\}$.
First, if $\step \leq \step_{\max}/2$, then \eqref{eq:main-result-non-asymptotic} can be written more explicitly as follows
\begin{equation}
\label{eq:main-result-non-asymptotic-constant}
\frac{1}{T} \sum_{k=0}^{\totstep-1} \PE [ \superlyap(\prm_k) ] \leq  \Bterm + \frac{2 \Vinit + \Liplyap \eta_0 T \step^2 }{ \step T  \{ \rholyap - \bsf_1\} } ,
\end{equation}
where
\begin{align}
\label{eq:definition-biais}
\Bterm  & \eqdef {2 \bsf_0 } / ({ \rholyap - \bsf_1}) \eqsp, \qquad \Vinit  \eqdef   \PE[\lyap(\prm_0)] - \lyap_\star.
\end{align}
Not surprisingly, the first term, $\Bterm$, in the  \rhs\ of \eqref{eq:main-result-non-asymptotic-constant} which is related to the bias of the random oracle $\Hg$, cannot be made small by any choice of $\step$. \revisionupdates{Indeed, according to its definition, $b_0$ (thus $\Bterm$) scales with $\tau_0$ (see \eqref{eq:definiton-bsf-0}), and $\tau_0$ corresponds (see  \Cref{assum:field}\ref{item:field:boundbias}) to  the part of the bias that \textit{does not} scale with the $\hg$ or $W$. When  $\tau_0\neq 0$, observing the oracle does not allow to find an exact 0 of the field $\hg$: indeed, the oracle could for example be equal to  $\hg+\tau_0$  and thus the \SA~scheme would converge to the roots of $\hg+\tau_0$, that are  distinct from those of $\hg$.}

On the other hand, the second term in the {\rhs} of \eqref{eq:main-result-non-asymptotic-constant}  mixes the dependence on the initial conditions, the bias and variance of the stochastic oracle. It can be adjusted according to the time horizon $\totstep$ and the different parameters of the problem by a suitable choice of $\step$.
The function  $\step \mapsto 2 \Vinit / \step + \Liplyap \eta_0 \totstep \step$ is minimized on $\ocint{0,\step_{\max}/2}$ by setting
\begin{equation}
\label{eq:optimal-step}
\step_{\totstep} \eqdef  (2 \Vinit / (\eta_0 \Liplyap \totstep) )^{1/2} \wedge (\gamma_{\max}/2).
\end{equation}
\begin{corollary}
\label{coro:main-result-quantitative}
Assume \Cref{assum:field,assum:lyapunov} and \Cref{assumNA:uniform-bound-bias,assumNA:bias-is-small}.
Then, for any  $\totstep \geq 1$, setting  $\step_{k+1} \eqdef \step_{\totstep}$ for $k \in \{0,\dots,\totstep-1\}$ we get
\begin{multline*}
\label{eq:conclusion-simplified}
\frac{1}{\totstep} \sum_{k=0}^{\totstep-1} \PE[ \superlyap(\prm_k) ]  \leq  \Bterm+ \frac{ 2 \, \sqrt{2 \Vinit \eta_0 \Liplyap}}{\sqrt{\totstep}  \{\rholyap - \bsf_1 \}}
\vee   \frac{8 \Vinit}{ \step_{\max} T  \{\rholyap - \bsf_1\}} \eqsp.
\end{multline*}
\end{corollary}
In the unbiased case ($\Bterm = 0$), this yields an ${\cal O}(1/\sqrt{T})$ convergence rate for {\SA}. 

\vspace{.1cm}
\noindent \textbf{Lower bounds and $\epsilon$-approximate stationarity.}
In general non-convex optimization, it is intractable to find a global minima of functions or even to test if a point is a local minimum or a high-order saddle point. As a remedy, the most common approach by far is to consider $\epsilon$-approximate stationarity. Our goal is to find a point $\prm_R \in \rset^d$ with
\begin{equation}\label{eq:def-eps-stat}
\PE[\superlyap(\prm_R) ] \leq \epsilon \eqsp,
\end{equation}
where the expectation is taken over the randomness in both the mean-field oracle and the query $R$.
The use of stationarity as a convergence criterion dates back to the early days of nonlinear optimization (see \cite{wright1999numerical}). Recent years have seen rapid development of a body of work that studies non-convex optimization through the lens of non-asymptotic convergence rates to $\epsilon$-stationary points \cite{ghadimi2013stochastic,carmon2017convex,fang2018spider,lei2017non}.

We will discuss how to choose the constant step size $\step$ and the total number of iterations $\totstep$ to guarantee $\epsilon$-approximate stationarity when USO is satisfied. 
\Cref{coro:main-result-quantitative} implies
\begin{corollary}\label{cor:complexity}
Assume \Cref{assum:field}, \Cref{assum:lyapunov}-\Cref{assumNA:uniform-bound-bias}, \Cref{assumNA:bias-is-small} and USO (i.e., $\boundbias_0=\boundbias_1=0$).
Then, for $\epsilon>0$, the number of iterations to guarantee an $\epsilon$ approximate stationary point \eqref{eq:def-eps-stat} is lower bounded by
\begin{equation}\label{eq:totstep}
    \totstep(\epsilon) =   \frac{ 8 \Vinit \eta_0 \Liplyap}{\epsilon^2 \rholyap^2} \vee   \frac{8 \Vinit}{ \step_{\max} \epsilon \rholyap } \eqsp .
\end{equation}
\end{corollary}
\begin{proof}
This complexity is obtained by choosing the query $R$ to be uniform over $\{0,\dots,\totstep-1\}$, and  $\totstep$ such that the upper bound in \Cref{coro:main-result-quantitative} is at most $\epsilon$. Under USO, $\bsf_1 =0$ and this bound reads  $\PE[ \superlyap(\prm_R) ]  \leq  \nicefrac{ 2\sqrt{2 \Vinit \eta_0 \Liplyap}}{(\sqrt{\totstep} \rholyap)}
\vee   \nicefrac{8 \Vinit}{(\step_{\max} \totstep \rholyap )}$. It is easily checked that \eqref{eq:totstep} ensures \eqref{eq:def-eps-stat}.
\end{proof}
When $\totstep \geq \totstep(\epsilon)$, the algorithm which returns $\prm_R$ when $R$ is a uniform random variable on $\{0, \dots, \totstep-1\}$, satisfies the $\epsilon$-approximate stationarity condition $\PE\left[ \superlyap(\prm_R)\right] \leq \epsilon$.
The upper bound  in~\Cref{cor:complexity} shows that  there are two regimes depending on the value of $\epsilon$ \wrt $\nofrac{\step_{\max}   \eta_0 \Liplyap}{\rholyap} =2 \eta_0/\eta_1$.
\begin{enumerate}[label=\alph*),nosep]
\item  \label{item:complexity-low} In the high-precision regime where $\epsilon \in \ocint{0,2\eta_0/\eta_1}$, 
\begin{equation}
\label{eq:expression-Tepsilon-HP}
\totstep(\epsilon) = \frac{ 8 \Vinit \Liplyap }{ \rholyap^2} \; \frac{ \eta_0}{\epsilon^2},
\end{equation}
achieved with a constant step size $\step(\epsilon)= \frac{\rholyap \epsilon}{2 \eta_0 \Liplyap}$.
\item  \label{item:complexity-high} In the low-precision regime where $\epsilon  > \nicefrac{2 \eta_0}{\eta_1} $,
\begin{equation}
\label{eq:expression-Tepsilon-LP}
\totstep(\epsilon) = \frac{4  \Vinit \Liplyap }{\rholyap^2} \frac{\eta_1}{\epsilon},
\end{equation}
achieved with a constant stepsize $\step= \step_{\max}/2$.
\end{enumerate}

\vspace{.1cm}
\noindent \textbf{Last iterate or ``Random'' iterate?}
The standard analysis holds only for random stopping, or, when $\superlyap$ is convex, for a convex combination of the iterates.
Most practitioners just use the final iterate of {\SA} instead of randomly selecting a solution $\prm_{R_\totstep}$ from
$\{\prm_k\}_{k=0}^{\totstep-1}$. In the case of {\SG} for convex functions (the mean field is the gradient of a smooth convex function) and of a stochastic oracle without bias, \cite{jain2019making,jain2021making} (see also \cite{shamir2013stochastic}) show that with a clever choice of piecewise constant
step, it is possible to obtain for the last iteration the same decay $1/\sqrt{T}$ as for the randomly stopped estimator (or the averaged one in the case where $\superlyap$ is convex, see \Cref{rem:convex-superlyap}). Unfortunately, this approach is strongly linked to the use of a gradient algorithm and does not extend to {\SA} under the general assumptions we consider.

Another option would be to output a solution $\hat{\prm}_\totstep$ such that
\begin{equation} \label{eq:SA-end-bestiterate}
 \superlyap(\hat{\prm}_\totstep)=\min_{k=0, \ldots, \totstep-1} \superlyap(\prm_k)\eqsp.
\end{equation}
In the case of {\SG} for smooth functions and unbiased oracles, the  lower bounds reported for example in \cite{arjevani2022lower}  show that there is no hope of improving the speed in $1/\sqrt{T}$.
Using \eqref{eq:SA-end-bestiterate} requires additional computational effort for $\{ \superlyap(\prm_k) \}_{k=0}^{\totstep-1}$. Since in most practical cases $ \superlyap(\prm_k)$ cannot be computed exactly, estimation using Monte Carlo simulations would introduce approximation errors and raise robustness issues. For {\SG}, \cite[Section~6.1.1.2]{lan2020first} describes a two-stages procedure which runs $S$ times the optimization procedure from the same initial condition. This procedure provides a comparable theoretical guarantee, but is rarely used in practice.

\vspace{.1cm}
\noindent \textbf{Faster Rate.}
It is possible to improve \Cref{theo:main-result-quantitative} if \Cref{assum:field,assum:lyapunov} hold with $\superlyap =  \lyap$; note that, since $\superlyap \geq 0$ and is null on $\Lambda_\lyap$, this implies $\lyap_\star = 0$. We note that this setting applies to the stochastic gradient algorithms with strongly convex objective function and a special case of TD(0) learning. 
Starting from \Cref{lem:main-result-quantitative}, taking the expectation of both sides of \eqref{eq:condition-lyap}, we get
\begin{equation}
\label{eq:condition-lyap-1}
\PE[\superlyap(\prm_{k+1})]  \leq   \lambda_{k+1} \,  \PE[\superlyap(\prm_k)] + b_{k+1}
\end{equation}
where  
\begin{align}
\label{eq:definition-lambda-k}
\lambda_{k}& \eqdef   1 - \step_{k} ( \rholyap -  \bsf_1 ) + \step^2_{k} \Liplyap \eta_1/2,  \\
\label{eq:definition-b-k}
b_{k} & \eqdef  \step_{k} \bsf_0 +  \step_{k}^2 \Liplyap \eta_0/2 \,.
\end{align}
A straightforward induction shows that for any $k \in \nset$,
\begin{equation}
\label{eq:expansion}
\textstyle \PE[\superlyap(\prm_{k})] \leq \Lambda_{1:k} \, \PE[\superlyap(\prm_0)] + \sum_{j=1}^{k} \Lambda_{j+1:k} \,  b_{j},
\end{equation}
where we have set $\Lambda_{k+1:k}\eqdef 1$ and for $1 \leq j \leq k$,
$\Lambda_{j:k} \eqdef  \prod_{i=j}^k \lambda_i$.
Plugging \eqref{eq:definition-b-k} into \eqref{eq:expansion}, we obtain for any $k \in \nset$,
\begin{multline*}
\PE[\superlyap(\prm_{k+1})]
\leq  \Lambda_{1:k+1} \PE[\superlyap(\prm_0)]  \\
\textstyle + \frac{ \Liplyap \eta_0}{2} \sum_{j=1}^{k+1} \step^2_{j} \Lambda_{j+1:k+1}
+ \bsf_0 \sum_{j=1}^{k+1} \step_{j} \Lambda_{j+1:k+1} \eqsp.
\end{multline*}
\Cref{lem:summation} provides sharp estimates of $\sum_{j=1}^{k+1} \step^{\ell}_{j} \Lambda_{j+1:k+1}$, for $\ell\in\{1,2\}$, which leads to the following conclusions.
\begin{tcolorbox}[boxsep=2pt,left=4pt,right=4pt,top=3pt,bottom=3pt]
\begin{theorem}
\label{theo:main-result-quantitative-fast}
Assume \Cref{assum:field}, \ref{assum:lyapunov} are satisfied with $\lyap=\superlyap$.
 Assume in addition \Cref{assumNA:uniform-bound-bias}, \ref{assumNA:bias-is-small}, and  that the stepsize sequence $\sequence{\step}[k][\nset]$ is non-increasing and chosen such that, for any $k \in \nset$, $\step_{k} \leq \step_{\max}/2$ and
\begin{equation}
\label{eq:condition-step}
{\step_{k}} / {\step_{k+1}} \leq 1 + \gamma_{k+1} ({\rholyap-\bsf_1}) / 4 \eqsp.
\end{equation}
Then, for any  $k \geq 0$, we get
\begin{equation}
\PE[ \superlyap(\prm_{k}) ] \leq \Lambda_{1:k} \,  \PE[\superlyap(\prm_0)]
 + \frac{2 \Liplyap \eta_0 }{\rholyap-\bsf_1} \step_{k} + \Bterm\eqsp.
\end{equation}
\end{theorem}
\end{tcolorbox}
\begin{proof}
Under \eqref{eq:condition-step} and $\step_{k+1} \leq \step_{\max}/2$, we have $\lambda_{k+1} \leq 1 - \step_{k+1} ( \rholyap - \bsf_1)/2$. Using \Cref{lem:summation} in the Appendix, with $a \eqdef (\rholyap - \bsf_1)/2$ and $b \eqdef (\rholyap-\bsf_1)/4$, concludes the proof.
\end{proof}
Condition \eqref{eq:condition-step} encompasses constant stepsize policies which are common
in the literature. Diminishing stepsize  sequences  are also a popular choice, e.g., the sequence
\begin{equation}
\label{eq:diminishing-step}
\step_{k+1} \eqdef \tilde{\step} /( k +1 + \totstep_0)^{\beta} \quad \text{for $\beta \in \ocint{0,1}$},
\end{equation}
satisfies \eqref{eq:condition-step}  by choosing appropriately  $\tilde{\step}$ and $\totstep_0$.
The following corollary is obtained by setting  $\beta=1$.
\begin{corollary}
\label{coro:main-result-quantitative-fast}
Assume \Cref{assum:field,assum:lyapunov} are satisfied with $\lyap=\superlyap$.  Assume in addition \Cref{assumNA:uniform-bound-bias,assumNA:bias-is-small}.  Let $\tilde{\step} \geq 6/( \rholyap - \bsf_1)$ and $\totstep_0 \geq 2 \tilde{\gamma}/\gamma_{\max}$; and
suppose that  $\step_{k+1} \eqdef \tilde{\step}/(k+1+\totstep_0)$, for $k \in \{0, \dots, \totstep-1 \}$. Then  it holds
\begin{align}
& \PE[ \superlyap(\prm_\totstep) ] \leq \nonumber \\
& \Bterm + \left( \frac{ \totstep_0 }{ \totstep + \totstep_0 } \right)^{ \frac{\tilde{\gamma} (\rholyap - \bsf_1)}{ 2 } } \PE\left[ \superlyap(\prm_0) \right] + \frac{2 \Liplyap \eta_0 \tilde{\gamma} }{(\totstep+\totstep_0)(\rholyap-\bsf_1)}
\eqsp. \label{eq:conclusion-simplified-bis}
\end{align}
\end{corollary}
\Cref{coro:main-result-quantitative-fast} provides a bound  on the \emph{last iterate} of {\SA}. In the upper bound, The second term is  ${\cal O}(\totstep^{-\alpha})$  for some $\alpha \geq 3$, and the third term is  ${\cal O}(\totstep^{-1})$.  As discussed in \cite{nemirovski2009robust,moulines2011non},  for {\SG} this upper bound is tight up to factors independent of $\totstep$.

\subsection{\revisionupdates{Application to the examples}}
\label{subsec:finite-time-bounds-examples}
\subsubsection{Stochastic Gradient Method}
\label{subsec:SGD}
We apply the results of the previous sections to the case of {\SG}. As there are many results in the literature on this subject, \eg \cite{nemirovski2009robust,juditsky2011first,juditsky2011first-ii,bubeck2015convex,lan2020first} and the references therein which offer excellent discussions, the only interest in the following results is to show that we can find ``state-of-the-art'' results by simply choosing  $\lyap$ and $\superlyap$.
\begin{proposition}
Assume \Cref{assum:SGD}, the sequence $\sequence{\step}[k][\nset^*]$ satisfies $0 < \step_k \leq 2/ \Lip{\nabla F}$. The following holds for {\SG} \eqref{eq:definition-SGD},
\begin{enumerate}[leftmargin=*]
\item \label{item:result:SGD} For any $T \in \nset$,
\begin{multline}
\label{eq:main-result-non-asymptotic-SGD}
\sum_{k=0}^{\totstep-1} \frac{( 2 \step_{k+1} - \Lip{\nabla F}\step^2_{k+1} )}{\sum_{\ell=0}^{\totstep-1}  ( 2 \step_{\ell+1} - \Lip{\nabla F} \step^2_{\ell+1}  )} \PE[ \| \nabla F(\prm_k) \|^2 ] \\
\leq \frac{\PE[F(\prm_0)] - F_\star  +   \Lip{\nabla F} (M^2/n) \sum_{k=0}^{\totstep-1} \step^2_{k+1}  }{\sum_{\ell=0}^{\totstep-1}  ( 2 \step_{\ell+1} - \Lip{\nabla F} \step^2_{\ell+1})}.
\end{multline}
\item \label{item:result:SGD:convex} If \Cref{assum:SGD:convex} also holds, then for any $T \in \nset$,
\begin{equation}
\PE[ F(\bar{\prm}_{\totstep})] - F_\star
 \leq \frac{\PE[\|\prm_0 - \prm_\star\|^2]  + (M^2/n) \sum_{k=0}^{\totstep-1} \step^2_{k+1}  }{\sum_{\ell=0}^{\totstep-1}  ( 2 \step_{\ell+1}  - \Lip{\nabla F} \step^2_{\ell+1}   )}
\end{equation}
where $\bar{\prm}_{\totstep}
= \sum_{k=0}^{\totstep-1} \frac{( 2 \step_{k+1}  - \Lip{\nabla F}\step^2_{k+1} \eta  )}{\sum_{\ell=0}^{\totstep-1}  ( 2 \step_{\ell+1} - \Lip{\nabla F} \step^2_{\ell+1}  )} \prm_k$.
\item \label{item:result:SGD:strongly-convex} If \Cref{assum:SGD:strongly-convex} also holds and  $\step_k < \mu / \Lip{\nabla F}^2 $, $\nofrac{\step_k}{\step_{k+1}} \leq 1 + (\nofrac{\mu}{4}) \step_{k+1}$.
Then, for any $k \in \nset$,
\begin{multline}
\PE[\| \prm_k - \prm_\star\|^2] \leq \frac{\Lip{\nabla  F}^2 (M^2/n)}{\mu} \step_k^2 \\
+ \prod_{l=1}^k \left(1-2 \mu \step_l+2 \Lip{\nabla F}^2 \step_l^2\right) \PE[\| \prm_0 - \prm_\star \|^2] \eqsp.
\end{multline}
\end{enumerate}
\end{proposition}
Notice that item \ref{item:result:SGD} retrieves the conclusions of \cite[Theorem~6.1]{lan2020first}; item \ref{item:result:SGD:convex} is classical where precise results are given in \cite[Section~2] {nemirovski2009robust} but these results were known since \cite{nemirovskij1983problem}, see \cite{shalev2009stochastic,bubeck2015convex}; item \ref{item:result:SGD:strongly-convex} has a long history, see \cite[Section~2]{nemirovski2009robust}.

\subsubsection{Compressed {\SA}}
\label{subsub:exrate:CSGD} 
We describe convergence results obtained for compressed {\SA}. We focus on the case of a stochastic gradient field \eqref{eq:SGD-finite-sum}, which has been the most studied. The proofs for this section are given in \Cref{app:proofGS}.
We first analyze the Gauss-Southwell update~\eqref{eq:gauss_south}, which has been studied~\cite{nutini2015coordinate}, and is also known as the steepest-coordinate descent~\cite{boyd2004convex}.
Our results can be naturally extended to the compressed algorithm GD for any compressor that satisfies \Cref{assum:SGD:contractive}. 
\begin{proposition}\label{prop:rateGS}
Assume \Cref{assum:SGD} and consider a constant sequence $\sequence{\step}[k][\nset^*]$ chosen as $ \step_k  = \nicefrac{1}{8 d \Lip{ \nabla F } }$. The following holds for the iterates of the Gauss-Southwell algorithm~\eqref{eq:gauss_south}.
\begin{align}
\nonumber
&\frac{1}{T} \sum_{k=0}^{\totstep-1} \PE[ \| \nabla F(\prm_k) \|^2 ]  \leq  \frac{32 d^2 \Lip{ \nabla F }(\PE[F(\prm_0)] - F_\star)   }{ T } \eqsp. 
\end{align}
\end{proposition}
This result can be compared with \cite[Eq.~(7)]{nutini2015coordinate} for Gauss-Southwell, \cite[Theorem 5.3]{deklerk2020} for GD with relatively bounded errors, both of which give a result under an additional assumption of strong convexity. A similar result, also given in \cite[Theorem J.3]{gorbunov2020linearly} for an error-compensated version of the compressed algorithm, shows an increase in rate by a factor $\oneminusomgbiased^{-1} $. Note that our result has a dependence on $d^2$ that is suboptimal. A refined version of~\Cref{lem:main-result-quantitative} is given in \Cref{app:RS-forGD-general}, which allows us to obtain a bound scaling with $d$. However, the corresponding \Cref{lem:main-result-quantitative_2} only allows us to improve convergence for (i) gradient methods (stochastic or not), (ii) with bias, and (iii) when we choose $\lyap = F$. Since our interest is much more general, we omit the refinements of this more detailed analysis.

Second, we consider the case of unbiased compression as in~\eqref{eq:cSGD} applied to~a stochastic gradient descent update. This case has been extensively studied in the communication-constrained distributed optimization community, starting with \cite{alistarh_qsgd_2017} and with several successors \cite{horvoth2022natural}. Combining~\Cref{lem:csgd-var2} and \Cref{theo:main-result-quantitative} gives the following result:
\begin{proposition}\label{prop:compSGDunbiased}
Assume \Cref{assum:SGD} and consider unbiased-compressed \SG, i.e. \Cref{eq:cSGD} with a $\compressor$ satisfying~\Cref{assum:SGD:URBV} and the sequence $\sequence{\step}[k][\nset^*]$ satisfying $0 < \step_k \le  1 /(\Liplyap (2\omgunbiased+1))  $. The following holds for any $T \in \nset$,
\begin{align}
\label{eq:main-result-non-asymptotic-SGD}
& \sum_{k=0}^{\totstep-1} \frac{( 2 \step_{k+1} - \Lip{\nabla F}\step^2_{k+1} )}{\sum_{\ell=0}^{\totstep-1}  ( 2 \step_{\ell+1} - \Lip{\nabla F} \step^2_{\ell+1}  )} \PE[ \| \nabla F(\prm_k) \|^2 ] \\
& \leq \frac{\PE[F(\prm_0)] - F_\star  +  (M^2/n) \Lip{\nabla F} (1+\omgunbiased)  \sum_{k=0}^{\totstep-1} \step^2_{k+1}  }{\sum_{\ell=0}^{\totstep-1}  ( 2 \step_{\ell+1} - \Lip{\nabla F} \step^2_{\ell+1})}. \nonumber
\end{align}
\end{proposition}
The unbiased compression leads to a multiplicative increase of the noise level and a reduction of the maximum stepsize by a factor $1+2\omgunbiased$; see \cite[Theorem~9]{horvoth2022natural}. Third, we apply \Cref{theo:main-result-quantitative} to STE using \eqref{lem:csgd-var3} with deterministic rounding \eqref{eq:deterministic-rounding}.
\begin{proposition}\label{prop:STE}
 Assume \Cref{assum:SGD}
and consider the STE compression {\SG} algorithm~\eqref{eq:binary-connect}. 
Assume $\compressor$ is $Q_d$ in \eqref{eq:deterministic-rounding} with a quantization step $\Delta < 1/ (\Lip{\nabla F} \sqrt{d})$. Assume that the sequence $\sequence{\step}[k][\nset^*]$  is constant $\sequence{\step}[k]=\step/\sqrt{\totstep}  $, with $0 < \step \leq 1/ \Lip{\nabla F}$. 
The following holds for the iterates:
\begin{align*}
\nonumber
 &\frac{1}{\totstep}\sum_{k=0}^{\totstep-1} \PE[ \| \nabla F(\prm_k) \|^2]\leq 
  {2 \Lip{{\nabla F}} \sqrt{d} \Delta  } 
 \\
 & \qquad+\frac{2(\PE[\lyap(\prm_0)] - \lyap_\star) +    ( M^2/n + 3  \Lip{{\nabla F}} \sqrt{d} \Delta) \step^2 }{\sqrt{\totstep}   \step /2} 
\end{align*}
\end{proposition}
This rate can be compared to \cite[Theorem~3]{li2017training} for \texttt{BinaryConnect}, although we do not assume a bounded domain: we have a convergence rate of $1/\sqrt{\totstep}$ up to a fixed threshold proportional to $\sqrt{d}\Delta$.
Finally, we consider the case of low precision {\SG } in terms of \eqref{eq:de-sa}. 
Combining \Cref{lem:csgd-var4} with \eqref{eq:main-result-non-asymptotic-constant} (which holds under conditions of \Cref{theo:main-result-quantitative}) gives:
\begin{proposition}\label{prop:compSGDLinear}
Assume \Cref{assum:SGD} and consider compressed {\SG} \eqref{eq:de-sa} with constant stepsize $\step_k = \bar{\step}$. If $\compressor$ satisfies~\Cref{assum:SGD:UniformQuant} and
\begin{equation} \label{eq:main-result-non-asymptotic-LinearQ-step}
\bar{\step} + \linearunbiased \sqrt{d} < { 2 } / { \Lip{\nabla F} }.
\end{equation}
Then the following holds for any $T \in \nset$,
\begin{align}
\frac{1}{\totstep} \sum_{k=0}^{\totstep-1} \PE[ \| \nabla F (\prm_k) \|^2 ] & \leq \frac{2 (\PE[F(\prm_0)] - F_\star)}{ \bar{\step} T } + \bar{\step} \,  \frac{M^2\Lip{\nabla F} }{n} \nonumber \\
& + \sqrt{d} \linearunbiased \Lip{\nabla F} \left( 3 + {M^2}/ {n}\right). \label{eq:main-result-non-asymptotic-LinearQ}
\end{align}
\end{proposition}
Eq.~\eqref{eq:main-result-non-asymptotic-LinearQ-step} holds provided that $\linearunbiased \leq 2 / \Lip{\nabla F} \sqrt{d}$. In the case of random quantization, this constraints the resolution of quantizer. Moreover, inserting $\bar{\step} = {\cal O}( 1/ \sqrt{T} )$ shows that the first two terms in \eqref{eq:main-result-non-asymptotic-LinearQ} are ${\cal O}( 1/ \sqrt{T} )$, while the last term does not vanish with $T$. This shows that the compressed {\SG} Algorithm converges to a ${\cal O}( \linearunbiased \sqrt{d} )$ approximate stationary solution, i.e., similar to \cite[Theorem 2]{Zheng19}.

\subsubsection{Stochastic EM algorithms}
\label{subsec:stochastic-EM}
The stochastic EM algorithms are used to find a root of the mean-field $\hg$ defined in \eqref{eq:EM:meanfield}. We specialize our results to the Mini-batch EM introduced in \Cref{sec:EM:verifH:onlineEM} and to the SAEM algorithm introduced in \Cref{sec:EM:verifH:SAEM}, under the assumptions \Cref{assumEM:expo,assumEM:Mstep,assumEM:supplementary}.

\paragraph{Mini-batch EM}  This case is an example of USO (see \Cref{def:USO}) where
$\boundbias_0=\boundbias_1 =0$. Using \Cref{sec:EM:verifH}, it is easily checked that $\eta_0  =  \nicefrac{\boundvarem_0}{\lbatchem}$,  $\eta_1  = 1+ \nicefrac{\boundvarem_1}{\lbatchem}$, $\step_{\max} =  \nicefrac{2 v_{\min}}{ (\Liplyap \eta_1) }$, and $  \omega_k  = 2 v_{\min} - \step_k \Liplyap \eta_1$.
Let $\totstep > 0$ be the total number of iterations. Assume that, for all $k \in \{0,\dots,\totstep-1\}$,   $\step_{k+1} < \nicefrac{2 v_{\min}}{ (\Liplyap  \eta_1)}$. If the Mini-batch EM is stopped at a random iteration index $R_{\totstep}$ whose distribution is given in \eqref{eq:def:Rtotstep},  then $\PE[\|\hg(\prm_{R_\totstep})\|^2]$  is upper bounded  by \Cref{theo:main-result-quantitative} as:
\begin{proposition}\label{prop:batchEM:quantitativebounds}
Assume \Cref{assumEM:expo,assumEM:Mstep,assumEM:supplementary,assumEM:var:onlineEM}. Let $\sequence{\prm}[k][\nset]$ be the {\SA} sequence with random oracle \eqref{eq:onlineEM:oracle}.  Then for any $\totstep \in \nset$,
\begin{align*}
\PE\left[ \|\hg(\prm_{R_\totstep})\|^2\right] \leq \frac{2 \Vinit +   \Liplyap \boundvarem_0 \sum_{k=1}^{\totstep} \step^2_{k} / \lbatchem }{ \sum_{k=1}^{\totstep}  \step_{k}  (2 v_{\min} - \step_{k} \Liplyap (1+\boundvarem_1/ \lbatchem))},
\end{align*}
where $R_\totstep$ is the uniform random variable on $\{0, \ldots, T-1\}$ and $\Vinit$, $\Liplyap$, $\vminem$ and $\boundvarem_\ell$ are defined by \eqref{eq:definition-biais}, \Cref{assumEM:supplementary} and \Cref{assumEM:var:onlineEM}.
\end{proposition}
Assume constant step sizes on $\{0,\dots,\totstep-1\}$:  
\[
\step_{k+1} = \step \eqdef  \left(\frac{2 \Vinit}{\boundvarem_0 \Liplyap} \frac{\lbatchem}{\totstep} \right)^{1/2} \wedge \frac{v_{\min}}{\Liplyap (1+\boundvarem_1/ \lbatchem)},
\]
and that $R_{\totstep}$ is the uniform random variable on $\{0, \ldots, T-1\}$.
Then, using \Cref{coro:main-result-quantitative},  we get the following  upper bound
\begin{equation}
\label{eq:complexity-EM}
\PE\left[ \|\hg(\prm_{R_\totstep})\|^2\right] \leq  \frac{ 2 \sqrt{2 \Vinit \boundvarem_0 \Liplyap}}{\sqrt{\lbatchem \, \totstep} v_{\min}}
\vee   \frac{8 \Vinit \Liplyap (1+\boundvarem_1/ \lbatchem)}{ 2 T  v_{\min}^2}.
\end{equation}

Let $\epsilon \in \ooint{0, 2 \boundvarem_0/\boundvarem_1}$. We discuss how to select the number of iterations $\totstep$, the step size $\step$  and the size of the mini-batch $\lbatchem$ so that $\PE\left[ \|\hg(\prm_{R_\totstep})\|^2\right]  \leq \epsilon$.
Assume first that 
\begin{equation}
\label{eq:constraint-epsilon-EM}
0 < \epsilon \leq {2 \boundvarem_0} / ({\lbatchem +  \boundvarem_1}) \eqsp.
\end{equation}
In this case (see \eqref{eq:expression-Tepsilon-HP}), the number of iterations needed to guarantee an $\epsilon$-approximate stationary point is
\begin{equation}
\label{eq:bound-complexity-EM-HP}
\totstep(\epsilon, \lbatchem) \eqdef  \frac{ 8 \Vinit \Liplyap }{ v_{\min}^2} \; \frac{\boundvarem_0}{\lbatchem\epsilon^2} \eqsp,
\end{equation}
which grows as $\epsilon^{-2}$ as $\epsilon \downarrow 0^+$.  The  bound in \eqref{eq:bound-complexity-EM-HP} is achieved by taking a constant step size $\step(\epsilon,\lbatchem) \eqdef \nicefrac{v_{\min}  \lbatchem \epsilon}{(2\boundvarem_0 \Liplyap)}$. We observe that $\totstep(\epsilon,\lbatchem)$ is inversely proportional to the mini-batch size $\lbatchem$, while the step size is proportional to $\lbatchem$. Increasing $\lbatchem$ allows more aggressive step sizes to be used and the number of iterations to be reduced accordingly. 

It is interesting to study the impact of the choice of $\lbatchem$ on the computational complexity. Note that the computational cost of mini-batch EM depends on two factors: the evaluation of the functions $\barsem_i$ for the current mini-batch $\State_{k+1}$, and the cost of updating the parameter by calling the optimization map $\mapem$. The cost of evaluating the stochastic oracle is $\lbatchem \cost_{\barsem} + \cost_{\mapem}$, and  after $\totstep(\epsilon,\lbatchem)$ iterations, this cost is
\begin{equation} \label{eq:EM:complexity:case1}
\cost(\epsilon, \lbatchem) \eqdef  \frac{ 8 \Vinit \Liplyap \boundvarem_0}{ v_{\min}^2 \epsilon^2}\cost_{\barsem}   \, \left( 1 +   \frac{\cost_{\mapem}}{\cost_{\barsem}} \frac{1}{ \lbatchem} \right) \eqsp.
\end{equation}
If $\cost_{\mapem}$ is negligible \wrt $\cost_{\barsem}$, there is no clear incentive to take a mini-batch larger than $1$  (to see the interest in using a larger mini-batch, we would have to go much further in evaluating computational cost by considering the possibility of parallelization or multithreading, etc.). However, if $\cost_{\mapem}$ is taken into account - which is often the case, since the evaluation of $\mapem$ requires solving an optimization program - then it becomes interesting to increase the size of the mini-batch.
Since this discussion assumes that \eqref{eq:constraint-epsilon-EM} holds, the maximum batch size is (up to appropriate roundings):
\[
\lbatchem(\epsilon)  \eqdef \nofrac{2 \boundvarem_0}{\epsilon} - \boundvarem_1  \eqsp.
\]
This ``optimal'' mini-batch size is inversely proportional to the  accuracy $\epsilon$, consistent with the fact that using aggressive strategies to reduce the number of iterations is a win.  It is interesting to note that the step
size $\step(\epsilon,\lbatchem(\epsilon))$ is equal to $\vminem(1- \epsilon \boundvarem_1/(2 \boundvarem_0))/\Liplyap > \vminem / (2 \Liplyap)$.

Assume now that $\lbatchem$ is such that $\epsilon \geq \nofrac{2 \boundvarem_0}{(\lbatchem + \boundvarem_1)}$.
Using \eqref{eq:expression-Tepsilon-LP}, the total number of iterations needed to guarantee an $\epsilon$-approximate stationary point is 
\begin{equation}
\label{eq:bound-complexity-EM-LP}
\totstep(\epsilon, \lbatchem) \eqdef  \frac{4  \Vinit \Liplyap }{\vminem^2} \frac{1+ \nofrac{\boundvarem_1}{\lbatchem}}{ \epsilon  } \,,
\end{equation}
which grows as $\epsilon^{-1}$ as $\epsilon \downarrow 0^+$.  The  bound in \eqref{eq:bound-complexity-EM-LP} is achieved by taking a constant step size $\step(\epsilon,\lbatchem) \eqdef \step_{\max}/2 = \vminem  \lbatchem/ ( \Liplyap (\lbatchem+\boundvarem_1))$. After $\totstep(\epsilon,\lbatchem)$ iterations, the total cost of evaluating the stochastic oracles is
\[
\cost(\epsilon, \lbatchem) \eqdef  \frac{4  \Vinit \Liplyap }{\vminem^2 \epsilon}  (\lbatchem+\boundvarem_1) \cost_{\barsem}   \left( 1+  \frac{  \cost_{\mapem}}{\cost_{\barsem} }  \frac{1}{\lbatchem}\right).
\]
This quantity is minimized by $\lbatchem(\epsilon) \eqdef (2 \boundvarem_0 /\epsilon - \boundvarem_1) \vee  \lbatchem^\star$,
\[
\lbatchem^\star\eqdef \sqrt{\boundvarem_1 \cost_{\mapem}/\cost_{\barsem}}.
\]
 When $\epsilon \leq 2 \boundvarem_0/ (\lbatchem^\star+\boundvarem_1)$, then $\lbatchem(\epsilon) = 2 \boundvarem_0 /\epsilon - \boundvarem_1$ and $\cost(\epsilon, \lbatchem(\epsilon))$ is equal to the previous case (see \eqref{eq:EM:complexity:case1}); otherwise, the cost  $\cost(\epsilon, \lbatchem^\star)$ is lower. See the summary of the cost for minibatch EM in Table~\ref{tab:EMtable}.

\paragraph{SAEM with exact sampling (SAEM-ES)}
Consider the case when conditionally to the past, the random variables $\{Z_{i,k+1}^j, 1 \leq j \leq \nbrMC\}$ are  sampled from the distribution $\pi_i(z_i; \mapem(\prm_k))$ for all $i \in \{1,\ldots, n\}$; and  $\{Z_{i,k+1}^j, 1  \leq j \leq \nbrMC, i \in \{1, \ldots, n\} \}$ are independent.  Assume  \Cref{assumEM:var:SAEM:exact}.

We have (see \Cref{sec:EM:verifH:SAEM}) $\eta_0  =  \nicefrac{\boundvarem_0}{(n \nbrMC)}$, $\eta_1  = 1+\nicefrac{\boundvarem_1}{(n \nbrMC)}$, $\bsf_0  = \bsf_1 =0$,
$\step_{\max}  = 2 \nicefrac{v_{\min}}{( \Liplyap \eta_1)}$ and $\omega_k  = 2 v_{\min} - \step_k \Liplyap \eta_1$. From \Cref{theo:main-result-quantitative}, we obtain the following result.
\begin{proposition}\label{prop:SAEMiid:quantitativebounds}
Assume \Cref{assumEM:expo,assumEM:Mstep,assumEM:supplementary}. Let $\sequence{\prm}[k][\nset]$ be the {\SA} sequence with random oracle \eqref{eq:oracle-SAEM}.  Assume in addition \Cref{assumEM:var:SAEM:exact}. Then for any $\totstep \in \nset$,
\begin{align*}
\PE\left[ \|\hg(\prm_{R_\totstep})\|^2\right] \leq \frac{2 \Vinit +   \Liplyap \boundvarem_0 \sum_{k=1}^{\totstep} \step^2_{k} / (n \nbrMC) }{ \sum_{k=1}^{\totstep}  \step_{k}  (2 v_{\min} - \step_{k} \Liplyap (1+\boundvarem_1/ (n \nbrMC)))},
\end{align*}
where $R_\totstep$ is the uniform random variable on $\{0, \ldots, T-1\}$ and $\Vinit$, $\Liplyap$, $\vminem$ and $\boundvarem_\ell$ are defined by \eqref{eq:definition-biais}, \Cref{assumEM:supplementary} and \Cref{assumEM:var:SAEM:exact}.
\end{proposition}

We discuss the case of constant step sizes $\step_k = \step$.  By \Cref{coro:main-result-quantitative},  we choose
\[
\step \eqdef  \left(\frac{2 \Vinit}{\boundvarem_0 \Liplyap} \frac{n \nbrMC}{\totstep} \right)^{1/2} \wedge \frac{v_{\min}}{\Liplyap (1+\boundvarem_1/ (n \nbrMC))},
\]
and obtain the upper bound
\[
\PE\left[ \|\hg(\prm_{R_\totstep})\|^2\right] \leq  \frac{ 2 \sqrt{2 \Vinit \boundvarem_0 \Liplyap}}{\sqrt{n \nbrMC \, \totstep} v_{\min}}
\vee   \frac{8 \Vinit \Liplyap (1+\boundvarem_1/ (n \nbrMC))}{ 2 T  v_{\min}^2},
\]
where $R_{\totstep}$ is the uniform random variable on $\{0, \ldots, T-1\}$.\vspace{.1cm}

Let $\epsilon \in \ooint{0, 2 \boundvarem_0/\boundvarem_1}$. We discuss how to select the number of iterations $\totstep$, the step size $\step$  and the number $\nbrMC$ of  Monte Carlo samples  so that $\PE\left[ \|\hg(\prm_{R_\totstep})\|^2\right]  \leq \epsilon$.
Assume that $\nbrMC$ satisfies
\begin{equation}
\label{eq:constraint-epsilon-SAEMiid}
0 < \epsilon \leq {2 \boundvarem_0} / ({n \nbrMC +  \boundvarem_1}) \eqsp.
\end{equation}
In this case (see \eqref{eq:expression-Tepsilon-HP}), the total number of iterations needed to guarantee an $\epsilon$-approximate stationary point is 
\begin{equation}
\label{eq:bound-complexity-SAEMiid-HP}
\totstep(\epsilon, \nbrMC) \eqdef \frac{ 8 \Vinit \Liplyap }{ v_{\min}^2} \; \frac{\boundvarem_0}{n \nbrMC \epsilon^2} \eqsp,
\end{equation}
which grows as $\epsilon^{-2}$ as $\epsilon \downarrow 0^+$.  The  bound in \eqref{eq:bound-complexity-SAEMiid-HP} is achieved by taking a constant step size $\step(\epsilon,\nbrMC) \eqdef \nicefrac{v_{\min}  n \nbrMC \epsilon}{(2\boundvarem_0 \Liplyap)}$. The minimal number of iterations $\totstep(\epsilon,\nbrMC)$ is inversely proportional to the Monte Carlo sample size $\nbrMC$, while the step size is proportional to $\nbrMC$. Increasing $\nbrMC$ allows more aggressive step sizes to be used and the number of iterations to be reduced accordingly.  As in the Mini-batch EM, let us evaluate the computational cost to understand the impact of this choice on the computational complexity. The computational cost depends on the number of calls to the  optimization map $\mapem$ (with cost $\cost_\mapem$) and the  cost of approximating each of the $n$ expectations $\barsem_i$ by a Monte Carlo sum with $\nbrMC$ terms (each term has a cost $\cost_{\MC}$). The cost of evaluating one oracle is
 $\cost_{\mapem} + n \nbrMC  \cost_{\MC}$, and  after $\totstep(\epsilon, \nbrMC)$ iterations, this cost is
\begin{equation} \label{eq:SAEMiid:complexity:case1}
\cost(\epsilon, \nbrMC) \eqdef \frac{8 \Vinit \Liplyap \boundvarem_0}{ v_{\min}^2 \epsilon^2}\cost_{\mapem}   \, \left(\frac{1}{ n \nbrMC}  + \frac{\cost_{\MC}}{\cost_{\mapem}} \right) \eqsp.
\end{equation}
Increasing the number  $\nbrMC$ of Monte Carlo samples reduces the cost due to  the optimization step $\mapem$ but will have no impact on the Monte Carlo cost. Since this discussion assumes \eqref{eq:constraint-epsilon-SAEMiid}, the maximal Monte Carlo sample size is given (up to appropriate roundings) by $
\nbrMC(\epsilon)  \eqdef n^{-1} ( \nofrac{2 \boundvarem_0}{\epsilon} - \boundvarem_1)$:
the "optimal" sample size is inversely proportional to the  accuracy $\epsilon$. It is interesting to note that the step size $\step(\epsilon,\nbrMC(\epsilon))$ associated with this choice of sample size is $\vminem(1- \epsilon \boundvarem_1/(2 \boundvarem_0))/\Liplyap > \vminem / (2 \Liplyap)$.

Assume now that $\nbrMC$ is such that $\epsilon \geq \nofrac{2 \boundvarem_0}{(n \nbrMC + \boundvarem_1)}$.
Using \eqref{eq:expression-Tepsilon-LP}, the total number of iterations needed to achieve an $\epsilon$-approximate stationary point is lower bounded by
\begin{equation}
\label{eq:bound-complexity-SAEMiid-LP}
\totstep(\epsilon,\nbrMC) \eqdef  \frac{4  \Vinit \Liplyap }{\vminem^2} \frac{1+ \nofrac{\boundvarem_1}{n \nbrMC}}{ \epsilon  } \,,
\end{equation}
which grows as $\epsilon^{-1}$ as $\epsilon \downarrow 0^+$.  The  bound in \eqref{eq:bound-complexity-SAEMiid-LP} is achieved by taking a constant step size $\step(\epsilon,\nbrMC) \eqdef \step_{\max}/2 = \vminem  n \nbrMC/ ( \Liplyap (n \nbrMC+\boundvarem_1))$. After $\totstep(\epsilon,\nbrMC)$ iterations, the total cost of evaluating the stochastic oracles is
\[
\cost(\epsilon, \nbrMC) \eqdef  \frac{4  \Vinit \Liplyap }{\vminem^2 \epsilon}  (n \nbrMC+\boundvarem_1) \, \cost_{\mapem}  \left( \frac{1}{n \nbrMC} + \frac{\cost_{\MC}}{\cost_{\mapem} } \right).
\]
This quantity is minimal with $\nbrMC(\epsilon) \eqdef \{n^{-1}(2 \boundvarem_0 /\epsilon - \boundvarem_1)\} \vee  \nbrMC^\star$ where
\[
\nbrMC^\star\eqdef  n^{-1}\sqrt{ \boundvarem_1 \cost_{\mapem}/\cost_{\MC}}.
\]
 When $\epsilon \leq 2 \boundvarem_0/ (n \nbrMC^\star+\boundvarem_1)$, then $\nbrMC(\epsilon) = n^{-1}( 2 \boundvarem_0 /\epsilon - \boundvarem_1)$ and  $\cost(\epsilon, \nbrMC(\epsilon))$ is equal to the previous case (see \eqref{eq:SAEMiid:complexity:case1}); otherwise, the cost  $\cost(\epsilon,  \nbrMC^\star)$ is lower. See the summary of the cost for SAEM-ES in Table~\ref{tab:EMtable}.

\paragraph{SAEM with self-normalized Importance Sampling (SAEM-IS)}
\label{sec:SAEM:complexity}
Consider the case of SAEM-IS:  conditionally to the past, the random variables $\{Z_{i,k+1}^j, 1  \leq j \leq \nbrMC, 1 \leq i \leq n \}$ are independent, and for all $i \in \{1, \ldots, n\}$ and $j$, the distribution of $Z_{i,k+1}^j$ is $\tilde \pi_{i}(z_i; \mapem(\prm_k))$.
Assume  \Cref{assumEM:var:SAEM:IS}. From \Cref{sec:EM:verifH:SAEM}, we have for any Monte Carlo batch size $\nbrMC$ large enough: $\eta_0 = \nicefrac{\boundvarem_0 }{\nbrMC}$, $\eta_1 = 1 + \nicefrac{\boundvarem_1}{\nbrMC}$, $\bsf_0 = \nicefrac{\clyap_{\bsf}}{\nbrMC}$ and $\bsf_1$ small enough so that $\vminem - \bsf_1 \geq \vminem/2$. The exact expressions of $\eta_0, \eta_1, \bsf_0$ and $\bsf_1$ in terms of $s_\star$, $\clyap_{\chi,\ell}$, $n$, $\nbrMC$ and $\vmaxem$
are given in \Cref{sec:appendix:SAEM}.
From \Cref{theo:main-result-quantitative}, we obtain:
\begin{proposition}\label{prop:SAEMis:quantitativebounds}
Assume \Cref{assumEM:expo,assumEM:Mstep,assumEM:supplementary}. Let $\sequence{\prm}[k][\nset]$ be the {\SA} sequence with random oracle \eqref{eq:oracle-SAEM}.  Assume in addition that  for all $i \in \{1,\ldots, n\}$ and $k \geq 0$, the random variables $\{Z_{i,k+1}^j, j \geq 1\}$ are  independent and sampled from the distribution $\pi_i(z_i; \mapem(\prm_k))$ and \Cref{assumEM:var:SAEM:IS} holds. Then for any $\totstep \in \nset$,
\begin{align*}
\PE\left[ \|\hg(\prm_{R_\totstep})\|^2\right] & \leq \frac{2 \clyap_{\bsf}  \sum_{k=1}^{\totstep} \step_{k}/\nbrMC}{\sum_{\ell=1}^{\totstep} \step_{\ell} \, (2 v_{\min} - \step_{\ell} \Liplyap (1+\boundvarem_1/ (n \nbrMC)))}
\\
& + \frac{2 \Vinit +   \Liplyap \boundvarem_0 \sum_{k=1}^{\totstep} \step^2_{k} / (n \nbrMC) }{ \sum_{k=1}^{\totstep}  \step_{k}  (2 v_{\min} - \step_{k} \Liplyap (1+\boundvarem_1/ (n \nbrMC)))},
\end{align*}
where $R_\totstep$ is the uniform random variable on $\{0, \ldots, T-1\}$ and $\Vinit$, $\Liplyap$ and $\vminem$ are defined by \eqref{eq:definition-biais} and \Cref{assumEM:supplementary}  respectively.
\end{proposition}
Let us discuss how to choose a constant step size $\step$,  the total number of iterations $\totstep$, and the Monte Carlo batch size $\nbrMC$ to satisfy the $\epsilon$-approximate stationary condition. SAEM-IS is not a USO algorithm: the stochastic oracles are biased approximations of the mean field. Consequently, $\Bterm \neq 0$ in \Cref{coro:main-result-quantitative} as this term does not depend on the step size $\step$ nor on the number of iterations $\totstep$. Nevertheless, $\nbrMC$ goes to infinity, we have $\bsf_0 \to 0$, $\bsf_1 \to 0$, which implies that  $\Bterm \to 0$. As such, $\Bterm$ can be made small by a convenient choice of $\nbrMC$.

From \Cref{coro:main-result-quantitative}, we have
\begin{equation}\label{eq:SAEMis:quantitativebounds}
\PE\left[ \| \hg(\prm_{R_\totstep}) \|^2 \right] \leq \frac{4 \clyap_{\bsf}}{\vminem \nbrMC} +  \frac{4 \sqrt{2 \Vinit \boundvarem_0 \Liplyap}}{\sqrt{\nbrMC} \sqrt{\totstep} \vminem} \vee \frac{16 \Vinit \Liplyap \eta_1}{  \totstep  \vminem^2},
\end{equation}
where $R_{\totstep}$ is the uniform random  variable on $\{0, \cdots, \totstep-1\}$. Such an upper bound is obtained with a constant step size
\[
\step  \eqdef \left( \frac{2 \Vinit \nbrMC}{\boundvarem_0 \Liplyap \totstep} \right)^{1/2} \wedge  \frac{\vminem - \bsf_1}{\Liplyap (1+ \nicefrac{\boundvarem_1}{\nbrMC})}.
\]

Choose $\kappa \in \ooint{0,1}$ such that for any $\epsilon \in \ocint{0, 2 \nicefrac{\boundvarem_0}{\boundvarem_1}}$,
\begin{equation}\label{eq:SAEM:conditionkappa}
\frac{4 \clyap_{\bsf}}{(1-\kappa) \vminem } +  \epsilon\boundvarem_1 \leq \frac{2 \boundvarem_0}{\kappa}.
\end{equation}
Let $\epsilon \in \ocint{0, 2 \nicefrac{\boundvarem_0}{\boundvarem_1}}$.   First, we choose $\nbrMC$ such that $\Bterm = 4 \clyap_{\bsf} /(\vminem \nbrMC) \leq  (1-\kappa) \epsilon$. This yields
\begin{equation}\label{eq:SAEMis:killbias}
\revisionupdates{\nbrMC} \geq 4 \clyap_{\bsf}/ ((1-\kappa) \vminem \epsilon).
\end{equation}
Now we choose $\totstep$ and $\nbrMC$ to make the second term in \eqref{eq:SAEMis:quantitativebounds} smaller than $\kappa \epsilon$. As in the comments following \Cref{coro:main-result-quantitative}, we distinguish two regimes: the second term in \eqref{eq:SAEMis:quantitativebounds} is lower than $\kappa \epsilon$ if  the number of iterations $\totstep$ is larger than
\[
\frac{16 \Vinit \Liplyap(1+ \nicefrac{\boundvarem_1}{\nbrMC})}{\vminem^2 \kappa \epsilon} \vee \frac{32 \Vinit \boundvarem_0 \Liplyap}{\nbrMC \vminem^2 \kappa^2 \epsilon^2},
\]
and this bound, seen as a function of $\epsilon$, defines  two regimes depending on the value of $\epsilon$ \wrt\ $2 \boundvarem_0 /(\kappa (\nbrMC+ \boundvarem_1))$.

In  the high-precision regime where $\epsilon \in \ocint{0, \nicefrac{2 \boundvarem_0}{\kappa(\nbrMC + \boundvarem_1)}}$, the number of iterations $\totstep$ is lower bounded by
\[
\totstep(\epsilon,\nbrMC) \eqdef \frac{32 \, \Liplyap \Vinit \boundvarem_0}{\vminem^2 \nbrMC \kappa^2 \epsilon^2}.
\]
The step size is  $\step(\epsilon, \nbrMC) \eqdef  \kappa\epsilon \vminem \nbrMC /(4 \Liplyap \boundvarem_0)$. Increasing the number of Monte Carlo points allows more aggressive step sizes and decreases the number of iterations. Nevertheless, it has an impact on the computational cost of the algorithm. The cost, per iteration, is the sum of the optimization cost when computing $\mapem$ (it is denoted by $\cost_{\mapem}$) and the Monte Carlo cost $n \nbrMC \cost_{\MC}$ when approximating each of the $n$ expectations $\barsem_i$ with $\nbrMC$ draws ($\cost_{\MC}$ denotes the cost of one Monte Carlo draw). After $\totstep(\epsilon, \nbrMC)$ iterations, it is equal to
\begin{equation}\label{eq:SAEMis:complexity:case1}
\cost(\epsilon, \nbrMC) = \frac{32 \, \Liplyap \Vinit \boundvarem_0}{\vminem^2 \kappa^2\epsilon^2} \cost_{\mapem} \left( \frac{1}{\nbrMC} + \frac{\cost_{\MC}}{\cost_{\mapem}} n \right).
 \end{equation}
 There is a gain in increasing $\nbrMC$;
nevertheless, in this high-precision regime, $\nbrMC$ is upper bounded by  $2 \boundvarem_0 / (\kappa\epsilon) - \boundvarem_1$, and it also satisfies \eqref{eq:SAEMis:killbias}: the definition of $\kappa$ (see \eqref{eq:SAEM:conditionkappa}) allows the choice $\nbrMC(\epsilon ) = 2 \boundvarem_0 /(\kappa\epsilon) - \boundvarem_1$.

In the low-precision regime, which corresponds to $\epsilon \geq \nicefrac{2 \boundvarem_0}{\kappa(\nbrMC + \boundvarem_1)}$, the number of iterations is lower bounded by
\[
\totstep(\epsilon, \nbrMC) \eqdef \frac{16 \Vinit \Liplyap (1+ \nicefrac{\boundvarem_1}{\nbrMC})}{\vminem^2  \kappa \epsilon}
\]
and the step size  $\step(\epsilon, \nbrMC)$ is larger than $\vminem /(2 \Liplyap (1+ \nicefrac{\boundvarem_1}{\nbrMC}))$. The computational cost is
\[
\cost(\epsilon, \nbrMC) \eqdef \frac{16 \Vinit \Liplyap (1+ \nicefrac{\boundvarem_1}{\nbrMC})}{\vminem^2 \kappa \epsilon}\cost_{\mapem} \left(1 + \frac{\cost_{\MC}}{\cost_{\mapem}} n \nbrMC \right).
\]
It is minimized with $\nbrMC(\epsilon) = \nbrMC^\star \vee (2 \boundvarem_0 /(\kappa \epsilon) - \boundvarem_1)$ where
\[
\nbrMC^\star \eqdef \sqrt{\boundvarem_1 \cost_{\mapem} / (n \cost_{\MC})}.
\]
When $\kappa \epsilon \leq 2 \boundvarem_0/ (\nbrMC^\star+\boundvarem_1)$, then $\nbrMC(\epsilon) = 2 \boundvarem_0 /(\kappa\epsilon) - \boundvarem_1$ and the cost $\cost(\epsilon, \nbrMC(\epsilon))$ is equal to the previous case (see \eqref{eq:SAEMis:complexity:case1}); otherwise, the cost  $\cost(\epsilon,  \nbrMC^\star)$ is lower. See the summary of the cost for SAEM-IS in Table~\ref{tab:EMtable}.

\newcommand{\STAB}[1]{\begin{tabular}{@{}c@{}}#1\end{tabular}}
\begin{table}[t]
    \resizebox{\linewidth}{!}{\centering
    \begin{tabular}{cc |c}
    \toprule 
    \multicolumn{2}{c|}{Stochastic EM Algorithms} & Computational Cost\\
    \midrule 
    \multirow{2}{*}{\STAB{\rotatebox[origin=c]{90}{Mini-batch EM}}} & \makecell{High precision regime\\ $\epsilon \in \ocint{0, \frac{2 \boundvarem_0}{ \lbatchem^\star+\boundvarem_1} }$} & $\displaystyle \frac{ 8 \Vinit \Liplyap \boundvarem_0}{ v_{\min}^2}\frac{\cost_{\barsem}}{\epsilon}  \, \left( \frac{1}{\epsilon} +   \frac{\cost_{\mapem}}{\cost_{\barsem}} \frac{1}{2 \boundvarem_0 - \epsilon \boundvar_1} \right)$ \\
    & \makecell{Low precision regime\\ $\epsilon \in \coint{ \frac{2 \boundvarem_0}{ \lbatchem^\star+\boundvarem_1} , \frac{2 \boundvarem_0}{\boundvarem_1} }$} & $\displaystyle \frac{8  \Vinit \Liplyap \boundvarem_0}{\vminem^2}  \frac{\lbatchem^\star+\boundvarem_1}{2 \boundvarem_0} \frac{\cost_{\barsem}}{\epsilon}    \left( 1+  \frac{  \cost_{\mapem}}{\cost_{\barsem} }  \frac{1}{\lbatchem^\star}\right)$ \\
    \midrule 
    \multirow{2}{*}{\STAB{\rotatebox[origin=c]{90}{SAEM-ES}}} & \makecell{High precision regime\\ $\epsilon \in \ocint{0, \frac{ 2 \boundvarem_0} {n \nbrMC^\star + \boundvarem_1} }$} & $\displaystyle \frac{8 \Vinit \Liplyap \boundvarem_0}{ v_{\min}^2 \epsilon}\cost_{\mapem}   \, \left(\frac{1}{2 \boundvarem_0 - \epsilon \boundvarem_1}  + \frac{\cost_{\MC}}{\cost_{\mapem}} \frac{1}{\epsilon} \right)$ \\
    & \makecell{Low precision regime\\ $\epsilon \in \coint{ \frac{ 2 \boundvarem_0} {n \nbrMC^\star + \boundvarem_1} , \frac{2 \boundvarem_0}{\boundvarem_1} }$} & $\displaystyle \frac{8  \Vinit \Liplyap \boundvarem_0 }{\vminem^2 \epsilon} \, \frac{n \nbrMC^\star+\boundvarem_1}{2\boundvarem_0} \, \cost_{\mapem}  \left( \frac{1}{n \nbrMC^\star} + \frac{\cost_{\MC}}{\cost_{\mapem} } \right)$ \\
    \midrule
    \multirow{2}{*}{\STAB{\rotatebox[origin=c]{90}{SAEM-IS}}} & \makecell{High precision regime\\ $\epsilon \in \ocint{0, \frac{ 2 \boundvarem_0}{\kappa(\nbrMC^\star + \boundvarem_1)} }$} & $\displaystyle \frac{32 \, \Liplyap \Vinit \boundvarem_0}{\vminem^2 \kappa \epsilon} \cost_{\mapem} \left( \frac{1}{2\boundvarem_0 - \kappa \epsilon  \boundvarem_1} + \frac{\cost_{\MC}}{\cost_{\mapem}} \frac{n}{\kappa\epsilon} \right)$ \\
    & \makecell{Low precision regime\\ $\epsilon \in \coint{ \frac{ 2 \boundvarem_0}{\kappa(\nbrMC^\star + \boundvarem_1)} , \frac{2 \boundvarem_0}{\boundvarem_1} }$} & $\displaystyle \frac{16 \Vinit \Liplyap (\nbrMC^\star+ \boundvarem_1)}{\vminem^2 \kappa\epsilon}\cost_{\mapem} \left(\frac{1}{\nbrMC^\star} + \frac{\cost_{\MC}}{\cost_{\mapem}} n  \right)$ \\
    \bottomrule
    \end{tabular}
    }
    \caption{Computation costs to find $\epsilon$-approximate stationary points for stochastic EM. The optimal batch size $\lbatchem^\star$ (Mini-batch EM), the Monte Carlo sample size $\nbrMC^\star$ (SAEM-ES/SAEM-IS), and the step sizes can be found in \Cref{subsec:stochastic-EM}. The low-precision regime does not exist if $\boundvarem_1 = 0$.}\vspace{-.2cm}
    \label{tab:EMtable}
\end{table}

\subsubsection{TD Learning}
\label{subsec:TD-learning}
TD(0) is an example of USO (see \Cref{def:USO}) with
$\boundbias_0=\boundbias_1 =0$. Using the results of \Cref{subsec:TD-learning-checking}, it is easily checked that $\bsf_0 = \bsf_1 = 0$, $\eta_0  = 6(1 + 2 \| \valuefunc[\star]\|_{\Dtd_{\statdistMRP}}^2)$, $\eta_1 = 3 (1+\lambda)^2$, $\step_{\max} = 2(1-\lambda)/(3(1+\lambda)^2)$ and $\omega_k = 2(1-\lambda) - 3 \step_k (1+\lambda)^2$.  We obtain the following result from \Cref{theo:main-result-quantitative} and \Cref{rem:convex-superlyap}, which extends  \cite{bhandari2018finite}.

\begin{proposition}
\label{prop:TD:main-convergence}
Assume \Cref{assum:TD:stationary-policy,assum:TD:normed:calR,assum:TD:sampling} and $\sup_s \| \feature(s) \|\leq 1$.
Consider the  TD(0) sequence defined in \eqref{eq:iteration-TD0}. Set $\totstep > 0$, and let $\sequence{\step}[k][\nset]$ be a sequence such that $0 < \step_k <2 (1-\lambda)/(3(1+\lambda)^2)$. Set
\begin{equation*}
\bar{\prm}_{\totstep} \eqdef \sum_{k=1}^{\totstep} \frac{ \step_{k} ( 2 (1-\lambda) - 3 \step_{k} (1+\lambda)^2  )}{\sum_{\ell=1}^{\totstep} \step_{\ell}  ( 2  (1-\lambda) - 3 \step_{\ell} (1+\lambda)^2  )} \prm_k.
\end{equation*}
Then,
\begin{multline*}
\PE[\|\valuefunc[\bar{\prm}_{\totstep}]-\valuefunc[\star] \|_{\Dtd_{\statdistMRP}}^2]  \\
\leq \frac{\PE[\| \prm_0 - \prm_\star\|^2]   +  6  \{1 + 2  \| \valuefunc[\star]\|_{\Dtd_{\statdistMRP}}^2\}
\sum_{k=1}^{\totstep} \step^2_{k}  }{\sum_{\ell=1}^{\totstep}  \step_{\ell} ( 2  (1-\lambda) -  3\step_{\ell} (1+\lambda)^2 )},
\end{multline*}
where  $\Sigma_{\statdistMRP}$ and $\boundvar_0$ are given in \eqref{eq:covariance-feature} and \eqref{eq:TD:boundvar_0} respectively.  $\prm_\star$ is any solution of $\valuefunc[\star] = \Feature \prm_\star$.
\end{proposition}

\Cref{prop:TD:main-convergence} bounds the mean squared distance between the value function estimator under the averaged iterate $\bar{\prm}_{\totstep}$ and the fixed point to the projected Bellman equation \eqref{eq:projected-bellman}. The  strength of this result is that the step sizes and the bound do not depend on the condition number of the feature covariance matrix (under \Cref{assum:TD:full-rank}, $\vmintd>0$).
Let us first consider the case of constant stepsize policy.  Set
\[
\step_\totstep \eqdef  \sqrt{\frac{2 \PE[\|\valuefunc[\prm_{0}]-\valuefunc[\star] \|_{\Dtd_{\statdistMRP}}^2]}{6 (1 + 2 \| \valuefunc[\star]\|_{\Dtd_{\statdistMRP}}^2)\totstep}} \wedge \frac{1-\lambda}{3(1+\lambda)^2}.
\]
Using \Cref{coro:main-result-quantitative}, we get:
\begin{corollary}
\label{coro:TD:main-convergence}
Assume \Cref{assum:TD:stationary-policy,assum:TD:normed:calR,assum:TD:sampling}.
Then, for any  $\totstep \geq 1$, setting  $\step_{k+1} \eqdef \step_{\totstep}$ for $k \in \{0,\dots,\totstep-1\}$,
 we get
\begin{multline*}
\label{eq:conclusion-simplified}
\PE[\|\valuefunc[\bar{\prm}_{\totstep}]-\valuefunc[\star] \|_{\Dtd_{\statdistMRP}}^2] \\ \leq  \frac{ 2 \, \sqrt{12 \Vinit \{1 + 2  \| \valuefunc[\star]\|_{\Dtd_{\statdistMRP}}^2\}}}{\sqrt{\totstep}  (1-\lambda) }
\vee   \frac{12 \Vinit (1+\lambda)^2} {\totstep  (1-\lambda)^2} \eqsp
\end{multline*}
where $\Vinit= \PE[\|\valuefunc[\prm_{0}]-\valuefunc[\star] \|_{\Dtd_{\statdistMRP}}^2 ]$ and $\bar{\prm}_{\totstep} \eqdef \totstep^{-1} \sum_{k=0}^{\totstep-1} \prm_k$.
\end{corollary}
Using \Cref{lem:TD:upper-lower-spectrum-covariance}, we observe that the Lyapunov functions $\lyap$ and $\superlyap$ are equivalent, \ie\ for all $\prm \in \rset^d$, we get
\begin{equation}
\label{eq:equivalence-lyapunov}
2 \sqrt{\vmintd} \lyap(\prm) \leq \superlyap(\prm) \leq 2 \lyap(\prm) \eqsp.
\end{equation}
The conclusions of \Cref{lem:TD:lyap,lem:TD:bound-mean-field} may be rewritten as
\begin{align*}
&\pscal{\nabla \lyap(\prm)}{\hg(\prm)} \leq -2 \sqrt{\vmintd} (1-\lambda) \lyap(\prm) \eqsp, \\
&\| \hg(\prm) \|^2 \leq 2 (1+\lambda)^2 \lyap(\prm) \eqsp,
\end{align*}
showing that \Cref{assum:field} and \Cref{assum:lyapunov} are satisfied with $\superlyap = \lyap$, $\rholyap \eqdef  2 \sqrt{\vmintd} (1-\lambda)$, $\clyap_{\hg,0} \eqdef 0$, and $\clyap_{\hg,1} \eqdef  2 (1+\lambda)^2$, $\boundvar_0$ and $\boundvar_1$ given by \eqref{eq:TD:boundvar_0}.  

We can therefore apply \Cref{coro:main-result-quantitative-fast}; note that under   \Cref{assum:TD:stationary-policy,assum:TD:normed:calR,assum:TD:sampling,assum:TD:full-rank}, $\boundbias_0 = \boundbias_1=0$ and $\clyap_\lyap < \infty$, \Cref{assumNA:uniform-bound-bias,assumNA:bias-is-small} are verified and $\bsf_0 = \bsf_1 = 0$.
This yields the following result,  which gives a convergence rate $O(\totstep^{-1})$.
Note that under \Cref{assum:TD:full-rank}, there exists a unique $\prm_\star \in \rset^d$  such that $\valuefunc[\star] = \Feature \prm_\star$.
\begin{corollary}
\label{coro:TD:fast}
Assume \Cref{assum:TD:stationary-policy,assum:TD:normed:calR,assum:TD:sampling,assum:TD:full-rank}. Let $\totstep \geq 1$ and set $\step_{k+1} \eqdef \tilde \step/(k+1+\totstep_0)$  for $k = \{0, \ldots, \totstep -1\}$, where $\tilde \step > \nicefrac{3}{\sqrt{\vmintd} (1-\lambda)}$ and $\totstep_0 \geq 2 \tilde \step (1+\lambda)^2/   (\sqrt{\vmintd} (1-\lambda))$. Then
\begin{align*}
&\PE\left[ \|\prm_{\totstep} - \prm_\star \|^2   \right] \leq  \left(  \frac{\totstep_0}{\totstep+\totstep_0}\right)^{\tilde \step \sqrt{\vmintd} (1-\lambda)}  \PE\left[ \|\prm_0 - \prm_\star \|^2   \right] \\
& \qquad + \frac{12 \tilde \step}{(\totstep+\totstep_0)\sqrt{\vmintd} (1-\lambda)} \left(1 + \{\lambda^2 +1\} \| \valuefunc[\star] \|^2_{\Dtd_{\statdistMRP}}\right).
\end{align*}
\end{corollary}
In this case, however, the choice of stepsize $\step_k= \tilde{\step}/(k+\totstep_0)$,  depends on the minimal eigenvalue $\vmintd$ of the feature covariance matrix through the constants $\tilde{\step}$, $\totstep_0$.
These results are in the spirit of the \emph{robust {\SA}} introduced by \cite{nemirovski2009robust}.
With constant stepsize policy, after $\totstep$ iterations, the mean squared distance between the value function estimates under the averaged iterate and the fixed point to the projected Bellman equation decreases at the rate $O(\totstep^{-1 / 2})$. Of course, this is worse than the rate $O(\totstep^{-1})$ for the decreasing step size SA algorithm.
However, the expected error bounds of \Cref{coro:TD:main-convergence} are guaranteed regardless of the knowledge of $\vmintd$. The ${\cal O}(\totstep^{-1/2})$ bound still holds with any constant step size $\step_k= \tilde{\step}/\sqrt{\totstep}$, $k \in \{1,\dots,\totstep\}$ with $\tilde{\step} >0$, an error in the choice of $\tilde{\step}$  has a linear effect on the error bound in $\max \left\{\tilde{\step}, \tilde{\step}^{-1}\right\}$. This is to be  compared with a potentially catastrophic effect of an appropriate choice of the hyperparameters $\tilde{\step}$ and $\totstep_0$ in \Cref{coro:TD:fast}. These observation justifies the `robustness' of the method as "fine tuning" of the step sizes to the objective function is not necessary.

\section{Almost-sure convergence} \label{sec:asymptotic-convergence}
This section overviews the asymptotic convergence analysis of {\SA} scheme where we study the behavior of \eqref{eq:sa} when the optimization horizon tends to infinity ($k \to \infty$) and  we will use decreasing step sizes. At the first glance, these asymptotic convergence results may appear less powerful than the non-asymptotic bounds in \Cref{sec:nonasymp}, yet we emphasize that these results are presented in flavor of the \emph{almost-sure} convergence towards one of the equilibrium points set. In contrast, the non-asymptotic bounds only show the convergence towards a \emph{near-equilibrium} point within a finite number of iterations, and the results are often given in expectation. In fact, the first results in {\SA} were obtained in the almost-sure convergence framework in the pioneering works of \cite{robbins1951stochastic} and \cite{blum1954multidimensional}; see \cite[Chapter~1-2]{kushner2003stochastic} for a historical introduction. Nevertheless, both types of convergence results are equally important for our understanding of the {\SA} schemes.

\subsection{The ODE method}
\label{subsec:ODE-method}

A powerful method for establishing almost-sure convergence results is the so-called ordinary differential equation (ODE) method, which allows us to relate the almost-sure limit point of {\SA} schemes (see \eqref{eq:sa}) with the limiting sets of the flow of the autonomous ODE
\begin{equation}
\label{eq:ODE-method}
{\rmd \prm} / {\rmd t}=\hg(\prm).
\end{equation}
The key element is a detailed analysis of the \emph{flow} associated to vector field $\hg$. Let $\flow: \rset \times \rset^d \to \rset^d$ $(t,\dpy) \mapsto \flow[][](t,\dpy)= \flow[t][](\dpy)$, be a continuous function. The family $\flow= (\flow[t])_{t \in \rset}$ is called a flow of $\rset^d$, if $\flow[0]= \Id$ and for all $(t,s) \in \rset^2$, $\flow[t] \circ \flow[s]= \flow[t+s]$. $\flow$ is a semi-flow if we substitute the above $\rset$ by $\rset_+$. 
\revisionupdates{The forward orbit of $\dpy \in \rset^d$ is the set $\orbit^{+}(\dpy)=\left\{\flow[t](\dpy): t \geq 0\right\}$ and the orbit of $\dpy$ is $\orbit(\dpy)=\left\{\flow[t](\dpy): t \in \rset\right\}$.}

The continuous vector field $\hg$ on $\rset^d$ is said to have unique integral curve if there exists a flow $\flow[][\hg]: \rset \times \rset^d \to \rset^d$, $(t,\dpy) \mapsto \flow[][\hg](t,\dpy)= \flow[t][\hg](\dpy)$, which is differentiable with respect to $t$ satisfying, for all $t \in \rset$,
\[
{\rmd \flow[t][\hg]} / {\rmd t}= \hg(\flow[t][\hg]), \quad \flow[0][\hg](\dpy)=\dpy \eqsp.
\]
A point $\dprm_\star \in \rset^d$ is an \emph{equilibrium} if $\flow[t](\dprm_\star)=\dprm_\star$ for all $t \in \rset$. The set $\equilibrium(\flow[][\hg])$ of the equilibrium point of the flow $\flow[][\hg]$ coincides with roots of $\hg$:
\begin{equation}
\label{eq:root-equilibrium}
\equilibrium(\flow[][\hg]) = \{ \dprm_\star \in \rset^d, \hg(\dprm_\star)= 0 \} \,.
\end{equation}
A set $\Lambda \subset \rset^d$ is  \emph{invariant} (resp. \emph{positively invariant}) for the flow $\flow$ if $\flow[t](\Lambda) \subset \Lambda$ for all $t \in \rset$ (resp. for all $t \in \rset_+$). In this case, we denote by $\flow \vert \Lambda$ the flow (resp. semi-flow) restricted to $\Lambda$. If $\dprm_\star$ is an equilibrium point, then the set $\{\dprm_\star\}$ is invariant.

The ODE method was introduced by \cite{ljung1977analysis} and further refined in \cite{ljung1983theory}
(with recursive system identification in mind) and extensively studied thereafter; see, e.g., the books \cite{kushner2003stochastic}, \cite{borkar2009stochastic}, \cite{benveniste2012adaptive} for a comprehensive introduction and further references. The possibility of associating limit points of stochastic approximation procedure and a subset of the family of  invariant sets of the flow $\flow[][\hg]$  - which includes equilibrium points but might include more "complex " sets -  is the main motivation for the work of \cite{ljung1977analysis}, which extends older contributions on this subject. Previous results in this field have been established under simple and often unverifiable assumptions about the flow $\flow[][\hg]$, for example the existence of a single global "asymptotically stable" equilibrium \cite{robbins1951stochastic} or of a finite number of equilibrium points whose "basins of attraction" cover the whole space.
The success of the ODE method stems from the fact that many  results of the rich theory of dynamical systems are readily available; see \cite{benaim1999dynamics},  \cite[Chapter~4]{kushner2003stochastic}, and \cite[Chapter~2]{borkar2009stochastic}.

We will consider two situations. In the first, the most classical, the bias disappears asymptotically, i.e. $\lim_{k\to \infty} \boundbias_{\ell,k}= 0$, $\ell=0,1$. There is an extensive literature on this subject: we give below a very brief overview of these results, which are mainly inspired by \cite{benaim1996dynamical, benaim1999dynamics}; see also \cite{kushner2003stochastic,borkar2009stochastic}. We then extend these results to the case where the bias does not vanish but is bounded by a sufficiently small constant, by extending the results from \cite{tadic2011asymptotic,tadic2017asymptotic,ramaswamy2017analysis}. To keep the presentation concise, we report the results with only a sketch of proofs.
\subsection{Limit set of SA with vanishing bias}
\label{subsec:limit-sets}
We consider the following \emph{deterministic} sequence, which is a perturbed Euler discretization of the ODE \eqref{eq:ODE-method}
\begin{equation}
\label{eq:sa-sequence}
\dprm_{k+1}= \dprm_k + \step_{k+1} \{ \hg(\dprm_k) + \dnoise_{k+1} + \dbias_{k+1} \}.
\end{equation}
We  first strengthen the conditions on the stepsize sequence:
\begin{assumSA}
\label{assum:SA:stepsize}
$\sequence{\step}[k][\nset]$ is a non increasing sequence of positive numbers such that $\sum_{\revisionupdates{k=1}}^\infty \step_k=\infty$ and
$\lim _{k \rightarrow \infty} \step_k=0$.
\end{assumSA}
For the mean field $\hg$, we strengthen \Cref{assum:lyapunov} used in the non-asymptotic analysis, by assuming that the field $\hg$ is globally Lipschitz; it guarantees the existence and the uniqueness of the solutions of the ODE \eqref{eq:ODE-method} which can be extended to $\rset$: 
\begin{assumSA}
\label{assum:SA:vector-field}
The vector field $\hg$ is Lipschitz continuous, \ie\ for all $\dprm, \dprm' \in \rset^d$,
$\| \hg(\dprm) - \hg(\dprm') \| \leq \Lip{\hg} \| \dprm - \dprm' \|$.
\end{assumSA}
A key assumption in the analysis (and one of the most annoying one) is that the sequence is bounded.
\begin{assumSA}
\label{assum:SA:boundedness}
$\sup_{k \in \nset} \| \dprm_k \| < \infty$.
\end{assumSA}
 This may seem innocuous since we are interested in convergence, but the stability of {\SA} is a non-trivial issue in general. When \Cref{assum:SA:boundedness} is not satisfied - which unfortunately occurs in many practical examples - a classical approach is to project the iterates on a compact set. The projection then modifies the underlying dynamical system, which  becomes a differential inclusion; see \cite[Chapters~4-5]{kushner2003stochastic} for details. Another possibility is to project onto a growing sequence of compact sets; this has been advocated by \cite{andradottir1995stochastic} and recently studied in great detail (with Markovian noise) in \cite{andrieu2014markovian,fort:etal:2016}. Another solution is to reinitialize the sequence upon crossing a growing sequence of boundaries, as suggested by \cite{chen1986stochastic} and further worked out in \cite{andrieu2005stability}.
Technical details to prove stability under these modifications go beyond this survey.

Define $\timecum_0 \eqdef 0$ and for $k \in \nset^*$, $\timecum_k \eqdef \sum_{\ell=1}^k \step_\ell$. The "inverse" of $k \rightarrow \timecum_k$ is the map $m: \rset_{+} \rightarrow \nset$ defined by
$$
m(t) \eqdef \sup \left\{k \geq 0: t \geq \timecum_k\right\}.
$$
By construction, $m(\timecum_k)=k$.
We denote for $\totstep >0$, $\mathcal{I}_k(\totstep) \eqdef \{k+1, \ldots, m(\timecum_k+\totstep)\}$. We finally assume that the perturbations satisfy the following assumption.
\begin{assumSA}
\label{assum:SA:noise}
It holds $\lim_{k \to \infty} \dbias_k=0$ and that
for all $T>0$,
\[ \textstyle
\lim _{n \rightarrow \infty} \sup \left\{\left\|\sum_{\ell=n}^{k-1} \step_{\ell+1} \dnoise_{\ell+1}\right\|: k \in \mathcal{I}_n(T) \right\}=0 \eqsp.
\]
\end{assumSA}
This condition is often referred as the \emph{asymptotic rate of change} (see \eg\ \cite[Chapter~5, pp.137-138]{kushner2003stochastic}).
We compare the sequence $\sequence{\dprm}[k][\nset]$ with the flow induced by the vector field $\hg$.
For a sequence $\sequence{\dpy}[k][\nset]$ in $\rset^d$, we define the continuous time affine and piecewise constant interpolated functions $\dpY$, $\bar{\dpY}: \rset_{+} \to \rset^d$ respectively by
\begin{equation}
\label{eq:interpolated-function}
\dpY\left(\timecum_k+s\right)=\dpy_k+s \frac{\dpy_{k+1}-\dpy_k}{\timecum_{k+1}-\timecum_k} \text {, and } \bar{\dpY}\left(\timecum_k+s\right)=\dpy_k
\end{equation}
for all $k \in \nset$ and $0 \leq s<\step_{k+1}$.
We define, for $(t,T) \in \rset_+^2$,
\begin{equation}
\label{eq:definition-Delta}
\textstyle \Delta_{\dnoise}(t,T) \eqdef  \sup _{0 \leq \timeinc \leq T}\left\|\int_t^{t+\timeinc} \bar{\dNoise}(s)  \, \rmd s\right\|.
\end{equation}
Note that \Cref{assum:SA:noise} is equivalent to $ \lim_{t \rightarrow \infty} \Delta_{\dnoise}(t, T)=0$. 
With these notations, \eqref{eq:sa-sequence} writes, for $t \geq 0$,
\[
\textstyle \dPrm(t)-\dPrm(0)=\int_0^t \, \left( \hg(\bar{\dPrm}(s))+\bar{\dNoise}(s)  + \bar{\dBias}(s) \right)  \, \rmd s.
\]
The first key result for {\SA}  algorithms is the convergence of the linear interpolation to the ODE flow:
\begin{proposition}[\protect{\cite[Proposition~4.1]{benaim1999dynamics}}]
\label{prop:SA:asymptotic-trajectories}
Assume \Cref{assum:SA:stepsize,assum:SA:vector-field,assum:SA:boundedness,assum:SA:noise}. Then for all $T > 0$,
\begin{equation*}
\lim_{t \to \infty} \sup _{0 \leq \timeinc \leq T}\left\|\dPrm(t+ \timeinc)-\flow[\timeinc][\hg](\dPrm(t))\right\| =0.
\end{equation*}
\end{proposition}
Let $\delta>0$, $T>0$. A $(\delta, T)$-\emph{pseudo-orbit} for a flow $\flow$ from $\bm{a},\bm{b} \in \rset^d$ is a finite sequence of partial trajectories: there exist $N$ and time instants $\{t_i\}_{i=1}^N \subset \ocint{0,T}$ and  $\dpy_1, \dots, \dpy_N \in \rset^d$ satisfying
 $\| \dpy_0 - \bm{a} \| <\delta$, $\| \flow[t_j](\dpy_j) - \dpy_{j+1} \| \leq \delta$ for $j=0, \ldots, k-1$
and $\| \dpy_k  - \bm{b} \| \leq \delta$.
 We write $ \bm{a} \hookrightarrow^{\flow}_{\delta, T} \bm{b}$ if there exists a $(\delta, T)$-pseudo-orbit  for the flow $\flow$ from $\bm{a}$ to $\bm{b}$.
 We write $\bm{a} \hookrightarrow^{\flow} \bm{b}$ if $\bm{a} \hookrightarrow^{\flow}_{\delta, T} \bm{b}$ for every $\delta>0, T>0$. 
If $\bm{a} \hookrightarrow^{\flow} \bm{a}$ then $\bm{a}$ is a \emph{chain recurrent point} for the flow $\flow$. The set of chain-recurrent points of the flow $\flow$ is denoted $\chainrecurrent(\flow)$. It is easy to verify that $\chainrecurrent(\flow)$ is a closed positively invariant set and that $\equilibrium(\flow) \subset \chainrecurrent(\flow)$. When $\chainrecurrent(\flow)$ is compact, then it is invariant (see \cite[Theorem~5.5]{benaim1999dynamics}).
A subset $\Lambda$ is said \emph{internally chain-recurrent} if $\Lambda$ is a nonempty compact invariant set of which every point is chain-recurrent for the restricted flow $\flow \vert \Lambda$ (i.e., $\chainrecurrent(\flow \vert \Lambda)=\Lambda)$.

We now have all the essential notions to formulate the central result  for the convergence of {\SA} sequences.
In the form given below, the result is due to \cite{benaim1996dynamical}; see also \cite{benaim1999dynamics}.
Denote by  $\limset(\sequence{\dprm}[k][\nset])$ the limit set of the sequence  $\sequence{\dprm}[k][\nset]$:  $\dprm_\star \in \limset(\sequence{\dprm}[k][\nset])$,  if there exists a sequence $\sequence{n}[k][\nset]$, satisfying $\lim_{k \to \infty} n_k = +\infty$ and $\lim_{k \to \infty} \dprm_{n_k}= \dprm_\star$.  
\begin{theorem}[\protect{after \cite[Theorem~1.2]{benaim1996dynamical}}]
\label{theo:convergence-sa}
Let $\sequence{\dprm}[k][\nset]$ be the sequence given by \eqref{eq:sa-sequence}.
Assume \Cref{assum:SA:stepsize,assum:SA:vector-field,assum:SA:boundedness,assum:SA:noise}. 
Then $\limset(\sequence{\dprm}[k][\nset])$ is a connected  internally chain-recurrent set for the flow $\flow[][\hg]$.
\end{theorem}
It is shown in \cite[Theorem~1.3]{benaim1996dynamical} that this result cannot be improved, in the sense that any connected internally chain-recurrent set of the vector field $\hg$ is the limit set of a sequence \eqref{eq:sa-sequence} which satisfies the assumptions \Cref{assum:SA:stepsize,assum:SA:vector-field,assum:SA:boundedness,assum:SA:noise}; see \cite[Theorem~1.3]{benaim1996dynamical}. One must therefore be careful that the set of limit points of sequences \eqref{eq:sa-sequence} are not just the equilibrium of the vector field $\hg$, but can be a priori larger sets on which the vector field $\hg$ does not necessarily vanish. Thus, to obtain guarantees to converge only to equilibrium points, one must be able to guarantee that the only connected internally recurrent sets are included in $\{\prm \in \rset^d, \hg(\prm)=0\}$, which is, of course, possible but  typically requires additional assumptions beyond \Cref{assum:SA:vector-field}; see \cite[Chapter~2, Corollary~4]{borkar2009stochastic}.

The characterization of connected internally chain recurrent sets for the flow $\flow[][\hg]$ is simplified when the vector field $\hg$ is associated to a  \emph{Lyapunov function}. A continuous function $\lyap: \rset^d \to \rset$ is called a \emph{Lyapunov function} for $\Lambda$ if the function $t \in \rset_{+} \to \lyap(\flow[t][\hg](\dpy))$ is constant for $\dpy \in \Lambda$ and strictly decreasing for $\dpy \in \rset^d \setminus \Lambda$. If $\Lambda$ coincides with the equilibrium set of the flow associated with the vector field $\hg$, $\{ \dprm \in \rset^d, \hg(\dprm)=0 \}$, then $\lyap$ is  a \emph{strict Lyapunov function} and the flow $\flow[\hg]$ is said to be \emph{gradient-like}.
\begin{corollary}[\protect{\cite[Proposition~6.4, Corollary~6.6]{benaim1999dynamics}}]
\label{coro:convergence-sa} Let  $\sequence{\dprm}[k][\nset]$ be the sequence given by \eqref{eq:sa-sequence}.
Assume \Cref{assum:SA:stepsize,assum:SA:vector-field,assum:SA:boundedness,assum:SA:noise}. 
Assume in addition that
\begin{enumerate}[label=\alph*),nosep]
\item $\flow[][\hg]$ admits a Lyapunov function $\lyap$ for a compact invariant set $\Lambda$ of the flow $\flow[][\hg]$.
\item $\lyap(\Lambda)$ has an empty interior.
\end{enumerate}
Then $\limset(\sequence{\dprm}[k][\nset])$ is contained in $\Lambda$ and $\lyap$ is constant on $\limset(\sequence{\dprm}[k][\nset])$. If in addition $\lyap$ is a strict Lyapunov function for $\flow[][\hg]$, then $\limset(\sequence{\dprm}[k][\nset]) \subset \{ \hg = 0 \}$.
\end{corollary}
We conclude this short presentation with a practical characterization of Lyapunov functions which directly connects almost-sure convergence  with the assumption \Cref{assum:lyapunov} we used to establish the non-asymptotic bounds.
\begin{proposition}
\label{prop:characterization-lyapunov}
Assume that there exists a continuously differentiable function  $\lyap: \rset^d \to \rset$ such that, for any $\dprm \in \rset^d$, $\pscal{\nabla \lyap(\dprm)}{\hg(\dprm)} \leq 0$. Define
\begin{equation}
\label{eq:Lambda-lyap}
\Lambda_{\lyap} \eqdef \{\dprm \in \rset^d, \pscal{\nabla \lyap(\dprm)}{\hg(\dprm)}=0\}.
\end{equation}
Assume that
\begin{enumerate}[label=\alph*),nosep]
\item $\Lambda_\lyap$ is compact and invariant for $\flow[][\hg]$.
\item $\lyap(\Lambda_\lyap)$ has an empty interior.
\end{enumerate}
Then, every connected internally chain-recurrent set $L \subset \Lambda_\lyap$ is contained in $\Lambda_\lyap$ and $V$ restricted to $L$ is constant.  If in addition $\Lambda_{\lyap}= \{ \hg = 0 \}$ then $\lyap$ is a strict Lyapunov function.
\end{proposition}
\begin{proof}
Note indeed that, by the chain rule, we get,
\begin{equation}
\label{eq:chain-rule}
\frac{\rmd }{\rmd t} \lyap(\flow[t][\hg](\dpy))= \pscal{\nabla \lyap (\flow[t][\hg](\dpy))}{\hg (\flow[t][\hg](\dpy))}.
\end{equation}
If the set $\Lambda_\lyap$ is compact and invariant, then $\lyap$ is a Lyapunov function for $\Lambda_\lyap$.
Indeed, since $\Lambda_\lyap$ is invariant, then for any $\dpy \in \Lambda_\lyap$, $\flow[t][\hg](\dpy) \in \Lambda_\lyap$ for all $t \geq 0$. It follows from \eqref{eq:chain-rule} that $t \mapsto \lyap(\flow[t][\hg](y))$ is constant.  If $\dpy \in \rset^d \setminus \Lambda_{\lyap}$, then for all $t \in \rset_+$, $\flow[t][\hg](\dpy) \not \in \Lambda_\lyap$ (since $\Lambda_\lyap$ is invariant) and \eqref{eq:chain-rule} shows that $t \mapsto \lyap(\flow[t][\hg](\dpy))$ is strictly decreasing.
If the set $\Lambda_\lyap$ (see \eqref{eq:Lambda-lyap}) is compact and invariant for the flow $\flow[t][\hg]$ associated to the vector-field $\hg$, and if $\lyap(\Lambda_{\lyap})$ has an empty interior, then every connected internally chain recurrent set for the flow is included in $\Lambda_{\lyap}$.
\end{proof}
\noindent A direct proof of \Cref{theo:convergence-sa} for the field $\hg$ satisfying the assumptions of \Cref{prop:characterization-lyapunov} is in \cite{delyon1999convergence} and refined by \cite{andrieu2005stability}. 

\vspace{.1cm}
\noindent \textbf{Almost-sure convergence.} Consider the sequence $\sequence{\prm}[k][\nset]$ defined by \eqref{eq:sa}. We first show that under \Cref{assum:field}-\Cref{assum:lyapunov} and \Cref{assum:SA:stepsize}-\Cref{assum:SA:vector-field}, the assumptions \Cref{assum:SA:boundedness}-\Cref{assum:SA:noise} are satisfied with probability 1.
To apply the results above, we rewrite \eqref{eq:sa} as:
\begin{align*}
\prm_{k+1} = \prm_k + \step_{k+1} \{\hg(\prm_k) + \Noise_{k+1} + \Bias_{k+1} \},
\end{align*}
where we have set
\begin{align}
\label{eq:definition-Noise}
\Noise_{k+1}&\eqdef \Hg(\prm_k,\State_{k+1})-\CPE{\Hg(\prm_k,\State_{k+1})}{\mcf_k} \eqsp,\\
\Bias_{k+1}& \eqdef \CPE{\Hg(\prm_k,\State_{k+1})}{\mcf_k} - \hg(\prm_k) \eqsp.
\end{align}
By construction $\sequence{\Noise}[k][\nset]$ is a martingale difference sequence adapted to the filtration $\mcf=\sequence{\mcf}[k][\nset]$, where for $k \in \nset$, $\mcf_k \eqdef \sigma(\prm_0, \State_{\ell}, \ell=1,\dots,k)$.
We can then use \Cref{theo:convergence-sa} to establish the almost-sure  convergence of the sequence $\sequence{\prm}[k][\nset]$. The proof is in \Cref{app:proof-ODE}.
\begin{tcolorbox}[boxsep=2pt,left=4pt,right=4pt,top=3pt,bottom=3pt]
\begin{theorem} 
\label{theo:almost-sure-cvgce}
Assume that $\sum_{k=0}^\infty \step_k = \infty$, $\sum_{k=0}^\infty \step_k^2 < \infty$,
\begin{enumerate}[label=\roman*),nosep]
\item \label{item:theo:almost-sure-cvgce:noise} \Cref{assum:field} holds with $\lim_k \clyap_\lyap\boundbias_{0,k} = \lim_{k \to \infty} \clyap_\lyap \boundbias_{1,k} =0$ and  $\sum_{k=0}^\infty  \step_k \clyap_\lyap  \sqrt{\boundbias_{0,k}} < \infty$, where $\clyap_\lyap$ is  in \eqref{eq:definition:clyap-lyap}.
\item \label{item:theo:almost-sure-cvgce:lyap} \Cref{assum:lyapunov} holds with  the function $\lyap$ satisfying  $\lim_{\| \dprm \| \to \infty} \lyap(\dprm) = \infty$.
\end{enumerate}
Then, with probability one,  the sequence $\sequence{\prm}[k][\nset]$ converges, $\lim_{k \to \infty} \lyap(\prm_k)$ exists and  $\limset(\sequence{\prm}[k][\nset])$ is a connected internally chain recurrent set of the vector-field $\hg$. If in addition, $\Lambda_\lyap \eqdef  \{ \dprm \in \rset^d, \pscal{\nabla\lyap(\dprm)}{\hg(\dprm)=0}\}$ is compact and invariant for $\flow[][\hg]$ and $\lyap(\Lambda_\lyap)$ has an empty interior, then with probability 1, $\limset(\sequence{\prm}[k][\nset]) \subset \Lambda_{\lyap}$. Finally, if $\PE\left[ \lyap(\prm_0)\right]<\infty$, then $\sup_k \PE\left[\lyap(\prm_k) \right]< \infty$.
\end{theorem}
\end{tcolorbox}

\subsection{Limit set of SA sequences with bounded bias}
\label{subsec:asymptotic-convergence-bias}
We now briefly discuss the behavior of the stochastic approximation when the bias is bounded but does not tend to 0. There are few results in the literature, and all of these results were obtained for stochastic gradient algorithms; see \cite{tadic2011asymptotic,tadic2017asymptotic,ramaswamy2017analysis,ramaswamy2018stability} . However, the analysis of the proofs shows that these results also hold for general stochastic approximation algorithms. Presenting these results in details would lead us to introduce a large number of delicate mathematical concepts. Unlike the rest of this tutorial, we present the results in a more informal manner. The proofs  are based on the use of differential inclusion methods, which generalize ODEs. The idea is thus to see the bias term as a bounded perturbation of the "unbiased" dynamic; these results build on the work of  \cite{benaim2005stochastic,benaim2006stochastic,benaim2012perturbations,davis2020stochastic,majewski2018analysis}.
Basically, we replace the recursion \eqref{eq:sa-sequence} by
\begin{equation}
\label{eq:sa-sequence-di}
\dprm_{k+1}= \dprm_k + \step_{k+1} \{   \bm{g}_{k+1} + \dnoise_{k+1}  \},
\end{equation}
where $\bm{g}_{k+1} \in \DIhg(\dprm_k) :=\ball{\hg(\dprm_k)}{\boundbias_0}$, with $\ball{\dpy}{\delta} := \{ \dpy', \|\dpy'-\dpy \|\leq \delta\}$.
To analyse such sequence, we should replace the ODE by a differential inclusion (DI). $\DIhg: \rset^d \rightrightarrows \rset^d$ is a set valued map in the sense that for each $\dpy \in \rset^d$, we have that $\DIhg(\dpy)$ is a subset of $\rset^d$ (in this case, a ball of radius $\delta$, centered at $\hg(\dprm)$). To determine the limit points of \eqref{eq:sa-sequence-di},
we introduce the DI
\begin{equation}
\label{eq:differential-inclusion}
{\rmd \dprm(t)} / {\rmd t}  \in \DIhg(\dprm(t)) .
\end{equation}
We say that an absolutely continuous curve (a.c.) $\dprm: \rset_+ \rightarrow \rset^d$ is a \emph{solution} if  \eqref{eq:differential-inclusion} holds for almost every $t \in \rset_+$.

Various notions of continuity exist for set valued maps. The one that will be important for us is the notion of upper semicontinuity.
A set valued map $\operatorname{H}: \rset^d \rightrightarrows \rset^d $ is upper semi continuous (u.s.c.) at $\dprm \in \rset^d$ if for every $U$, a neighborhood of $\operatorname{H}(\dprm)$, there is $\delta>0$ such that
$$
\|\dprm'-\dprm\| \leq \delta \Longrightarrow \operatorname{H}(\dprm') \subset U .
$$
Under \Cref{assum:SA:vector-field}, the set-valued map $\DIhg$ is upper semi continuous (u.s.c.) and it follows from
\cite[Chapter~2.1]{aubin2012differential} (see also \cite[Chapter~4.1]{clarke2008nonsmooth}) that
\begin{proposition}
Assume \Cref{assum:SA:vector-field}. For every $\dpy \in \rset^d$, there exists a solution to \eqref{eq:differential-inclusion} such that $\dprm(0)= \dpy$.
\end{proposition}
\noindent Denote by $\SolutionDI$ the set of solutions. The \emph{set-valued semi-flow} $\flow[][\DIhg]$ associated to the set-valued map $\DIhg$ is defined by
\begin{equation}
\label{eq:set-valued-semiflow}
\flow[t][\DIhg](\dpy)= \{ \dprm(t), \dprm \in \SolutionDI, \dprm(0)= \dpy \},~t \in \rset_+ \eqsp.
\end{equation}
As above, we denote by $\dPrm$ the linear interpolation of the sequence $\dpy$; see \eqref{eq:interpolated-function}.
Let $\sequence{s}[k][\nset]$ be an increasing sequence of positive numbers such that $\lim_{k \to \infty} s_k= \infty$.
Denote by $\dPrm_k(t)= \dPrm(s_k+t)$, for $t \in \rset^+$.
The following result follows from \cite[Theorem~3.2]{majewski2018analysis} (see also \cite[Theorem~4.2]{benaim2012perturbations}).
It shows that the limits points  of shifted in time linear interpolation converge to solutions of the differential inclusion \eqref{eq:differential-inclusion}. The formulation is weaker than the one obtained for the ODE in  \eqref{prop:SA:asymptotic-trajectories} (the statement differs from \cite[Theorem~4.2]{benaim2005stochastic} which is wrong as it is stated).
\begin{proposition}
\label{prop:SA:convergence-interpolated}
Assume \Cref{assum:SA:vector-field,assum:SA:stepsize,assum:SA:boundedness,assum:SA:noise}. For any increasing sequence $\sequence{n}[k][\nset]$
there exists a subsequence $\sequence{\tilde{n}}[k][\nset] \subset \sequence{n}[k][\nset]$ and an absolutely continuous function $\dPrm_\infty$ such that for any $T > 0$,
\[
\lim_{k \to \infty} \sup_{0 \leq \timeinc \leq T} \| \dPrm_{\tilde{n}_k}(\timeinc)- \dPrm_{\infty}(\timeinc) \|=0 \eqsp.
\]
Moreover, $\dPrm_\infty$ is a solution of the differential inclusion \eqref{eq:differential-inclusion}.
\end{proposition}
A set $\Lambda \subset \rset^d$ is \emph{invariant}  if for every $\dpy \in \Lambda$ there exists a  trajectory $\dprm$ contained in $\Lambda$ with $\dprm(0)=\dpy$. See \cite[Section~3.2]{benaim2005stochastic} for further discussion.
Let $\Lambda \subset \rset^d$ be an invariant set of \eqref{eq:differential-inclusion} and consider $\flow[][\DIhg] \Vert \Lambda$, the set-valued flow $\flow[][\DIhg]$ restricted to $\Lambda$. For $\bm{a}, \bm{b} \in \rset^d$, we write $\bm{a} \hookrightarrow \bm{b}$ if for every $\delta>0$ and $T>0$ there exists an integer $n \in \nset$, solutions $\dpy_1, \ldots, \dpy_n$ to \eqref{eq:differential-inclusion}, and real numbers $t_1, t_2, \ldots, t_n$ greater than $T$ such that $\dpy_i(s) \in \Lambda$ for all $0 \leq s \leq t_i$ and  $i=1, \ldots, n$,
 $\left\|\dpy_i(t_i)-\dpy_{i+1}(0)\right\| \leq \delta$, for all $i=1, \ldots, n-1$, and
 $\left\|\dpy_1(0)-\bm{a}\right\| \leq \delta$ and $\left\|\dpy_n(t_n)- \bm{b}\right\| \leq \delta$.
The sequence $\left(\dpy_1, \ldots, \dpy_n\right)$ is called an $(\delta, T)$ chain in $\Lambda$ from $\bm{a}$ to $\bm{b}$  for the differential inclusion \eqref{eq:differential-inclusion}.

A point $\bm{a}$ is chain recurrent if $\bm{a} \hookrightarrow \bm{a}$. The set of chain-recurrent point is denoted $\chainrecurrent(\flow[][\DIhg])$. We now have all the necessary notions to formulate an analogue of \Cref{theo:convergence-sa}.
\begin{theorem}
\label{theo:convergence-di}
Assume \Cref{assum:SA:stepsize,assum:SA:vector-field,assum:SA:boundedness,assum:SA:noise}.
Then $\limset(\dprm)$ is a connected  internally chain-recurrent set for the flow $\flow[][\DIhg]$.
\end{theorem}
The last step consists in linking the internally chain-recurrent sets of the $\flow[][\hg]$, associated to the ODE \eqref{eq:ODE-method}, to those of the set-valued flow $\flow[][\DIhg]$, associated to the differential inclusion \eqref{eq:differential-inclusion}. For this purpose we use a general result, \cite[Theorem~3.1]{benaim2012perturbations} on perturbations of set-valued dynamical systems. This key result is used in the proofs of \cite{ramaswamy2017analysis} and \cite{tadic2017asymptotic}.
\begin{theorem}
\label{theo:convergence-di}
Assume \Cref{assum:SA:stepsize,assum:SA:vector-field,assum:SA:boundedness,assum:SA:noise}. Let $\mathcal{V}$ be an open neighborhood of $\chainrecurrent(\flow[][\hg])$. Then, there exists $\boundbias_0^{\max}$ such that, for all $\boundbias_0 \in \coint{0,\boundbias_0^{\max}}$,
$\chainrecurrent(\flow[][\DIhg]) \subset \mathcal{V}$.
\end{theorem}
In words, this means that the boundary sets of perturbed recursions are in a neighborhood of the boundary sets of unperturbed recursions, if the perturbation is small enough. The weakness of this result is that it is non-quantitative. It can be made more precise, in the case where the mean field is the gradient of a sufficiently regular function; see \cite[Theorem~2.1]{tadic2017asymptotic}.
We finally give a simplified version of the previous result which is based on \Cref{coro:convergence-sa}:
\begin{corollary}
\label{coro:convergence-di}
Assume that $\lyap$ is a strict Lyapunov function for $\flow[][\hg]$. Then, for all open neighborhood $\mathcal{V}$ of the equilibrium set $\{ \hg = 0\}$, there exists $\boundbias_0^{\max}$ such that, for all $\boundbias_0 \in \coint{0,\boundbias_0^{\max}}$, $\limset(\dprm) \subset \mathcal{V}$.
\end{corollary}

\section{Variance Reduction}\label{sec:vr}
Lastly we review on a recent advance in {\SA}, namely the \emph{variance reduction} technique. 
These results are motivated by relevant applications in ML and SP, and involve slight modifications to the basic {\SA} scheme \eqref{eq:sa}. Below, we shall introduce the main algorithm and the general theoretical results. 
The proofs are postponed to \Cref{sec:RV:proofslemma}.

\label{subsec:variance-reduction}
Variance reduction in {\SA} aims to provide a sequence of iterates
$\sequence{\prm}[k][\nset]$ having smaller variance than a plain {\SA}
scheme. We describe a general variance reduction
technique for non-gradient {\SA} when the mean field $\hg$ is a finite
sum
\begin{equation}\label{eq:VR:objective}
\textstyle \hg(\prm) \eqdef \frac{1}{n} \sum_{i=1}^n \hg_i(\prm), \quad \prm \in \rset^d.
\end{equation}
Originally, variance reduction techniques for {\SA} were proposed in the
stochastic gradient setting
(see the survey paper \cite{gower2020variance} and the classical references \cite{johnson2013accelerating,defazio2014saga,wang:etal:2017,nguyen:etal:2017,fang2018spider,wang:etal:2019}, see also \cite{shang:etal:2020,zhang:etal:2022,han:gao:2022,luo:etal:2022,twyman:etal:2023}).
These results were later  extended to non-gradient {\SA}, in a series of works targeting mostly the stochastic versions of the EM algorithm \cite{chen2018stochastic,karimi:etal:2019,fort2020stochastic,fort2021fast,fort2021geom,fort:moulines:2021,fort:moulines:2021b,dieuleveut2021federated,fort:moulines:2022}.

Variance reductions techniques are based on the use of
\emph{control variates} (see, for example. \cite[chapter 5]{asmussen:glynn:2007}). Given an unbiased estimator ${\bf U}$ of the unknown quantity $\PE [{\bf U}]$, a control variate is a centered random variable $\contvarem$, such that the variance of ${\bf U} +
\contvarem$ is less than the variance of ${\bf U}$; such a variate
yields an estimate ${\bf U}+ \contvarem$ of $\PE [{\bf U}]$ with a smaller variance than the original estimator ${\bf U}$. The difficulty is to design such a variable $\contvarem$ with minimal additional computational effort / memory footprint. The mechanism proposed in this section relies on the \emph{Stochastic Path-Integrated Differential EstimatoR} (SPIDER) proposed by \cite{nguyen:etal:2017} and later improved in \cite{fang2018spider,wang:etal:2019}. It has been shown that the SPIDER method offers  optimality guarantees  among other variance reduction methods in both the stochastic gradient setting and the stochastic EM setting (see e.g. \cite[Table~1]{wang:etal:2019} and \cite{arjevani2022lower} for the stochastic gradient case and \cite[section~6]{fort2020stochastic} for stochastic EM).

At iteration $(k+1)$, SPIDER defines an oracle for $\hg(\prm_k)$ by using the current estimate $\Hgrv_{k}$ of
$\hg(\prm_{k-1})$ as follows
\[ \textstyle 
\Hgrv_{k+1} \eqdef \frac{1}{\lbatchrv} \sum_{i \in \batchrv_{k+1}} \hg_i(\prm_k) + \contvarem_{k+1}
\]
where
\begin{equation}\label{eq:contvar:spider}
\textstyle  \contvarem_{k+1} \eqdef \Hgrv_{k} - \frac{1}{\lbatchrv} \sum_{i \in \batchrv_{k+1}} \hg_i(\prm_{k-1}),
\end{equation}
and $\batchrv_{k+1}$ is a set of indices of cardinality $\lbatchrv$ chosen randomly in $\{1, \ldots, n\}$ with or without replacement and independently of the history of the algorithm $\mcf_k$.  $\contvarem_{k+1}$ is a random variable chosen to be correlated with the naive oracle $\lbatchrv^{-1} \sum_{i \in \batchrv_{k+1}} \hg_i(\prm_k)$, the correlation relying essentially on the random minibatch $\batchrv_{k+1}$.  Since
$\CPE{\lbatchrv^{-1} \sum_{i \in \batchrv_{k+1}} \hg_i(\prm_{k-1})}{\mcf_k} = \hg(\prm_{k-1})$,
\eqref{eq:contvar:spider} shows that $\contvarem_{k+1}$ is the difference of two estimators of $\hg(\prm_{k-1})$; yet its
(conditional) expected value is not zero. Namely, we have
\[
\CPE{\contvarem_{k+1}}{\mcf_k} =  \Hgrv_{k} - \hg(\prm_{k-1});
\]
(see \Cref{lem:VR:bias}).
The \emph{bias} is tamed by resetting  the control variate:
every $\kin$ iterations, the control variate is set to zero.

\begin{algorithm}[t]
  \KwData{ an initial value $\prm_{\mathrm{init}}$, the number of inner loops $\kin$ and outer loops $\kout$, a stepsize sequence $\step_{t,k+1}$ for $t=1, \ldots, \kout$ and $k=0,\ldots, \kin-1$}
\KwResult{A $\rset^d$-valued sequence $\prm_{t,k+1}$, $t=1, \ldots, \kout$ and $k=0,\ldots, \kin-1$.}
$\prm_{0,\kin} = \prm_{\mathrm{init}}$; \\
\For{$t=1, \ldots, \kout$}{$\prm_{t,0} = \prm_{t-1,\kin}$ and $\prm_{t,-1}= \prm_{t-1,\kin}$ \label{line:algo:SPIDERwithinSA:initinner} \;
Sample $\batchrv_{t,0}$, of size $\lbatchrv$, in $\{1, \ldots, n\}$ with or without replacement \;
Set $\Hgrv_{t,0} = \hg(\prm_{t,0})$  \label{line:algo:SPIDERwithinSA:initoracle} \;
\For{$k=0, \ldots, \kin -1$}{Sample $\batchrv_{t,k+1}$, of size $\lbatchrv$, in $\{1, \ldots, n\}$ with or without replacement \;
$\Hgrv_{t,k+1} = \Hgrv_{t,k} +\lbatchrv^{-1} \sum_{i \in \batchrv_{t,k+1}} \left\{ \hg_i(\prm_{t,k}) - \hg_i(\prm_{t,k-1}) \right\}$ \;
 $\prm_{t,k+1} = \prm_{t,k} + \step_{t,k+1} \Hgrv_{t,k+1}$;}}
  \caption{SA-SPIDER \label{algo:SPIDERwithinSA}}
  \end{algorithm}

The \emph{SA-SPIDER} algorithm is given by \Cref{algo:SPIDERwithinSA}. The iteration index is a pair $(t,k)$, where $t$ counts the number of control variate updates (\emph{outer loops} number) and $k$ is the number of {\SA} updates since last reset (\emph{inner loops} number). 

The following intermediate results show how the bias of the stochastic mean field
$\Hgrv_{\cdot}$, its conditional variance and the $L^2$-moment of the error $\Hgrv_{t,k+1} - \hg(\prm_{t,k})$ evolve along the inner
loops. Define the filtration
\[
\mcf_{t,k} \eqdef  \sigma\left( \prm_{\mathrm{init}}, \batchrv_{\tau,\kappa}, 1 \leq \tau \leq t, 1 \leq \kappa \leq k \right),
\]
for all $t \in \{1, \ldots, \kout\}$ and $k \in \{1, \ldots, \kin \}$. Assume
\begin{assumVR} \label{assumVR:lipschitz:hi}
For all $i \in \{1, \ldots, n\}$, there exists $\Liphrv_i \geq 0$ such that for all $\prm, \prm' \in \rset^d$,
$\| \hg_i(\prm) - \hg_i(\prm') \| \leq \Liphrv_i \| \prm - \prm' \|$.
\end{assumVR}
\begin{lemma}\label{lem:VR:bias}
Consider the iterates from \Cref{algo:SPIDERwithinSA}. For any $t \in \{1, \ldots, \kout\}$ and $k \in \{0, \ldots, \kin-1\}$, it holds
  \begin{align}
  \CPE{\Hgrv_{t,k+1}}{\mcf_{t,k}} - \hg(\prm_{t,k})     &= \Hgrv_{t,k} - \hg(\prm_{t,k-1}) \label{item:lem:VR:bias:1} \\
  \CPE{\Hgrv_{t,k+1}- \hg(\prm_{t,k})}{\mcf_{t,0}}     &= \Hgrv_{t,0} - \hg(\prm_{t,-1}) = 0. \label{item:lem:VR:bias:2}
\end{align}
\end{lemma}
Equation \eqref{item:lem:VR:bias:1} shows that at each inner iteration, the oracle $\Hgrv_{t,k+1}$ is a biased approximation of the mean field $\hg(\prm_{t,k})$ and the bias propagates  along inner iterations.  Conditionally to the initialization of the inner iteration, the bias is equal to the error $\Hgrv_{t,0} - \hg(\prm_{t,-1})$ (see \eqref{item:lem:VR:bias:2}). The strategy used in \Cref{line:algo:SPIDERwithinSA:initinner} and \Cref{line:algo:SPIDERwithinSA:initoracle}  of \Cref{algo:SPIDERwithinSA} implies that $\Hgrv_{t,0} - \hg(\prm_{t,-1})=0$. In that sense, we say that at the beginning of each inner loop, the bias is canceled  and the control variate is reset. 

\begin{lemma}\label{lem:VR:varL2error}
Assume \Cref{assumVR:lipschitz:hi} and set  $\Liphrv^2 \eqdef n^{-1} \sum_{i=1}^n \Liphrv_i^2$.  For any
 $t \in \{1, \ldots, \kout\}$ and $k \in \{0, \ldots, \kin-1\}$, it
 holds 
\[
 \CPE{ \|\Hgrv_{t,k+1}
 - \CPE{\Hgrv_{t,k+1}}{\mcf_{t,k}} \|^2}{\mcf_{t,k}} \leq \frac{\Liphrv^2}{\lbatchrv} \step_{t,k}^2 \|\Hgrv_{t,k}\|^2, \\
\]
and with the convention that $\step_{t,0} \eqdef 0$, it holds
\begin{multline*}
\CPE{ \| \Hgrv_{t,k+1} - \hg(\prm_{t,k}) \|^2 }{\mcf_{t,k}} \leq \| \Hgrv_{t,k} - \hg(\prm_{t,k-1}) \|^2   \\ + ( {\Liphrv^2} / {\lbatchrv} ) \step^2_{t,k} \| \Hgrv_{t,k}\|^2.
\end{multline*}
\end{lemma}
\noindent When $\step_{t,k+1} \leq \step_{t,k}$ for any $k \geq 1$,  \Cref{lem:VR:varL2error} implies that 
\begin{align} \label{eq:VR:controlMSE}
 &\step_{t,k+1}   \CPE{ \| \Hgrv_{t,k+1} - \hg(\prm_{t,k}) \|^2 }{\mcf_{t,0}}   \\  & \leq  \frac{2 \Liphrv^2}{\lbatchrv} \sum_{\ell=1}^{k} \step^3_{t,\ell}  \CPE{ \| \Hgrv_{t,\ell}  - \hg(\prm_{t,\ell-1})  \|^2 
 +  \| \hg(\prm_{t,\ell-1})  \|^2 }{\mcf_{t,0}}. \nonumber
\end{align}
In the simple case when  the stepsize sequence is constant $(\step_{t,k} = \step)$ and the mean field is bounded $(\clyap_{\hg,1} =0)$, the summed MSE along the inner loop iterations satisfies
\begin{multline*}
\left(1 - 2  \step^2 \frac{\Liphrv^2 \kin}{\lbatchrv}  \right) \sum_{k=0}^{\kin-1} \CPE{ \| \Hgrv_{t,k+1} - \hg(\prm_{t,k}) \|^2 }{\mcf_{t,0}}   \\ \leq  2  \step^2 ({\Liphrv^2 \kin^2} / {\lbatchrv})    \clyap_{\hg,0}.
\end{multline*}
This inequality illustrates the benefit of variance reduction: the cumulated MSE can be set arbitrarily small by a convenient choice of the learning rate $\step$. Such a property remains true when the stepsize sequence is not constant and the mean field $\hg$ is not bounded; it will be a key ingredient for the convergence analysis provided in \Cref{theo:VR:nonasymptotic} below.

The following lemma is an analogue of the Robbins-Siegmund Lemma for obtaining non-asymptotic bounds.
\begin{lemma}\label{lem:VR:robbins-siegmund-ameliore}
Assume \Cref{assum:field}-\ref{item:clyap-hg} and \Cref{assum:lyapunov}.
For any  $t \in \{1, \ldots, \kout\}$ and $k \in \{0, \ldots, \kin-1\}$, it holds
\begin{align*}
&\nonumber \CPE{\lyap(\prm_{t,k+1})}{\mcf_{t,0}}  \leq  \CPE{\lyap(\prm_{t,k})}{\mcf_{t,0}}   \\
&\nonumber - \step_{t,k+1} \{ \nu_{t,k+1}  \CPE{\superlyap(\prm_{t,k})}{\mcf_{t,0}}  \\
& \qquad \qquad + \mu_{t,k+1}  \,  \CPE{\| \Hgrv_{t,k+1} - \hg(\prm_{t,k})\|^2}{\mcf_{t,0}} \} \nonumber \\
& + \step_{t,k+1} \operatorname{a} \,  \CPE{\| \Hgrv_{t,k+1} - \hg(\prm_{t,k})\|^2}{\mcf_{t,0}} + \step^2_{t,k+1} \Liplyap \clyap_{\hg,0}.
\end{align*}
where  $\nu_{t,k+1}  \eqdef   \rholyap/2 - \clyap_{\hg,1} \step_{t,k+1} \Liplyap$, $\mu_{t,k+1}  \eqdef \rholyap/2 - \step_{t,k+1} \Liplyap$   and $\operatorname{a}  \eqdef  (\clyap_\lyap^2 \rholyap^{-1} + \rholyap)/2$.
\end{lemma}
We will use \eqref{eq:VR:controlMSE} to show that the term $\CPE{\| \Hgrv_{t,k+1} - \hg(\prm_{t,k})\|^2}{\mcf_{t,0}}$ in the \rhs\ is negligible \wrt the term $\mu_{t,k+1}  \,  \CPE{\| \Hgrv_{t,k+1} - \hg(\prm_{t,k})\|^2}{\mcf_{t,0}}$ in the LHS. 
We can  establish a non-asymptotic
convergence result.
\begin{theorem} \label{theo:VR:nonasymptotic}
Assume \Cref{assum:field}-\ref{item:clyap-hg}, \Cref{assum:lyapunov}   and \Cref{assumVR:lipschitz:hi}. Let $\kin$,
$\kout$ be positive integers and $\prm_{\mathrm{init}} \in \rset^d$.
Let $\{\step_{t,k+1}, 1 \leq t \leq \kout, 0 \leq k \leq \kin-1\}$ be a stepsize sequence satisfying  for all $k \geq 1$: $\step_{t,k+1} \leq \step_{t,k}$,
$(1 \vee \clyap_{\hg,1}) \, \step_{t,k} \, \lambda_{t,k} \in \ooint{0, \rholyap/2}$
where $ \lambda_{t,k}\eqdef \Liplyap + \step_{t,k} \Liphrv^2     \left\{ \clyap_\lyap^2 \rholyap^{-1} +  \rholyap  \right\}\kin /\lbatchrv$.  
Consider the sequence given by \Cref{algo:SPIDERwithinSA}.  Then
\begin{align}
& \textstyle \sum_{t=1}^\kout \sum_{k=1}^\kin \step_{t,k} \left(\frac{\rholyap}{2} - \clyap_{\hg,1} \, \step_{t,k}  \, \lambda_{t,k} \right) \PE\left[\superlyap(\prm_{t,k-1}) \right] \nonumber \\
& \textstyle + \sum_{t=1}^\kout \sum_{k=1}^\kin \step_{t,k} \left(\frac{\rholyap}{2} -  \step_{t,k}  \, \lambda_{t,k} \right) \PE\left[\|\Hgrv_{t,k} - \hg(\prm_{t,k-1}) \|^2 \right] \nonumber \\
&\leq \PE\left[\lyap(\prm_{\mathrm{init}})\right] - \lyap_\star  +  \clyap_{\hg,0}  \operatorname{B}^{\mathrm{vr}}, \label{eq:SPIDER-SA-bound}
\end{align}
where 
\[
\operatorname{B}^{\mathrm{vr}} \eqdef \Liplyap  \sum_{t=1}^\kout \sum_{k=1}^{\kin}\step^2_{t,k} +
 \frac{\Liphrv^2 \kin }{\lbatchrv}
 \left(\frac{\clyap_\lyap^2}{\rholyap} + \rholyap \right) \sum_{t=1}^\kout \sum_{k=1}^{\kin} \step^3_{t,k}.
\]
\end{theorem}
\Cref{theo:VR:nonasymptotic} controls $\superlyap(\prm_{t,k})$ and the quadratic error $\| \Hgrv_{t,k+1} - \hg(\prm_{t,k})\|^2$ along the $\kin \kout$ iterations. Set
\[
\totstep \eqdef \kin \kout. 
\]
We observe the following consequences of \Cref{theo:VR:nonasymptotic}.\vspace{.1cm}

\noindent \textbf{Random stopping.}
The LHS in \Cref{theo:VR:nonasymptotic} can be viewed as the expected value of $\superlyap(\prm_{R_{\totstep }^\mathrm{vr}})$ and $\|\Hgrv_{\tilde R_{\totstep }^\mathrm{vr}+1} - \hg(\prm_{\tilde R_{\totstep }^\mathrm{vr}}) \|^2 $  where $R_{\totstep }^\mathrm{vr}$ and $\tilde R_{\totstep }^\mathrm{vr}$ are random variables taking values in $\{1, \ldots, \kout\} \times \{0, \ldots, \kin-1\}$, independent of the $\sequence{\prm}[k][\nset]$, and with probability mass functions (see \Cref{subsec:finite-time-bounds})
\begin{align*} 
& \PP ( R_{\totstep }^{\mathrm{vr}} = (t,k)) \propto  \step_{k+1} \left( {\rholyap} / {2} - \clyap_{\hg,1} \, \step_{t,k+1}  \, \lambda_{t,k+1} \right) \\
& \PP ( \tilde R_{\totstep }^{\mathrm{vr}} = (t,k)) \propto \step_{k+1}\left( {\rholyap} / {2} -  \step_{t,k+1}  \, \lambda_{t,k+1} \right).
\end{align*}

\begin{remark}
\label{rem:VR:convex-superlyap} Similar to \Cref{rem:convex-superlyap}. When the function $\superlyap$ is convex, a stronger convergence result can be derived. Define the convex combination of the iterates
\[
\bar{\prm}_{\totstep }^{\mathrm{vr}} \eqdef  \sum_{t=1}^{\kout} \sum_{k=0}^{\kin-1} \frac{ \step_{t,k+1} \, \omega_{t,k+1}^{\mathrm{vr}}}{\sum_{t'=1}^{\kout} \sum_{k'=0}^{\kin-1} \step_{t',k'+1} \, \omega_{t',k'+1}^{\mathrm{vr}}}\eqsp,
\]
where $\omega_{t,k+1}^{\mathrm{vr}} \eqdef \rholyap/2 -    \clyap_{\hg,1} \, \step_{t,k+1} \, \lambda_{t,k+1}$. We have from \Cref{theo:VR:nonasymptotic} that $\PE\left[ \superlyap(\bar{\prm}_{\totstep }^{\mathrm{vr}})\right]$ is upper bounded by the RHS of \eqref{eq:SPIDER-SA-bound} divided by $\sum_{t'=1}^{\kout} \sum_{k'=0}^{\kin-1} \step_{t',k'+1} \, \omega_{t',k'+1}^{\mathrm{vr}}$.
\end{remark}

\noindent \textbf{Constant stepsize.} Assume $\step_{t,k+1}=\step$ for each $t \in \{1, \ldots, \kout\}$ and $k \in \{0, \dots, \kin-1\}$.  The assumptions of 
\Cref{theo:VR:nonasymptotic} are satisfied by choosing $\step \in \ooint{0, \step_{\max}^{\mathrm{vr}}}$ where 
\[
\step_{\max}^{\mathrm{vr}} \eqdef \frac{\rholyap \Liplyap \lbatchrv}{2 L^2 (\clyap_\lyap^2+ \rholyap^2)\kin}\left(\{1+ 2 \frac{L^2 (\clyap_\lyap^2+ \rholyap^2) \kin}{\Liplyap^2 (1 \vee \clyap_{\hg,1})\lbatchrv } \}^{1/2} - 1\right).
\]
This yields
\begin{align*}
\label{eq:main-result-non-asymptotic-constant}
& \textstyle \frac{1}{\totstep } \sum_{t=1}^{\kout}\sum_{k=0}^{\kin-1} \left( \PE [ \superlyap(\prm_{t,k}) ]  + \PE\left[ \| \Hgrv_{t,k+1} - \hg(\prm_{t,k} \|^2 \right]\right)  \\
& \leq \frac{\Delta_1 }{ \step \totstep  \{\rholyap/2 - \step \lambda(\step) (1 \vee \clyap_{\hg,1})\} }  \\
& \quad + \step \, \clyap_{\hg,0}   \left( \frac{\Liplyap  + \step_{\max}^{\mathrm{vr}} \Liphrv^2 
 \left(\clyap_\lyap^2 + \rholyap^2 \right)  \kin/(\rholyap \lbatchrv)}{\rholyap/2 - \step \lambda(\step) (1 \vee \clyap_{\hg,1})} \right) \eqsp,
\end{align*}
where
\begin{align*}
\lambda(\step) & \eqdef  \Liplyap +   \step \Liphrv^2     \left\{ \clyap_\lyap^2 \rholyap^{-1} +  \rholyap  \right\} \kin /\lbatchrv, \\
\Delta_1 & \eqdef   \PE[\lyap(\prm_{\mathrm{init}})] - \lyap_\star.
\end{align*}
Contrary to \Cref{theo:main-result-quantitative} (see \Cref{coro:main-result-quantitative}), the \rhs\ can be made small by a clever choice of $\step$: even if the oracle $\Hgrv_{t,k+1}$ is a biased estimator of $\hg(\prm_{t,k})$ (see \Cref{lem:VR:bias}), the SPIDER variance reduction is able to manage this bias.   

The terms in the \rhs\  can be adjusted according to the total number of iterations $\totstep$ and the different parameters of the problem by a suitable choice of $\step$. When $\step \leq \step_{\max}^{\mathrm{vr}}/2$, then $\rholyap/2 - \step \lambda(\step) (1 \vee \clyap_{\hg,1}) \geq \rholyap/4$; in addition, the function  $\gamma \mapsto  \alpha / \gamma + \beta \gamma$ is minimized on $\ocint{0,\step_{\max}^{\mathrm{vr}}/2}$  at the point $\sqrt{\alpha/\beta} \wedge (\step_{\max}^{\mathrm{vr}}/2)$ when $\alpha, \beta >0$. Therefore, set
\[
\Delta_2 \eqdef 2\Liplyap \rholyap  + \step_{\max}^{\mathrm{vr}} \Liphrv^2 
 \left(\clyap_\lyap^2 + \rholyap^2 \right)  \kin/\lbatchrv ,
\]
and 
$\step_{\totstep }^{\mathrm{vr}}  \eqdef  \step_{\max}^{\mathrm{vr}}/2$ if $\clyap_{\hg,0} =0$ and otherwise
\begin{equation}
\label{eq:VR:optimal-step}
\step_{\totstep }^{\mathrm{vr}} \eqdef  \sqrt{2\Delta_1 \rholyap  / \totstep \Delta_2}  \wedge ( {\step_{\max}^{\mathrm{vr}}} / {2} ).
\end{equation}
\begin{corollary}[of \Cref{theo:VR:nonasymptotic}]
\label{coro:VR:main-result-quantitative}
Setting  $\step_{t,k+1}= \step_{\totstep}^{\mathrm{vr}}$ for $t \in \{1, \ldots, \kout\}$ and $k \in \{0,\ldots,\kin-1\}$ we get
\begin{multline*}
\frac{1}{\totstep }  \sum_{t=1}^{\kout}\sum_{k=0}^{\kin-1} \left( \PE [ \superlyap(\prm_{t,k}) ]  + \PE\left[ \| \Hgrv_{t,k+1} - \hg(\prm_{t,k} \|^2 \right]\right)  \\
\leq   {4 \Delta_1 } / ({ \step_{\totstep}^{\mathrm{vr}} \totstep  \rholyap})  + 2 \step_{\totstep}^{\mathrm{vr}} \, \clyap_{\hg,0}  {\Delta_2} / {\rholyap^2}\eqsp.
\end{multline*}
\end{corollary}
Choose $\kin =\lbatchrv$ and let us comment the rate of the \rhs\ when $\totstep \to +\infty$. When $\clyap_{\hg,0} = 0$,  the \rhs\ decreases at the rate $O(1/\totstep)$, while when $\clyap_{\hg,0} >0$, the rate is $O(1/\sqrt{\totstep})$.

\vspace{.1cm}
\noindent \textbf{$\epsilon$-Approximate Stationarity.}
Using \Cref{coro:VR:main-result-quantitative}, we analyze the complexity of SA-SPIDER to reach $\epsilon$-approximate stationarity.  First, observe that  the \rhs\ in
\Cref{coro:VR:main-result-quantitative} is an upper bound of  $\PE\left[\superlyap(\prm_{R_{\totstep }^\mathrm{vr}})\right]$ where $R_{\totstep }^\mathrm{vr}$ is a uniform random variable on $\{1, \ldots, \kout\} \times \{0, \ldots, \kin-1\}$.

Consider first the case when  $\clyap_{\hg,0} =0$; then $\step_{\totstep}^{\mathrm{vr}} = \step_{\max}^{\mathrm{vr}}/2$ and is independent of the accuracy $\epsilon$. Choose $\kout = {\cal O}(1/(\sqrt{n} \epsilon))$ and $\kin = \lbatchrv = \lceil\sqrt{n}\rceil$; it means that the total number of calls to one of the functions $\hg_i$ (see \eqref{eq:VR:objective}) is equal to $n$. Then  the total number of iterations to reach an $\epsilon$-approximate stationary point  is ${\cal O}(1/\epsilon)$ and the total number of calls to one of the functions $\hg_i$  is $n \kout + 2 \kout \kin \lbatchrv = {\cal O}(\sqrt{n}/\epsilon)$. Such a complexity analysis retrieves earlier results established in specific settings of SA-SPIDER:  {\SG} for non convex optimization \cite[Theorem 2]{fang2018spider}, \cite[Theorem 1]{wang:etal:2019}, the stochastic EM algorithms \cite[Section 3]{fort2020stochastic} and for more general SA-based root-finding problems  \cite[Corollary 4.3]{fort:moulines:2022}. 

Consider now the case  $\clyap_{\hg,0}>0$. Choose again $\kin = \lbatchrv  = \lceil n \rceil$, and  $\kout$ large enough so that 
\[
\totstep \geq \frac{8 (1 +\clyap_{\hg,0})^2 \Delta_1 \Delta_2}{\rholyap^3 \epsilon^2} \vee \frac{8 \Delta_1 \rholyap}{\Delta_2 (\step_{\max}^{\mathrm{vr}})^2}.
\]
Then  $\step_{\totstep}^{\mathrm{vr}} = \sqrt{2 \Delta_1 \rholyap}/\sqrt{\totstep \Delta_2}$ and SA-SPIDER reaches an $\epsilon$-approximate stationary point in $\totstep = {\cal O}(1/\epsilon^2)$ iterations, and by calling ${\cal O}(\sqrt{n}/\epsilon^2)$ functions $\hg_i$. To our best knowledge, this is the first complexity analysis of SA-SPIDER when $\clyap_{\hg,0} >0$.

\section{Conclusions}
We have overview state-of-the-art results for {\SA} scheme with a focus on its use as a general stochastic non-gradient algorithm common to statistical learning. We have proposed 
a general theoretical framework based on the designs of flexible Lyapunov function and unified the modern asymptotic, non-asymptotic convergence results for {\SA} schemes. We illustrated the applications of our techniques to {\SG}, compressed {\SA}, stochastic EM and TD learning; as well as presenting how the recent variance reduction technique can be adopted to {\SA}. We studied the effects of \emph{bias} in {\SA} updates caused by the non-gradient nature in popular designs.  

Our findings shed lights on how to design stochastic algorithms with nice convergence properties. In particular, we illustrated how to construct the stochastic random field in {\SA} from fixed point equation of the statistical learning problem, and to tame with the bias in {\SA} resulted from the design.

\bibliographystyle{IEEEtran}
\bibliography{ref_SA}

\clearpage
\newpage

\appendix
This document presents the proofs or detailed calculations that have been skipped in the main paper. Readers can find the following content: \Cref{subsec:inequality} shows two elementary inequalities for numerical sequences. \Cref{subsec:proof:csgd-checking} and \Cref{subsec:proof:TD-learning-checking} check the assumptions for the examples on compressed {\SA} and TD(0) learning. \Cref{subsec:proof:finite-time-bounds} gives a variant of \Cref{lem:main-result-quantitative} for the finite-time bounds of {\SA}. \Cref{app:proofGS} and \Cref{sec:appendix:SAEM} show how to obtain finite-time bounds for compressed {\SG} methods and SAEM algorithms. \Cref{app:proof-ODE} gives missing proof for the asymptotic convergence. \Cref{sec:RV:proofslemma} provides proofs for the variance reduced {\SA} algorithm.

\subsection{Table of notations}
\revisionupdates{In the following, we provide complete notation tables, that aggregate the notations of the paper, with references to the points where each notation is introduced. 
\Cref{tab:all_not_SGD,tab:all_not_CSA,tab:all_not_EM,tab:all_not_TD} respectively aggregate all notations used in the  application to the four examples, namely SGD, compressed \SA, EM, and TD learning.}

\ \\ \ \\


\begin{table}[h]
    \caption{Summary of notations used in the analysis of SGD, in \Cref{subsec:SGDintro,subsec:SGD-checking,subsec:SGD}.}
    \label{tab:all_not_SGD}
\resizebox{\linewidth}{!}{\begin{tabular}{p{0.4\linewidth}p{0.8\linewidth}p{0.25\linewidth}}\toprule
Notation & Object & Def.~in   \\
\midrule
&&{\Cref{subsec:SGDintro}}\\
\midrule
$F$ & function to be minimized &   \\
$\nabla F(\cdot)$ &gradient of  $F$ and mean-field $\hg$  &   \\
$n$ & number of observations &   \\
$(Z_1, \dots, Z_n)$ & observations &   \\
$\rho$ &  distribution of $\State_{k+1}$ &   \\
$\loss$ & loss & \\
$\lbatchgd$ & batch  size &  \\
$f_i \eqdef \ell (\prm, Z_i)$ & loss function on obs. $i$ & \\
\midrule
&&{\Cref{subsec:SGD-checking}}\\
\midrule
 $\Lip{\nabla f_i}$ & Smoothness of $f_i$ & \Cref{{assum:SGD}}.\ref{item:smoothF}\\
 $\Lip{\nabla f}$ & Smoothness of $f$ &  \\
$M$ & Uniform bound on $\| \nabla f_i (\prm) - \nabla F(\prm)\|$ & \Cref{{assum:SGD}}.\ref{item:SGD:boundedF}\\
$\mu >0$ & strong-convexity  modulus of $F$ & \Cref{assum:SGD:strongly-convex}\\
\midrule
\bottomrule
\end{tabular}
}
\end{table}


\begin{table}[h]
    \caption{Summary of notations used in the analysis of Compressed and modified \SA, in \Cref{subsec:ex-variantsSGD,subsec:csgd-checking,subsub:exrate:CSGD}.}
    \label{tab:all_not_CSA}
\resizebox{\linewidth}{!}{\begin{tabular}{p{0.4\linewidth}p{0.8\linewidth}p{0.25\linewidth}}
\toprule
Notation & Object & Def.~in   \\
\midrule
&&{\Cref{subsec:ex-variantsSGD}}\\
\midrule
$F$ & function to be minimized &   \\
$j_{k+1} \in \{1, \ldots, d\}$ & chosen coordinate in the $k$-th  iteration &   \cref{eq:gauss_south} \\
$\{\bm e_1, \ldots, \bm e_d\}$ & the canonical basis of $\rset^d$ &    \\
$\nabla_j F$ & $j$-th coordinate of the gradient  &     \\
$\compressor: \rset^d \times \sf U \to \rset^d$& Compression operator &   \\
$\sf U$ & general state space  &  \\
$\mu_U$ & distribution of $\sf U$ &  \\
$\rm{Rand}_h, \rm{Top}_h$ & Ex. of sparsification-based compression operator & \Cref{rem:expl-compr}.1 \\
${\rm h} $ & Sparsification parameter for $\rm{Rand}_h, \rm{Top}_h$  & \Cref{rem:expl-compr}.1  \\
$Q_d,  Q_s$ & Ex. of quantization-based compression operator &   \Cref{rem:expl-compr}.2 \\
$\Delta$ & quantization resolution &   \Cref{rem:expl-compr}.2\\
\midrule
&&{\Cref{subsec:csgd-checking}}\\
\midrule
$\omgbiased $ & Contractivness of $\mathcal C$ & \Cref{assum:SGD:contractive}, \cref{eq:comp_contractive} \\
$\zeta_1, \zeta_2 \in \bar \rset_+$ & technical constants (from Young inequality) & \Cref{lem:csgd-var1}\\ 
$\omgunbiased $ & Relative bound on variance of $\mathcal C$ & \Cref{assum:SGD:URBV}, \cref{eq:URVB} \\
$\omguniform $ & Uniform bound on variance of $\mathcal C$ & \Cref{assum:SGD:UnifBoundedVar}, \cref{eq:unif_bounded_C} \\
\bottomrule
\end{tabular}
}
\end{table}
\ \\ \ \\


\begin{table}[h]
    \caption{Summary of notations used in the analysis of EM, in \Cref{sec:EM:introduction,sec:EM:verifH,subsec:stochastic-EM}.}
    \label{tab:all_not_EM}
    \resizebox{\linewidth}{!}{\begin{tabular}{p{0.5\linewidth}p{0.75\linewidth}p{0.2\linewidth}}\toprule
Notation & Object & Def.~in   \\
\midrule
&& \Cref{sec:EM:introduction}\\
\midrule
$F(\prmo) \eqdef - \log \int_{\widetilde \Zset} p(z; \prmo) \tilde \mu(\rmd z)$ & Intractable objective function (to be minimized) & \cref{eq:EM:objective:prmo-2} \\
$\prmo$&  parameter of the original optimization problem & \\
$p(z; \prmo)=\prod_{i=1}^{n} p_i(z_i;\prmo)$& Product form; $z\in \widetilde{\Zset}$ and $z_i \in \Zset$ & \\
$\mu$ &  sigma-finite measure on the measurable set $\Zset$ & \\
$\pi_i(z_i; \prmo) \eqdef p_i(z_i;\prmo)/ \int_\Zset p_i(u; \prmo) \mu(\rmd u)$ &  A distribution on $\Zset$  & \\
$\QEM[\prmo'] $ & EM surrogate function tangent at $\prmo'$ & \cref{eq:EM:iter:theta}, \eqref{eq:EM:QEM} \\
\midrule
\multicolumn{3}{c}{-- In the maximum likelihood  context --}\\
\midrule
$Y_i$ & observation $\# i$ & \\
$z_i$  & latent variable  $\# i$& \\
$p_i(z_i; \prmo)$
& joint probability of the observation $Y_i$ and the latent variable $z_i$ for a given value of the parameter $\prmo$; & \\
$ g_i(\prmo) \eqdef   \int_\Zset p_i(z_i; \prmo) \, \mu(\rmd z_i) $&   likelihood of $Y_i$  &    \cref{eq:EM:likelihoodY} \\
$ \pi_i(z_i; \prmo)$ & posterior distribution of the latent variable $z_i$ given the observation $\# i$ when the value of the parameter is $\prmo$ & \cref{eq:EM:posterior} \\
\midrule
\multicolumn{3}{c}{In the exponential family context}\\
\midrule
$\sem_i$ & {sufficient statistics} associated to the observation $\# i$ & \Cref{assumEM:expo} \\
$\barsem_i(\prmo) \eqdef \int_\Zset \sem_i(z_i) \pi_i(z_i; \prmo) \mu(\rmd z_i)$&  expectation of the sufficient statistic & \cref{eq:EM:def:bars}\\ 
$\barsem(\prmo) \eqdef n^{-1} \sum_{i=1}^n \barsem_i(\prmo)$ & the mean value of the $n$ functions $\barsem_i$ & \\
$\QEM[\prmo'](\cdot) \eqdef \pscal{\barsem(\prmo')}{\phiem(\cdot)} - \psiem(\cdot)$ & the specific form of the EM surrogate function  & \\
$\phiem,  \psiem$ & functions  on $\rset^d$, parameterizing  the family  of surrogate functions& \Cref{assumEM:expo} \\
$\mapem $ & optimization map  in EM& \Cref{assumEM:Mstep} \\
$\lbatchem$&  batch size in Mini-batch EM  & \cref{eq:onlineEM:oracle}\\
$\rm m$ &  number of Monte Carlo samples  in SAEM& \cref{eq:oracle-SAEM}, \eqref{eq:self-normalized} \\
\midrule
&& \Cref{sec:EM:verifH}\\
\midrule
$F_\star$ & Uniform lower  bound on $F$& \Cref{assumEM:supplementary}\\
 $\Bem(\prm)$ &$d \times d$ p.d.~matrix, s.t.~$\nabla \lyap(\prm) = - \Bem(\prm) \hg(\prm)$ & \Cref{assumEM:supplementary}\\
$\vminem \leq \vmaxem$ & constants characterizing the conditioning of $\Bem$ & \Cref{assumEM:supplementary}\\
$ \boundvarem_{0}, \boundvarem_{1} \in \rset_+$ & control on the averaged sufficient statistic & \Cref{assumEM:var:onlineEM,assumEM:var:SAEM:exact}\\
$ s_\star$ & Uniform bound on $\| \sem_i(z) \|$  & \Cref{assumEM:var:SAEM:IS}\\
\bottomrule
\end{tabular}
}
\end{table}

\ \\ \ \\


\begin{table}[h]
    \caption{Summary of notations used in the analysis of TD learning, in \Cref{sec:definition-TD-learning,subsec:TD-learning-checking,subsec:TD-learning}.}
    \label{tab:all_not_TD}
\resizebox{\linewidth}{!}{\begin{tabular}{p{0.4\linewidth}p{0.8\linewidth}p{0.25\linewidth}}\toprule
Notation & Object & Def.~in   \\
\midrule
&& \Cref{sec:definition-TD-learning}\\ 
\midrule
 $\policy$ & policy in a Markov Decision Process&   \\ 
 $( \stateMRP, \kerMRP, \rewardMRP, \lambda )$ & Markov Reward Process (MRP)  &   \\ 
  $\stateMRP= \{s_1,\dots,s_n\}$ &  state-space   &   \\ 
$\kerMRP$ & $n \times n$ state transition matrix of the probability of transition &   \\ 
$\rewardMRP(s, s')$  & reward function  &   \\ 
$\lambda \in \ooint{0,1}$ &  discount factor  &   \\ 
$\rewardMRP(s)$ & expected instantaneous reward from state $s$ &   \\ 
 $\statdistMRP$& unique stationary distribution&    \Cref{assum:TD:stationary-policy} \\ 
 $\valuefunc$ &  value function  of the MRP &     \cref{eq:valuefunc_def}\\ 
 $\sequence{S}[k][\nset]$ & Markov chain  started at $S_0=s$, with Markov kernel $\kerMRP$  &   \\ 
$\valuefunc[\prm](s) \eqdef \feature(s)^{\top} \prm$  & linear approximation of $\valuefunc$ & \\ 
$\feature(s) \in \rset^d$ & \emph{feature vector} for the state $s \in \stateMRP$ &  \\ 
 $\prm \in \rset^d$ & parameter vector to be estimated &  \\ 
 $\Feature$ & $n \times d$ feature matrix&  \\ 
\midrule
&&\Cref{subsec:TD-learning-checking}\\
\midrule
$\lambda \in (0,1)$ &  contraction modulus &\Cref{lem:contraction-projected-bellman} \\
$\valuefunc[\star]$ &  unique value function  in  $\Span(\Feature)$ which solves the fixed point of the projected Bellman eq. &\Cref{lem:contraction-projected-bellman} \\
$\vmintd$  & minimal eigenvalue of $\Feature^\top \Dtd_{\statdistMRP} \Feature$ & \Cref{lem:TD:upper-lower-spectrum-covariance}\\
$\bSigma_{\statdistMRP}$ & {feature covariance matrix}&\Cref{eq:covariance-feature}\\
\bottomrule 
\end{tabular}
}
\end{table}

\subsection{Elementary inequalities} \label{subsec:inequality}
\begin{lemma}
\label{lem:summation}
Let $a>0$ and $\sequence{\step}[k][\nset]$ be a sequence such that $\step_k < 1/a$ for any $k \geq 1$. Then, for any integer $k \geq 1$,
\begin{equation}
\label{eq:sum}
\sum_{j=1}^{k+1} \step_j \prod_{l=j+1}^{k+1}\left(1-\step_l a\right)=\frac{1}{a}\left\{1-\prod_{l=1}^{k+1}\left(1-\step_l a\right)\right\}
\end{equation}
If in addition, for some $0 < b < a$ and all $k \geq 1$, $\step_k / \step_{k+1} \leq 1+b \step_{k+1}$, then
\begin{equation}
\label{eq:sum-for-squares}
\sum_{j=1}^{k+1} \step_j^2 \prod_{l=j+1}^{k+1}\left(1-\step_l a\right) \leq \step_{k+1}/ (a-b) \eqsp.
\end{equation}
\end{lemma}
\begin{proof}
Consider first \eqref{eq:sum}.
Let us denote $u_{j: k+1} \eqdef \prod_{l=j}^{k+1}\left(1-\step_l a\right)$. Then, for $j \in\{1, \ldots, k+1\}, u_{j+1: k+1}-u_{j:k+1}=$ $a \step_j u_{j+1:k+1}$. Hence,
\begin{align}
\sum_{j=1}^{k+1} \step_j \prod_{l=j+1}^{k+1}\left(1-\step_l a\right)
&= \step_{k+1}+ \sum_{j=1}^{k} \step_j u_{j+1:k+1} \notag \\
&=\step_{k+1}+ \frac{1}{a} \sum_{j=1}^{k}\left(u_{j+1: k+1}-u_{j: k+1}\right) \notag \\
&=a^{-1}\left(1-u_{1: k+1}\right). \label{eq:auxlemma_statement1}
\end{align}
Consider now \eqref{eq:sum-for-squares}.
\begin{align*}
\sum_{j=1}^{k+1} \step_j^2 \prod_{l=j+1}^{k+1}\left(1-\step_l a\right)&=\step_{k+1} \sum_{j=1}^{k+1}  \frac{\step_j}{\step_{k+1}} \step_j \prod_{l=j+1}^{k+1}\left(1-\step_l a\right) \\
&= \step_{k+1} \sum_{j=1}^{k+1} \step_j \prod_{l=j+1}^{k+1} \frac{\step_{l-1}}{\step_l}\left(1-\step_l a\right)
\end{align*}
where in the last equality, we used
\[
\frac{\step_j}{\step_{k+1}} = \frac{\step_k}{\step_{k+1}} \frac{\step_k-1}{\step_{k}} \cdots \frac{\step_j}{\step_{j+1}}.
\]
Note that, since $\step_{l-1} / \step_{l} \leq 1 + b \step_{l}$ for $l \geq 2$, we have
$$
\frac{\step_{l-1}}{\step_l}\left(1-\step_l a\right) \leq\left(1+b \step_l\right)\left(1-\step_l a\right) \leq 1-(a-b) \step_l.
$$
Substituting into the above inequality yields
$$
\sum_{j=1}^{k+1} \step_j^2 \prod_{l=j+1}^{k+1}\left(1-\step_l a\right) \leq \step_{k+1} \sum_{j=1}^{k+1} \step_j \prod_{l=j+1}^{k+1}\left(1-\step_l (a-b)\right).
$$
Applying \eqref{eq:auxlemma_statement1} implies 
$\sum_{j=1}^{k+1} \step_j^2 \prod_{l=j+1}^{k+1}\left(1-\step_l a\right) \leq \step_{k+1}  \frac{1}{a-b}$ and thus the lemma.
\end{proof}

\subsection{Proofs of \Cref{subsec:csgd-checking}} 
\label{subsec:proof:csgd-checking}
Throughout the proof, we will use the shorthand notations 
\begin{align*}
\bY_{k+1} &\eqdef \Hg(\prm_k, \State_{k+1}), \\
\bZ_{k+1} & \eqdef \compressor(\Hg(\prm_k, \State_{k+1}), \bU_{k+1}).
 \end{align*}
 Observe that by \Cref{assum:field} for the oracle $\bY_{k+1}$,
 \begin{equation} \label{eq:bias:mcg}
 \| \CPE{\bY_{k+1}}{\mcg_k}  - \hg(\prm_k)\|^2 \leq \boundbias_{0,k} + \boundbias_{1,k} \superlyap(\prm_k).
 \end{equation}
 For the proofs of \Cref{lem:csgd-var1} and \Cref{lem:csgd-var2}, define the filtrations $\{\mcg_{k+\nicefrac{1}{2}}, k \in \nset\}$ by 
\begin{align}\label{eq:filtration_X_then_U}
\mcg_{k+\nicefrac{1}{2}} & \eqdef \sigma\left( \prm_0, \State_1, \bU_1, \cdots, \State_k, \bU_k, \State_{k+1} \right).
\end{align}
 \Cref{assum:SGD:contractive} and \Cref{assum:SGD:URBV} claim
 \begin{equation} \label{eq:proof:csgd:CSG1}
\CPE{\| \bZ_{k+1} - \bY_{k+1}\|^2}{\mcg_{k+\nicefrac{1}{2}}} \leq \omgbiased \| \bY_{k+1}\|^2,
 \end{equation}
 and \Cref{assum:SGD:URBV}
 \begin{equation}\label{eq:proof:csgd:CSG2}
 \CPE{\bZ_{k+1}}{\mcg_{k+ \nicefrac{1}{2}}} =\bY_{k+1}. 
 \end{equation}

For the proof of \Cref{lem:csgd-var3}, we define the filtrations $\{\mcg_{k+\nicefrac{1}{2}}, k \in \nset\}$ by 
\begin{align}\label{eq:filtration_U_then_X}
\mcg_{k+\nicefrac{1}{2}} & \eqdef \sigma\left( \prm_0, \bU_1,  \State_1, \cdots, \bU_k, \State_k,  \bU_{k+1} \right).
\end{align}
 
\begin{proof}[Proof of \Cref{lem:csgd-var1}]
Let $k \geq 0$.  
 We  first prove \Cref{assum:field}-\ref{item:field:boundbias} for the compressed oracle $\bZ_{k+1}$. We write, {for any $\zeta_1 > 0$} 
\begin{multline*}
 \| \CPE{\bZ_{k+1}}{\mcg_{k}} - \hg(\prm_{k} ) \|^2 
 \leq {(1+\zeta_1)} \| \CPE{ \bZ_{k+1} -\bY_{k+1} }{ \mcg_k }\|^2 \\
 +  {(1+\zeta_1^{-1})} \|  \CPE{ \bY_{k+1} }{ \mcg_k }- \hg(\prm_k) \|^2.
\end{multline*}
The second term is upper bounded by \eqref{eq:bias:mcg}. For the first one, by convexity of $\| \cdot \|^2$ and \eqref{eq:proof:csgd:CSG1}, it holds that
\begin{align*}
& \| \CPE{ \bZ_{k+1} - \bY_{k+1}}{ \mcg_k } \|^2  \leq    \CPE{ \| \bZ_{k+1} - \bY_{k+1} \|^2 }{ \mcg_k }  \\
& \leq \CPE{ \CPE{ \| \bZ_{k+1}-  \bY_{k+1} \|^2 }{ \mcg_{k+\nicefrac{1}{2}} } }{ \mcg_k } \leq \omgbiased \ \CPE{ \| \bY_{k+1} \|^2 }{ \mcg_k }.
\end{align*}
We further observe that, by definition of the conditional expectation
\begin{align*}
 \CPE{ \| \bY_{k+1} \|^2 }{ \mcg_k }  & =  \| \CPE{\bY_{k+1}}{\mcg_k}\|^2  \\
 \qquad & + \CPE{ \| \bY_{k+1} -  \CPE{\bY_{k+1}}{\mcg_k}\|^2 }{ \mcg_k } \\
 & \leq  (1+\zeta_2)  \| \CPE{\bY_{k+1}}{\mcg_k} - \hg(\prm_k)\|^2  \\
 & + (1+\zeta_2^{-1}) \| \hg(\prm_k)\|^2 \\
 \qquad & + \CPE{ \| \bY_{k+1} -  \CPE{\bY_{k+1}}{\mcg_k}\|^2 }{ \mcg_k }, 
\end{align*}
for any $\zeta_2>0$.
This yields, by using \Cref{assum:field},
\begin{multline} \label{eq:boundY}
 \CPE{ \| \bY_{k+1} \|^2 }{ \mcg_k }  \leq  (1+\zeta_2) \left( \boundbias_{0,k} + \boundbias_{1,k}\superlyap(\prm_k) \right) \\
 + (1+\zeta_2^{-1}) \left( \clyap_{h,0}+ \clyap_{h,1} \superlyap(\prm_k) \right)  + \boundvar_0 +\boundvar_1 \superlyap(\prm_k).
\end{multline}
Hence, \Cref{assum:field}-\ref{item:field:boundbias} holds for $\bZ_{k+1}$ with the constants
\begin{multline}
\boundbias_{\ell,k; \compressor} \eqdef  ( (1+\zeta_1) +(1+\zeta_2)  (1+\zeta_1^{-1})\omgbiased) \boundbias_{\ell,k} + \\ (1+\zeta_2^{-1})(1+\zeta_1^{-1}) \omgbiased \clyap_{h,\ell} + (1+\zeta_1^{-1}) \omgbiased \boundvar_\ell. \nonumber
\end{multline}
We now prove \Cref{assum:field}-\ref{item:field:conditional-variance} for the compressed oracle $\bZ_{k+1}$. Using again the convexity of $\| \cdot \|^2$, it holds
\begin{align*}
& \CPE{\| \bZ_{k+1} - \CPE{\bZ_{k+1}}{\mcg_{k}}\|^2}{\mcg_{k}} \\
 &= \CPE{\CPE{\| \bZ_{k+1} - \CPE{\bZ_{k+1}}{\mcg_{k}}\|^2}{\mcg_{k+\nicefrac{1}{2}}}}{\mcg_{k}} \\
 & \leq  \CPE{ \CPE{\| \bZ_{k+1} - \bY_{k+1} \|^2}{\mcg_{k+\nicefrac{1}{2}}}}{\mcg_{k}} 
 \end{align*}
using the fact that $\bY_{k+1}$ is $\mcg_{k+\nicefrac{1}{2}}$ measurable, and  the fact that for any real random variable $U$, we have $\PE[U] = \argmin_{c\in \rset} \PE[\|U-c\|^2]$.

By \eqref{eq:proof:csgd:CSG1}, we write 
\begin{multline*}
\CPE{\CPE{\| \bZ_{k+1} - \bY_{k+1} \|^2}{\mcg_{k+\nicefrac{1}{2}}}}{\mcg_k}\leq \omgbiased \CPE{\|  \bY_{k+1}\|^2}{\mcg_k};
\end{multline*}
the expectation in the \rhs\ is upper bounded by \eqref{eq:boundY}.
which yields the bound \eqref{eq:increase_var_contractive}.
\end{proof}

\begin{proof}[Proof of \Cref{lem:csgd-var2}]
We follow the same lines as in the proof of \Cref{lem:csgd-var1} except that we use \eqref{eq:proof:csgd:CSG2}.
Let $k \geq 0$.  
 We  first prove \Cref{assum:field}-\ref{item:field:boundbias} for the compressed oracle $\bZ_{k+1}$. We write 
\[
 \| \CPE{\bZ_{k+1}}{\mcg_{k}} - \hg(\prm_{k} ) \|^2  = \| \CPE{ \CPE{\bZ_{k+1}}{\mcg_{k+\nicefrac{1}{2}}}}{\mcg_{k}} - \hg(\prm_{k} ) \|^2;
 \]
 with   \eqref{eq:proof:csgd:CSG2} and \Cref{assum:field} for $\bY_{k+1}$, this yields
 \[
  \| \CPE{\bZ_{k+1}}{\mcg_{k}} - \hg(\prm_{k} ) \|^2 \leq \boundbias_{0,k} + \boundbias_{1,k} \superlyap(\prm_k).
 \]
 Hence  we have $\boundbias_{\ell,k; \compressor} \eqdef \boundbias_{\ell,k}$ for $\ell \in \{0,1\}$.

We now verify \Cref{assum:field}-\ref{item:field:conditional-variance} for $\bZ_{k+1}$. Using \eqref{eq:proof:csgd:CSG2}, we write 
\begin{align*}
& \CPE{ \| \bZ_{k+1} - \CPE{\bZ_{k+1}}{\mcg_{k}} \|^2}{\mcg_{k}}  \\
& =  \CPE{ \| \bZ_{k+1} - \CPE{\bY_{k+1}}{\mcg_{k}} \|^2}{\mcg_{k}} \\
& =  \CPE{ \| \bZ_{k+1} - \bY_{k+1} \|^2}{\mcg_{k}}   +  \CPE{ \| \bY_{k+1} - \CPE{\bY_{k+1}}{\mcg_{k}} \|^2}{\mcg_{k}} \\
& + 2 \CPE{\pscal{\bZ_{k+1} - \bY_{k+1}}{\bY_{k+1} - \CPE{\bY_{k+1}}{\mcg_{k}} }}{\mcg_k}.
\end{align*}
For the scalar product, we have
\begin{multline*}
\CPE{\pscal{\bZ_{k+1} - \bY_{k+1}}{\bY_{k+1} - \CPE{\bY_{k+1}}{\mcg_{k}} }}{\mcg_{k+\nicefrac{1}{2}}} \\
= \pscal{\CPE{\bZ_{k+1}}{\mcg_{k+\nicefrac{1}{2}}} - \bY_{k+1}}{\bY_{k+1} - \CPE{\bY_{k+1}}{\mcg_{k}}},
\end{multline*}
since $\bY_{k+1} \in \mcg_{k+\nicefrac{1}{2}}$ and $\mcg_k \subset \mcg_{k+\nicefrac{1}{2}}$. By \eqref{eq:proof:csgd:CSG2}, the scalar product is zero. Therefore, by \Cref{assum:field}-\ref{item:field:conditional-variance} for $\bY_{k+1}$
\begin{align*}
& \CPE{ \| \bZ_{k+1} - \CPE{\bZ_{k+1}}{\mcg_{k}} \|^2}{\mcg_{k}}  \\
& \leq  \CPE{ \| \bZ_{k+1} - \bY_{k+1} \|^2}{\mcg_{k}}   +  \CPE{ \| \bY_{k+1} - \CPE{\bY_{k+1}}{\mcg_{k}} \|^2}{\mcg_{k}} \\
& \leq \CPE{ \| \bZ_{k+1} - \bY_{k+1} \|^2}{\mcg_{k}}   +  \boundvar_0 + \boundvar_1 \superlyap(\prm_k).
\end{align*}
By  \eqref{eq:proof:csgd:CSG1}, we have
\[
\CPE{ \| \bZ_{k+1} - \bY_{k+1} \|^2}{\mcg_{k}}  \leq \omgunbiased \, \CPE{\| \bY_{k+1} \|^2}{\mcg_k}.
\]
Using \eqref{eq:boundY} with $\zeta_1 = \zeta_2 = 1$
yields our conclusion.
\end{proof}

\begin{proof}[Proof of \Cref{lem:csgd-var3}]
Set  $\bW_{k+1} \eqdef  \Hg( \compressor(\prm_k,\bU_{k+1}) , \State_{k+1} )$ and $\tilde \prm_{k+1} \eqdef \compressor(\prm_k,\bU_{k+1})$.
By  \Cref{assum:SGD:UnifBoundedVar} we have that $\prm_k $ is $ \mcg_{k}$-measurable,  that $\tilde \prm_{k+1} $ is $ \mcg_{k+\nicefrac{1}{2}}$-measurable (as defined in \eqref{eq:filtration_U_then_X}), and that 
\begin{align}
   \CPE{ \|\tilde \prm_{k+1} -\prm _k\|^2 }{\mcg_{k}} &=
   \CPE{ \| \compressor(\prm_k,\bU_{k+1}) -\prm _k\|^2 }{\mcg_{k}} \nonumber\\ 
   &\le \omguniform. \label{eq:squared_moment_unif_bound}
\end{align} 
We first prove  \Cref{assum:field}-\ref{item:field:boundbias} for $\bW_{k+1}$.
\begin{align*} 
& \| \CPE{ \bW_{k+1} }{\mcg_k} - \hg(\prm_k) \|^2    \\
& \qquad =  \| \CPE{  \bW_{k+1}  - \hg(\tilde \prm_k) }{\mcg_k}  + \CPE{   \hg(\tilde \prm_k)- \hg(\prm_k) }{\mcg_k} \|^2.
\end{align*}
Thus for any $\zeta \in \bar{\rset}_{+}$, we get:
\begin{align*} 
& \| \CPE{ \bW_{k+1} }{\mcg_k} - \hg(\prm_k) \|^2    \le  (1+\zeta) \| \CPE{  \bW_{k+1}  - \hg(\tilde \prm_k) }{\mcg_k} \|^2 \\ 
& \qquad \qquad +  (1+\zeta^{-1}) \| \CPE{   \hg(\tilde \prm_k)- \hg(\prm_k) }{\mcg_k} \|^2.
\end{align*}
As $\hg$ is $\Lip{\hg}$ Lipschitz, we have an upper bound for the second term: $  \hg(\tilde \prm_k)- \hg(\prm_k) \le \Lip{\hg} \|\tilde \prm_k -\prm_k\|$ almost surely, thus: 
\begin{align}
    \| \CPE{  \bW_{k+1}  - \hg(\tilde \prm_k) }{\mcg_k} \|^2 \le \Lip{\hg}^2 \omguniform .
\end{align}
Moreover, by Jensen inequality,
\begin{align*}
& \|\CPE{\CPE{\bW_{k+1}  - \hg(\tilde \prm_k)}{\mcg_{k+\nicefrac{1}{2}}}}{\mcg_k}\|^2 \\ 
& \qquad = \|\CPE{\CPE{ \Hg( \tilde \prm_k , \State_{k+1} )
  - \hg(\tilde \prm_k)}{\mcg_{k+\nicefrac{1}{2}}}}{\mcg_k} \|^2 \\
&   \qquad \le  \CPE{\|\CPE{ \Hg( \tilde \prm_k , \State_{k+1} )
  - \hg(\tilde \prm_k)}{\mcg_{k+\nicefrac{1}{2}}} \|^2}{\mcg_k}
\end{align*}
and by \Cref{assum:field}-\ref{item:field:boundbias} with $\boundbias_1=0$, $\|\CPE{ \Hg( \tilde \prm_k , \State_{k+1} ) - \hg(\tilde \prm_k)}{\mcg_{k+\nicefrac{1}{2}}} \|^2 \le \boundbias_{0,k} + \boundbias_{1,k} \superlyap(\tilde \prm_k) \le \boundbias_0$.
Overall, \Cref{assum:field}-\ref{item:field:boundbias} is satisfied with $\boundbias_{0,\compressor} = (1+\zeta) \boundbias_0  + (1+\zeta^{-1})\Lip{\hg}^2 \omguniform$ and $\boundbias_{1,\compressor} = 0$.

We now prove \Cref{assum:field}-\ref{item:field:conditional-variance} for $\bW_{k+1}$:
\begin{align*}
& \CPE{ \| \bW_{k+1} - \CPE{ \bW_{k+1} }{\mcg_k} \|^2 }{\mcg_k}  \\
& \qquad = \CPE{ \CPE{ \| \bW_{k+1} - \CPE{ \bW_{k+1} }{\mcg_{k+\nicefrac{1}{2}}} \|^2}{\mcg_{k+\nicefrac{1}{2}}} }{\mcg_k}  \\
& \qquad \quad + \CPE{ \|  \CPE{ \bW_{k+1} }{\mcg_{k+\nicefrac{1}{2}}} - \CPE{ \bW_{k+1} }{\mcg_k} \|^2 }{\mcg_k}  .
\end{align*}
By \Cref{assum:field}-\ref{item:field:conditional-variance}, with $\boundvar_1 =0 $
\begin{align*}
   &  \CPE{ \| \bW_{k+1} - \CPE{ \bW_{k+1} }{\mcg_{k+\nicefrac{1}{2}}} \|^2}{\mcg_{k+\nicefrac{1}{2}}}  \\
    & =  \CPE{ \| \Hg(\tilde \prm_k ,\State_{k+1}) - \CPE{ \Hg(\tilde \prm_k , \State_{k+1}) }{\mcg_{k+\nicefrac{1}{2}}} \|^2}{\mcg_{k+\nicefrac{1}{2}}} \\
    & \le \boundvar_0 + \boundvar_1 \superlyap(\tilde \prm_k) \le \boundvar_0. 
\end{align*}
Moreover, 
\begin{align*}
  &   \CPE{ \|  \CPE{ \bW_{k+1} }{\mcg_{k+\nicefrac{1}{2}}} - \CPE{ \bW_{k+1} }{\mcg_k} \|^2 }{\mcg_k}\\
    \le &  \CPE{ \|  \CPE{\Hg(\tilde \prm_k ,\State_{k+1})}{\mcg_{k+\nicefrac{1}{2}}} - \CPE{ \Hg(\prm_k, \State_{k+1}) }{\mcg_k} \|^2 }{\mcg_k} \\
    = &  \CPE{ \|  \CPE{\Hg(\tilde \prm_k ,\State_{k+1})}{\mcg_{k+\frac{1}{2}}} - \CPE{ \Hg(\prm_k, \State_{k+1}) }{\mcg_{k+\frac{1}{2}}} \|^2 }{\mcg_k} .
\end{align*}
By assumption,  $\CPE{ \Hg( \cdot , \State_{k+1}) }{\mcg_{k+\nicefrac{1}{2}}} $ is $\Lip{\PE{\Hg}}$-Lipschitz. We conclude using \eqref{eq:squared_moment_unif_bound}:
\begin{align*}
  &   \CPE{ \|  \CPE{ \bW_{k+1} }{\mcg_{k+\nicefrac{1}{2}}} - \CPE{ \bW_{k+1} }{\mcg_k} \|^2 }{\mcg_k} \le \Lip{\PE{\Hg}}^2 \omguniform. 
\end{align*}
We get the result by combining the two bounds above.
\end{proof}
\begin{proof}[Proof of \Cref{lem:csgd-var4}]
Under the constant stepsize assumption, the random field can be written as
\begin{equation} \notag
\begin{aligned}
& \widetilde{\Hg}( \prm_k, \bU_{k+1}, \State_{k+1} ) = \frac{1}{\bar{\step}} \left( \compressor( \prm_k + \bar{\step} \, \bY_{k+1} , \bU_{k+1} ) - \prm_k \right) .
\end{aligned}
\end{equation}
Notice that \Cref{assum:field}-\ref{item:field:boundbias} can be easily verified  since 
\begin{align*}
\CPE{\widetilde{\Hg}( \prm_k, \bU_{k+1}, \State_{k+1} ) }{\mcf_k} = \CPE{ \bY_{k+1} }{\mcf_k}.
\end{align*}

To verify \Cref{assum:field}-\ref{item:field:conditional-variance}, we proceed by
\begin{align*}
& \bar{\step}^2 \CPE{ \| \widetilde{\Hg}( \prm_k, \bU_{k+1}, \State_{k+1} ) -  \CPE{ \bY_{k+1} }{\mcf_k} \|^2 }{\mcf_k} \\
& = \CPE{ \| \compressor( \prm_k + \bar{\step} \bY_{k+1} , \bU_{k+1} ) - ( \prm_k + \bar{\step} \CPE{ \bY_{k+1} }{\mcf_k} )  \|^2 }{ \mcf_k }.
\end{align*}
As $\prm_k \in \linearset^d$, applying \Cref{assum:SGD:UniformQuant} gives
\begin{align*}
& \CPE{ \| \widetilde{\Hg}( \prm_k, \bU_{k+1}, \State_{k+1} ) -  \CPE{ \bY_{k+1} }{\mcf_k} \|_2^2 }{\mcf_k} \\
& \leq \CPE{ \| \bY_{k+1} - \CPE{ \bY_{k+1} }{\mcf_k} \|_2^2 }{\mcf_k} + \frac{ \linearunbiased }{ \bar{\step} } \CPE{ \| \bY_{k+1} \|_1 }{\mcf_k } \\
& \leq \boundvar_0 + \boundvar_1 \superlyap( \prm_k ) + \frac{ \linearunbiased }{ \bar{\step} } \CPE{ \| \bY_{k+1} \|_1 }{\mcf_k }.
\end{align*}
Lastly, we obtain the following chain
\begin{align*}
& \CPE{ \| \bY_{k+1} \|_1 }{\mcf_k } \leq \CPE{ \| \CPE{ \bY_{k+1} }{ \mcf_k } - \bY_{k+1} \|_1 }{\mcf_k} \\
& \quad \qquad \qquad \| \hg( \prm_k ) \|_1 + \| \CPE{ \bY_{k+1} }{ \mcf_k } - \hg( \prm_k ) \|_1 \\
& \leq \sqrt{d} \left( \sqrt{ \boundbias_{0,k} + \boundbias_{1,k}\superlyap(\prm_k) } + \sqrt{\clyap_{h,0}+ \clyap_{h,1} \superlyap(\prm_k)} \right) \\
& \quad + \sqrt{d} \sqrt{ \boundvar_0 + \boundvar_1 \superlyap( \prm_k ) } \\
& \leq \sqrt{d} \left( \frac{3 + \boundbias_{0,k} + \clyap_{h,0} + \boundvar_0 }{2}  + \frac{\boundbias_{1,k} + \clyap_{h,1} + \boundvar_1 }{2} \superlyap( \prm_k ) \right) ,
\end{align*}
where we have used $\sqrt{x} \leq \frac{1+x}{2}$, $x \geq 0$ in the last inequality.
Collecting terms from the above bounds lead to the desired terms in~\eqref{eq:increase_var_linearQuant}.
\end{proof}

\subsection{Proofs of \Cref{subsec:TD-learning-checking}}
\label{subsec:proof:TD-learning-checking}

\begin{proof}[Proof of \Cref{lem:contraction-projected-bellman}]
Let $\prm , \prm' \in \rset^d$.
Since a projection is a contraction.
\[
\left\|\Proj_{\statdistMRP} \bellman \valuefunc[\prm]-\Proj_{\statdistMRP} \bellman \valuefunc[\prm'] \right\|_{\Dtd_\statdistMRP} \leq \left\|\bellman \valuefunc[\prm]- \bellman \valuefunc[\prm'] \right\|_{\Dtd_\statdistMRP}.
\]
From the expression of $\bellman$ (see \eqref{eq:definition-bellman}), we write
\[
\bellman \valuefunc[\prm](s)- \bellman \valuefunc[\prm'](s) = \lambda  \sum_{s'} \kerMRP(s,s') \left( \valuefunc[\prm](s') - \valuefunc[\prm'](s')\right).
\]
The squared $\Dtd_\statdistMRP$-norm of the \rhs\ is equal to
\[
\lambda^2 \sum_s \statdistMRP(s) \left(\sum_{s'} \kerMRP(s,s') \left( \valuefunc[\prm](s') - \valuefunc[\prm'](s')\right) \right)^2.
\]
Since $\kerMRP$ is a transition kernel, we have the inequality
\begin{multline*}
\left(\sum_{s'} \kerMRP(s,s') \left( \valuefunc[\prm](s') - \valuefunc[\prm'](s')\right) \right)^2  \\
\leq \sum_{s'} \kerMRP(s,s')  \left(\left( \valuefunc[\prm](s') - \valuefunc[\prm'](s')\right) \right)^2.
\end{multline*}
Finally, since  $\statdistMRP \kerMRP = \statdistMRP$ by \Cref{assum:TD:stationary-policy}, we have
\begin{align*}
\sum_s \statdistMRP(s) & \sum_{s'} \kerMRP(s,s')  \left(\left( \valuefunc[\prm](s') - \valuefunc[\prm'](s')\right) \right)^2  \\
& = \sum_{s'}  \statdistMRP(s')  \left(\left( \valuefunc[\prm](s') - \valuefunc[\prm'](s')\right) \right)^2 \\
& = \|  \valuefunc[\prm] - \valuefunc[\prm']\|^2_{\Dtd_\statdistMRP}.
\end{align*}
This concludes the proof.
 \end{proof}
Set
\begin{equation}
\label{eq:definition-G}
\boldsymbol{\mathcal{G}}(s,s') \eqdef  \feature(s) \{\lambda \feature(s') - \feature(s)\}^{\top} \, ;
\end{equation}
we may rewrite $\Hg(\prm,(s,s'))$ in \eqref{eq:iteration-TD0} as
\begin{equation}
\label{eq:TD-field-equivalent-expression}
\Hg(\prm,(s,s'))= \feature(s)\rewardMRP(s,s') + \boldsymbol{\mathcal{G}}(s,s') \prm \eqsp.
\end{equation}
We denote
\[
X_0 \eqdef (S_0,S'_0),  \qquad \bar{\boldsymbol{\boldsymbol{\mathcal{G}}}}_\statdistMRP(s,s') \eqdef  \boldsymbol{\mathcal{G}}(s,s') - \PE_{\statdistMRP}[\boldsymbol{\mathcal{G}}(X_0)].
\]
Let $\prm_\star$  be such that $\valuefunc[\star] = \Feature \prm_\star$.  Since $\hg(\prm_\star)=0$ and  $\hg(\prm) = \PE_{\statdistMRP}\left[ \Hg(\prm, \State_0) \right]$, we get
 \begin{equation} \label{eq:TD:diff:hg}
\hg(\prm) = \hg(\prm) - \hg(\prm_\star) = \PE_{\statdistMRP}[\boldsymbol{\mathcal{G}}(X_0)] (\prm - \prm_\star).
 \end{equation}
 We also have
 \[
 \Hg(\prm,(s,s')) - \Hg(\prm_\star,(s,s'))  = \boldsymbol{\mathcal{G}}(s,s') \, \left(\prm - \prm_\star\right).
 \]
 This yields
\begin{multline}
\label{eq:decomposition-field-TD}
\Hg(\prm,(s,s')) - \hg(\prm)= \Hg(\prm_\star,(s,s')) \\ +  \bar{\boldsymbol{\mathcal{G}}}_\statdistMRP(s,s') (\prm - \prm_\star) \eqsp.
\end{multline}
Set
\begin{equation}\label{eq:TD:def:xi}
\xi_{\prm}(s) \eqdef (\prm - \prm_\star)^\top \feature(s) = \valuefunc[\prm](s)-\valuefunc[\prm_\star](s).
\end{equation}
From \eqref{eq:definition-G} and \eqref{eq:TD:diff:hg}, we have
\begin{align}
\hg(\prm) & = \PE_{\statdistMRP}[\boldsymbol{\mathcal{G}}(S_0,S'_0)]  \, ( \prm - \prm_\star)  \nonumber \\
& = \PE_{\statdistMRP}[ \feature(S_0) \{ \lambda \xi_{\prm}(S'_0) - \xi_{\prm}(S_0)\} ]. \label{eq:TD:hg:diff:xi}
\end{align}
\begin{proof}[Proof of \Cref{lem:TD:bound-mean-field}]
Using \eqref{eq:TD:hg:diff:xi}, we get that
\begin{align*}
\| \hg(\prm) \|^2
&= \| \PE_{\statdistMRP}[ \feature(S_0) \{\lambda \xi_{\prm}(S'_0) - \xi_{\prm}(S_0) \}]  \|^2  \\
&\leq \PE_{\statdistMRP}[ \| \feature(S_0) \|^2]   \ \PE_{\statdistMRP}[ \{\lambda \xi_{\prm}(S'_0) - \xi_{\prm}(S_0) \}^2].
\end{align*}
For the second term, we use
\[
\PE[(\lambda U+V)^2] \leq \lambda^2 \PE[U^2] + \PE[V^2] + 2 \lambda \sqrt{\PE[U^2] \, \PE[V^2]}
 \]
with $U \eqdef \xi_{\prm}(S'_0)$ and $V \eqdef \xi_{\prm}(S_0)$.
By  \Cref{assum:TD:stationary-policy},
\[
 \PE_{\statdistMRP}[  \xi_{\prm}^2(S_0')]  =  \PE_{\statdistMRP}[  \xi_{\prm}^2(S_0)]  = \| \valuefunc[\prm]-\valuefunc[\star]\| _{\Dtd_{\statdistMRP}}^2
\]
which implies that
\begin{align*}
& \PE_{\statdistMRP}[ \{\lambda \xi_{\prm}(S'_0) - \xi_{\prm}(S_0) \}^2]  \leq  (1+\lambda)^2   \| \valuefunc[\prm]-\valuefunc[\star]\| _{\Dtd_{\statdistMRP}}^2.
\end{align*}
This concludes the proof of the first inequality.
From \eqref{eq:decomposition-field-TD}  and \Cref{assum:TD:sampling}, we  obtain
\begin{align*}
\nonumber
&\CPE{\| \Hg(\prm_{k},\State_{k+1}) - \hg(\prm_{k})\|^2}{\mcf_{k}} \\
&\leq 2 \, \PE_{\statdistMRP}[ \| \Hg(\prm_\star,X_0) \|^2] \\
\nonumber
&+ 2 \, (\prm_{k} - \prm_\star)^{\top} \PE_{\statdistMRP}[ \bar{\boldsymbol{\mathcal{G}}}_\statdistMRP(X_0)  \bar{\boldsymbol{\mathcal{G}}}^\top_\statdistMRP(X_0) ] (\prm_{k} - \prm_\star) \eqsp.
\end{align*}
We first upper bound $\PE_{\statdistMRP}[ \| \Hg(\prm_\star,X_0) \|^2]$. Since $\lambda \leq 1$ and $|\rewardMRP(s,s')| \leq 1$ (see \Cref{assum:TD:normed:calR}), we get that, for all $\prm \in \rset^d$, and $(s,s')$, it holds that
\begin{align*}
\|\Hg(\prm,(s,s'))\|& \leq 1 +
\| \feature(s) \{\lambda \valuefunc[\star](s') + \valuefunc[\star](s) \}\|
\\
 &  \leq  1 +     \{ \lambda |\valuefunc[\star](s')| + |\valuefunc[\star](s)| \} \eqsp.
\end{align*}
We obtain
\[
\PE_{\statdistMRP}[ \| \Hg(\prm_\star,X_0) \|^2] \leq 3  \left( 1 +     \{ \lambda^2 +1 \} \| \valuefunc[\star]\|^2_{\Dtd_{\statdistMRP}} \right).
\]
Let us upper bound the second term. For any ${\bf u} \in \rset^d$, it holds that
\[
{\bf u}^{\top} \PE_{\statdistMRP}[ \bar{\boldsymbol{\mathcal{G}}}_\statdistMRP(X_0)  \bar{\boldsymbol{\mathcal{G}}}^\top_\statdistMRP(X_0) ] {\bf u} \leq {\bf u}^{\top} \PE_{\statdistMRP}[\boldsymbol{\mathcal{G}}(X_0) \boldsymbol{\mathcal{G}}^{\top}(X_0)] {\bf u} \eqsp,
\]
and
\[
{\bf u}^T \boldsymbol{\mathcal{G}}(s,s') \boldsymbol{\mathcal{G}}^{\top}(s,s') {\bf u} \leq \{ {\bf u}^{\top} \feature(s) \}^2  (1 + \lambda)^2 \eqsp.
\]
By combining these two inequalities, we finally obtain
\begin{multline*}
(\prm_{k} - \prm_\star)^{\top} \PE_{\statdistMRP}[ \bar{\boldsymbol{\mathcal{G}}}_\statdistMRP(X_0)  \bar{\boldsymbol{\mathcal{G}}}^\top_\statdistMRP(X_0) ] (\prm_{k} - \prm_\star) \\
\leq (1+\lambda)^2  \| \valuefunc[\prm_k]-\valuefunc[\star]\|^2_{\Dtd_{\statdistMRP}}
\end{multline*}
This concludes the proof.
\end{proof}
\begin{proof}[Proof of \Cref{lem:TD:lyap}]
It follows from \eqref{eq:TD:hg:diff:xi} that
\begin{align*}
& \pscal{\prm-\prm_*}{\hg(\prm)-\hg(\prm_\star)} \\
&\quad = \PE_{\statdistMRP}[\xi_{\prm}(S_0)\{\lambda \xi_{\prm}(S_0') - \xi_{\prm}(S_0)\}] \\
& \quad = \lambda  \PE_{\statdistMRP}[\xi_{\prm}(S_0) \xi_{\prm}(S_0') ] -  \PE_{\statdistMRP}[(\xi_{\prm}(S_0))^2].
\end{align*}
By using the Cauchy-Schwarz inequality, we have
\[
\PE_{\statdistMRP}[\xi_{\prm}(S_0) \xi_{\prm}(S'_0)] \leq \{\PE_{\statdistMRP}[\xi^2_{\prm}(S_0)] \}^{1/2} \{\PE_{\statdistMRP}[\xi^2_{\prm}(S'_0)] \}^{1/2}.
\]
Under \Cref{assum:TD:stationary-policy}, $\statdistMRP \kerMRP= \statdistMRP$,  which implies that for any function $g$,
\[
\PE_{\statdistMRP}[g(S'_0)] = \PE_{\statdistMRP}[g(S_0)].
\]
This yields
\[
\pscal{\prm-\prm_*}{\hg(\prm)-\hg(\prm_\star)} \leq - (1-\lambda) \PE_{\statdistMRP}[\xi^2_{\prm}(S_0)].
\]
The proof is concluded by using
\eqref{eq:TD:def:xi} and
\begin{align*}
\PE_{\statdistMRP}[\xi^2_{\prm}(S_0)]  &  =  \PE_{\statdistMRP}[ \left(\valuefunc[\prm](S_0)-\valuefunc[\star](S_0) \right)^2] \\
& =  \left\|\valuefunc[\prm]-\valuefunc[\star] \right\|_{\Dtd_{\statdistMRP}}^2.
\end{align*}
\end{proof}
\begin{proof}[Proof of \Cref{lem:TD:upper-lower-spectrum-covariance}]
Note first that
\begin{align*}
\prm^{\top} \Feature^\top \Dtd_{\statdistMRP} \Feature \prm
&= \sum_{s \in \stateMRP} \statdistMRP(s) \pscal{\feature(s)}{\prm}^2 \\
&\leq \sum_{s \in \stateMRP} \statdistMRP(s) \| \feature(s) \|^2 \| \prm \|^2 \leq \| \prm \|^2,
\end{align*}
where we have used that  $\sum_{s \in \stateMRP} \statdistMRP(s)=1$. On the other hand, under
\Cref{assum:TD:stationary-policy},  $\Dtd_{\statdistMRP}$ is full rank ($\Dtd_{\statdistMRP}$ is diagonal and all its diagonal entries are positive), which implies that  the minimal eigenvalue of $\Feature^\top \Dtd_{\statdistMRP} \Feature$ is positive. The result follows.
\end{proof}

\subsection{Proofs of \Cref{subsec:finite-time-bounds}} \label{subsec:proof:finite-time-bounds}
We hereafter give a slightly different statement of~\Cref{lem:main-result-quantitative} for the particular case in which $\nabla V= -\hg$. 
\label{app:RS-forGD-general}
\begin{assumption}
\label{assum:lyapunov_partcase}
The field $h$ and the function $\lyap$ in \Cref{assum:lyapunov} satisfy:
$\nabla \lyap =- \hg$. 
\end{assumption}
Then under~\Cref{assum:lyapunov}-\ref{item:rholyap},
the Lyapunov functions $\lyap$, $\superlyap$ satisfy, for any $\prm \in \rset^d$, 
\begin{equation}\label{eq:link-lyap-superlyap-reversed}
   - \| \nabla \lyap( \prm) \|^2 =\pscal{ \nabla \lyap( \prm)}{\hg(\prm)} \leq - \rholyap \superlyap(\prm) \eqsp.
\end{equation}
Define, for any $k\geq 0$:
{\small
\begin{align}
     \omega_{2,k+1} & \eqdef (\rholyap- \boundbias_{1}) - \step_{k+1} \Liplyap \boundvar_1  \label{eq:omega2k}\\
     \step_{2,\max} & \eqdef \min\left(\frac{1}{\Liplyap}  , \frac{(\rholyap- \boundbias_{1})}{ \Liplyap \boundvar_1}\right)  \label{eq:step2k}.
\end{align}
}
\begin{lemma}[Robbins-Siegmund type inequality for $\nabla V = -\hg$]
\label{lem:main-result-quantitative_2}
Assume \Cref{assum:field,assum:lyapunov}, \Cref{assumNA:uniform-bound-bias} and \Cref{assum:lyapunov_partcase}.
Then, for any $k \geq 0$, we have almost-surely
\begin{align}
\nonumber
\CPE{\lyap(\prm_{k+1})}{\mcf_k} &\le \lyap(\prm_k) -\frac{\step_{k+1}}{2} \omega_{2,k+1} \superlyap(\prm_k) \\
&\quad +\frac{\step_{k+1}}{2}  \boundbias_{0}
+ \frac{\step^2_{k+1} \Liplyap }{2} \boundvar_0,
\end{align}
\end{lemma}
\begin{proof}
Let $k \geq 0$. 
By \Cref{assum:lyapunov}-\ref{item:lyap},  we have
\begin{multline*}
\lyap(\prm_{k+1}) \leq \lyap(\prm_k) + \pscal{\nabla \lyap(\prm_k)}{\prm_{k+1}-\prm_k}
\\ + (\Liplyap/2) \|\prm_{k+1} - \prm_k \|^2 \,.
\end{multline*}
Computing the conditional expectation of both sides of this inequality yields
\begin{align*}
&\CPE{\lyap(\prm_{k+1})}{\mcf_k} \leq \lyap(\prm_k) \\
&\quad + \step_{k+1} \pscal{\nabla \lyap(\prm_k)}{\CPE{\Hg(\prm_k, \State_{k+1})}{\mcf_k}} + \\
&\quad + \step^2_{k+1} (\Liplyap /2) \,  \CPEs{\| \Hg(\prm_k,\State_{k+1}) \|^2}{\mcf_k} \,.
\end{align*}
We now use $2\pscal{a}{b} = \|a\|^2+\|b\|^2 -\|a-b\|^2$.
\begin{align*}
 & 2 \pscal{\nabla \lyap(\prm_k)}{\CPE{\Hg(\prm_k, \State_{k+1})}{\mcf_k}} \\ 
&= - 2\pscal{-\nabla \lyap(\prm_k)}{\CPE{\Hg(\prm_k, \State_{k+1})}{\mcf_k}} \\
&= - \|- \nabla \lyap(\prm_k)\|^2 - \|\CPE{\Hg(\prm_k, \State_{k+1})}{\mcf_k}\|^2 \\ 
&+ \|- \nabla \lyap(\prm_k) - \CPE{\Hg(\prm_k, \State_{k+1})}{\mcf_k}\|^2.
\end{align*}
Define $\bias_k \eqdef \CPE{\Hg(\prm_k, \State_{k+1})}{\mcf_k} - \hg(\prm_k)$. We have, using \Cref{eq:link-lyap-superlyap-reversed},
\begin{align*}
&\CPE{\lyap(\prm_{k+1})}{\mcf_k} \leq \lyap(\prm_k) - \frac{\rholyap \step_{k+1}}{2} \superlyap(\prm_k) \\
&\quad  
- \frac{\step_{k+1}}{2}\|\CPE{\Hg(\prm_k, \State_{k+1})}{\mcf_k}\|^2 +\frac{\step_{k+1}}{2} \|\bias_k\|^2 \\
&\quad + \step^2_{k+1} (\Liplyap /2) \,  \CPEs{\| \Hg(\prm_k,\State_{k+1}) \|^2}{\mcf_k} \,.
\end{align*}
Note first that, using \Cref{assum:field}-\ref{item:field:boundbias} we get
\begin{align*}
&\CPE{\lyap(\prm_{k+1})}{\mcf_k} \leq \lyap(\prm_k) - \frac{\rholyap \step_{k+1}}{2} \superlyap(\prm_k) \\
&\quad  
- \frac{\step_{k+1}}{2}\|\CPE{\Hg(\prm_k, \State_{k+1})}{\mcf_k}\|^2 +\frac{\step_{k+1}}{2}  (\boundbias_{0}+ \boundbias_{1} \superlyap(\prm_k)) \\
&\quad + \step^2_{k+1} (\Liplyap /2) \,  \CPEs{\| \Hg(\prm_k,\State_{k+1}) \|^2}{\mcf_k} \,.
\end{align*}
We compute a bias-variance decomposition and use \Cref{assum:field}-\ref{item:field:conditional-variance}:
\begin{align*}
&\CPEs{\| \Hg(\prm_k,\State_{k+1}) \|^2}{\mcf_k}
= \|  \CPEs{\Hg(\prm_{k},\State_{k+1})}{\mcf_{k}}\|^2 
\\
& \quad + \CPE{\| \Hg(\prm_{k},\State_{k+1}) - \CPEs{\Hg(\prm_{k},\State_{k+1})}{\mcf_{k}}\|^2}{\mcf_{k}} \\
&\leq \|  \CPEs{\Hg(\prm_{k},\State_{k+1})}{\mcf_{k}}\|^2   + \boundvar_0 + \boundvar_1 \superlyap(\prm_{k}).
\end{align*}
Overall, we get:
\begin{align*}
&\CPE{\lyap(\prm_{k+1})}{\mcf_k} \leq \lyap(\prm_k) \\& \quad - \left(\frac{\step_{k+1}}{2} (\rholyap- \boundbias_{1}) - \frac{\step^2_{k+1} \Liplyap }{2} \boundvar_1  \right) \superlyap(\prm_k) \\
&\quad  
- \left(\frac{\step_{k+1}}{2} -   \frac{\step^2_{k+1} \Liplyap }{2}  \right ) \|\CPE{\Hg(\prm_k, \State_{k+1})}{\mcf_k}\|^2  \\
&\quad +\frac{\step_{k+1}}{2}  \boundbias_{0}
+ \frac{\step^2_{k+1} \Liplyap }{2} \boundvar_0.
\end{align*}
We thus get:
\begin{align*}
&\CPE{\lyap(\prm_{k+1})}{\mcf_k} \leq \lyap(\prm_k) -\frac{\step_{k+1}}{2} \omega_{2,k+1} \superlyap(\prm_k) \\
&\quad +\frac{\step_{k+1}}{2}  \boundbias_{0}
+ \frac{\step^2_{k+1} \Liplyap }{2} \boundvar_0,
\end{align*}
with $\omega_{2,k}$,  $\step_{2,\max}$ as in \eqref{eq:omega2k}, \eqref{eq:step2k}.
\end{proof}

\subsection{Proofs of \Cref{subsub:exrate:CSGD}}
\label{app:proofGS}
\begin{proof}[Proof of \Cref{prop:rateGS}]
By \Cref{cor:GSascomprGD}, the field ${\Hg(\prm, \sim) }= \rm{Top_1}(\nabla F(\prm))$ satisfies \Cref{assum:field} with $(\clyap_{\hg,0}, \clyap_{\hg,1}) =(0,1)$, $(\boundbias_{0,k}, \boundbias_{1,k})=(0,1-\nicefrac{1}{d})$ and $(\boundvar_0, \boundvar_1)=(0,0)$, for  $\superlyap = \|\nabla F (\cdot)\|^2$). 

We consider $\lyap = F$ and thus have \Cref{assum:lyapunov} with $\clyap_\lyap=1$, $\rho=1$, $\Lip{\lyap} = \Lip{ \nabla F }$. We thus verify 
\Cref{eq:definiton-bsf-0,eq:definition-eta,eq:step-max,eq:definition-u-k} with:
\begin{align*}
\bsf_0 &\eqdef \clyap_{\lyap} \sqrt{\boundbias_0}/2 =  0\\
\bsf_1 &\eqdef \clyap_{\lyap} (\sqrt{\boundbias_0} /2 + \sqrt{\boundbias_1}) = \sqrt{1-\nicefrac{1}{d}} \leq 1-\frac{1}{2d}.\\
\eta_0& \eqdef 0 \\
\eta_1& \eqdef \boundvar_1 + \boundbias_1 +  \clyap_{\hg,1} + \sqrt{\clyap_{\hg,1}} \left( \sqrt{\boundbias_0 }+ \sqrt{\boundbias_1 }\right)  \\
\nonumber
& \qquad + \sqrt{\boundbias_1 } \left( \sqrt{\clyap_{\hg,0}} + \sqrt{\clyap_{\hg,1}}\right) \le 4 , \\
\step_{\max} & \eqdef  2 \{\rholyap- \bsf_1\}/(\Liplyap \eta_1) \geq \frac{1}{d \Liplyap \eta_1}\, ,   \\
\omega_k & \eqdef 2\{ \rholyap - \bsf_1\} - \step_{k} \Liplyap \eta_1 \geq \frac{1}{2d}  \eqsp.
\end{align*}
\Cref{prop:rateGS} is then a direct application of \Cref{theo:main-result-quantitative} with the constants above.
\end{proof}

For the proofs of the next two propositions, we note that the random field $\Hg$ satisfies \Cref{assum:field} with $\lyap = F$, $\superlyap = \|\nabla F(\cdot)\|^2$. 
Moreover, the constants are $(\clyap_{\hg,0}, \clyap_{\hg,1})=(0,1)$, $(\boundbias_{0,k}, \boundbias_{1,k})= (0,0)$, $(\boundvar_0, \boundvar_1) = (M^2/n, 0)$, and $\rholyap = 1$.

\begin{proof}[Proof of \Cref{prop:compSGDunbiased}]
By \Cref{lem:csgd-var2}, the compressed {\SG} in \eqref{eq:cSGD} uses a random field that satisfies \Cref{assum:field} with the same set of constants that inherit from $\Hg$ except for
 $\boundvar_{0; \compressor} \eqdef (1+\omgunbiased) \boundvar_{0} , \qquad  \boundvar_{1; \compressor} \eqdef  2 \omgunbiased .$
Thus we have \eqref{eq:definiton-bsf-0}
with $\bsf_0= \bsf_1 = 0$ and \eqref{eq:definition-eta}, \eqref{eq:step-max}, \eqref{eq:definition-u-k} with 
\begin{align*}
\eta_0& \eqdef (1+\omgunbiased) M^2/n , \qquad \eta_1 \eqdef 2 \omgunbiased  + 1 , \\
\step_{\max} & \eqdef  2 /(\Lip{\nabla F} \eta_1) = 2 /(\Lip{\nabla F} (2\omgunbiased+1)) \, ,   \\
\omega_k & \eqdef 2 - \step_{k} \Lip{\nabla F} \eta_1 \geq  1 \eqsp.
\end{align*}
Where the last equation holds as $\step_k\le \step_{\max}/2$. As before, \Cref{prop:compSGDunbiased} is then a direct application of \Cref{theo:main-result-quantitative} with the constants above.
\end{proof}

\begin{proof}[Proof of \Cref{prop:STE}]
For  $\compressor$ is $Q_d$, \Cref{assum:SGD:UniformQuant} is satisfied with $\omguniform=d\Delta^2$. Using \Cref{cor:STEfield}, we can apply the result of \Cref{theo:main-result-quantitative}, with  
$\boundbias_{0,k,\compressor} \eqdef \Lip{{\nabla F}}^2 d\Delta^2,  \qquad
\boundvar_{0; \compressor} \eqdef   \frac{M^2}{n}+  \, \Lip{\nabla F}^2 d\Delta^2.$
Thus
\begin{align*}
\bsf_0 &=   \Lip{{\nabla F}} \sqrt{d} \Delta/2 \le 1/2 \ , \quad\bsf_1 =  \Lip{{\nabla F}} \sqrt{d} \Delta/2 \le 1/2,\\
\eta_0& = \frac{M^2}{n}+ 2 \, \Lip{\nabla F}^2 d\Delta^2 + \Lip{{\nabla F}} \sqrt{d} \Delta\le \frac{M^2}{n}+ 3  \Lip{{\nabla F}} \sqrt{d} \Delta \\
\eta_1& \eqdef    1+ \left(\Lip{{\nabla F}} \sqrt{d} \Delta\right)\le 2  \\
\step_{\max} & \eqdef  2 \{1- \bsf_1\}/(\Lip{{\nabla F}} \eta_1)   \geq  1 /(2\Lip{{\nabla F}} ) \\
\omega_k & \eqdef 2\{ \rholyap - \bsf_1\} - \step_{k} \Liplyap \eta_1 \geq 1/2 \eqsp.
\end{align*}
\Cref{prop:STE} is then a direct application of \Cref{theo:main-result-quantitative} with the constants above.
\end{proof}

\begin{proof}[Proof of \Cref{prop:compSGDLinear}]
By \Cref{lem:csgd-var4}, the compressed {\SG} in \eqref{eq:de-sa} uses a random field \eqref{eq:de-sa-H} that satisfies \Cref{assum:field} with the same set of constants inherited from $\Hg$ except for
\begin{align*}
& \boundvar_{0; \compressor} \eqdef \frac{M^2}{n} + \frac{ \linearunbiased \sqrt{d} }{2 \bar{\step}} (3 + M^2/n) , \quad  \boundvar_{1; \compressor} \eqdef \frac{ \linearunbiased \sqrt{d} }{2 \bar{\step}}.
\end{align*}
We also have \eqref{eq:definiton-bsf-0}
with $\bsf_0= \bsf_1 = 0$ and \eqref{eq:definition-eta}, \eqref{eq:step-max} with 
\begin{align*}
\eta_0& \eqdef M^2 / n + (M^2/n) \frac{ \linearunbiased \sqrt{d} }{ 2 \bar{\step} } , \quad \eta_1 \eqdef 1 + \frac{ \linearunbiased \sqrt{d} }{2 \bar{\step}} , \\
\step_{\max} & \eqdef { 2 } / ( { \Lip{\nabla F} ( 1 + { \linearunbiased \sqrt{d} } / {(2 \bar{\step})} ) } )  \, .
\end{align*}
To apply \Cref{theo:main-result-quantitative}, we have to satisfy $\bar{\step} < \step_{\max}$, which is then equivalent to 
\begin{align*}
& \bar{\step} < \frac{ 2 }{ \Lip{\nabla F} ( 1 + \frac{ \linearunbiased \sqrt{d} }{2 \bar{\step}} ) } \Longleftrightarrow \bar{\step} + \linearunbiased \sqrt{d} < \frac{ 2 }{ \Lip{\nabla F}  }.
\end{align*}
This yields the stepsize condition in \eqref{eq:main-result-non-asymptotic-LinearQ-step}. 
Our bound in \eqref{eq:main-result-non-asymptotic-LinearQ} is then achieved by plugging the above $\eta_0$ into \eqref{eq:main-result-non-asymptotic-constant}.
\end{proof}

\subsection{Proofs of \Cref{sec:SAEM:complexity}}
\label{sec:appendix:SAEM}
From \Cref{sec:EM:verifH:SAEM}, it holds
\begin{align*}
\eta_0 & =  12 s_\star  \sqrt{\clyap_{\chi,0}} / \nbrMC + 144 s_\star^2 \clyap_{\chi,0} / \nbrMC^2  + 4 s_\star^2 \sqrt{\clyap_{\chi,0} + \clyap_{\chi,1}} / (n \nbrMC), \\
\eta_1 & = 1 +  12 s_\star \left( \sqrt{\clyap_{\chi,0}} + \sqrt{\clyap_{\chi,1}}  \right)/\nbrMC+ 144 s_\star^2 \clyap_{\chi,1} / \nbrMC^2\\
& \qquad + 4 s_\star^2 \sqrt{\clyap_{\chi,1}} / (n \nbrMC),  \\
\bsf_0  &= 6 s_\star \sqrt{v_{\max}} \sqrt{\clyap_{\chi,0}} / \nbrMC, \\
\bsf_1  & =6 s_\star \sqrt{v_{\max}} \left(\sqrt{\clyap_{\chi,0}}  + 2  \sqrt{\clyap_{\chi,1}} \right)/ \nbrMC,\\
\step_{\max} & = 2 \{ v_{\min} -\bsf_1\} /( \Liplyap \eta_1), \\
\omega_k  & = 2 \{ v_{\min} -\bsf_1  \} - \step_k \Liplyap \eta_1.
\end{align*}

\subsection{Proofs of \Cref{sec:asymptotic-convergence}} \label{app:proof-ODE}

\begin{proof}[Proof of \Cref{theo:almost-sure-cvgce}]
We first establish that the sequence $\sequence{\prm}[k][\nset]$ satisfies \Cref{assum:SA:boundedness} with probability 1. Define
\[
\bsf_{0,k} \eqdef \clyap_{\lyap} \sqrt{\boundbias_{0,k}}/2,  \quad \bsf_{1,k} \eqdef \clyap_{\lyap} (\sqrt{\boundbias_{0,k}} /2 + \sqrt{\boundbias_{1,k}}).
\]
It follows from the Robbins-Siegmund inequality (see \Cref{lem:main-result-quantitative})  that
\begin{align*}
\CPE{\lyap(\prm_{k+1})}{\mcf_k}  & \leq  \lyap(\prm_k)  +  \step_{k+1} \bsf_{0,k} + \step_{k+1}^2 \Liplyap  \eta_0/2 \\
& -    \step_{k+1} \{ \rholyap -  \bsf_{1,k} - \step_{k+1} \Liplyap \, \eta_1 /2 \} \, \superlyap(\prm_k).
\end{align*}
By $\sum_{k=0}^\infty \step_k = \infty$ and  $\sum_{k=0}^\infty \step_k^2 < \infty$, and \ref{item:theo:almost-sure-cvgce:noise},
\begin{equation}\label{eq:SAas:driftpositive}
\lim_{k \to \infty} \{\rholyap -  \bsf_{1,k} - \step_{k+1} \Liplyap \, \eta_1 /2 \}= \rholyap
\end{equation}
exists and is positive. Therefore,  there exists $k_0$ such that for all $k \geq k_0$,
\begin{multline}\label{eq:lyapineq:SAas:trunc}
\CPE{\lyap(\prm_{k+1})}{\mcf_k}  \leq  \lyap(\prm_k)  +  \step_{k+1} \bsf_{0,k} + \step_{k+1}^2 \Liplyap  \eta_0/2  \\
-    \step_{k+1} \{ \rholyap -  \bsf_{1,k} - \step_{k+1} \Liplyap \, \eta_1 /2 \} \, \left( \superlyap(\prm_k) \wedge C \right)   \end{multline}
where $C$ is a positive constant, and
\begin{equation}\label{eq:lyapineq:SAas}
\CPE{\lyap(\prm_{k+1})}{\mcf_k}  \!\leq \! \lyap(\prm_k)  +  \step_{k+1} \bsf_{0,k} + \step_{k+1}^2 \Liplyap  \eta_0/2 \,.
\end{equation}
Define for $k \geq k_0$,
\begin{equation}
\label{eq:def:supermartingale}
M_k \eqdef  \lyap(\prm_k) - \lyap_\star  +  \sum_{\ell=k+1}^\infty \rm{c}_{\ell},
\end{equation}
where $\rm{c}_{\ell+1} \eqdef  \step_{\ell+1} \bsf_{0,\ell} + (\Liplyap  \eta_0/2) \step_{\ell}^2$.
It is easily checked that for all $k \geq k_0$, $M_k \geq 0$  and $\CPE{M_{k+1}}{\mcf_k} \leq M_k$ almost-surely. Hence,
$\sequence{M}[\ell-k_0][\nset]$ is a non-negative supermartingale.
It follows from \cite[Theorems~II-2-7,II-2-9]{neveu1975discrete}  that $\sup_{k \in \nset} M_k < \infty$ and the sequence $\sequence{M}[k][\nset]$ converges with probability one. This implies that,
with probability one, $\sup_{k} \lyap(\prm_k) < \infty$ and $\lim_{k \to \infty} \lyap(\prm_k)$ exists.
By \ref{item:theo:almost-sure-cvgce:lyap}, this yields $\sup_{k \in \nset} \| \prm_k \| < \infty$ with probability one.
Hence,  \Cref{assum:SA:boundedness} holds with probability one.

Note two corollaries of the discussion above. First, from  \eqref{eq:SAas:driftpositive}, \eqref{eq:lyapineq:SAas:trunc}, $\sup_k \lyap(\prm_k) < \infty$, the stepsize condition $\sum_{k=0}^\infty \step_k = \infty$, $\sum_{k=0}^\infty \step_k^2 < \infty$, and \ref{item:theo:almost-sure-cvgce:noise}, we have  when $C \to +\infty$,
\[
\sup_k \superlyap(\prm_k) < \infty
\]
with probability one. Second,  it also follows from  \eqref{eq:lyapineq:SAas} that
\[
\sup_{k \in \nset} \PE[\lyap(\prm_k)] \!\leq \! \PE[\lyap(\prm_0)] + \sum_{k =0}^\infty \step_{k+1} \{ \bsf_{0,k} +  (\Liplyap  \boundvar_0/2) \gamma_{k+1} \};
\]
 the \rhs\ is finite by the stepsize conditions, \ref{item:theo:almost-sure-cvgce:noise}  and the condition $\PE[\lyap(\prm_0)]< \infty$.

We now consider \Cref{assum:SA:noise}. We proved that $\sup_k \superlyap(\prm_k)< \infty$ with probability one. Combined with \ref{item:theo:almost-sure-cvgce:noise}, this yields with probability one,
\begin{multline*}
\| \CPE{\Hg(\prm_{k},\State_{k+1})}{\mcf_{k}} - \hg(\prm_{k}) \|^2 \\
\leq \boundbias_{0,k} + \boundbias_{1,k} \superlyap (\prm_k) \to 0.
\end{multline*}
The asymptotic rate of change condition is a consequence of the convergence theorem  for square integrable martingales (see e.g.~\cite[Theorem 2.15]{hall2014martingale}): the series $\sum_{k=1}^\infty \step_k \Noise_{k+1}$ converges with probability one on the event $\sum_{k=1}^\infty \step^2_k \CPE{\|\Noise_{k+1}\|^2}{\mcf_k} < \infty$.  We thus have to prove that $\sum_k \step_{k+1}^2 \boundvar_{0} +  \sum_k \step_{k+1}^2 \boundvar_{1}  \superlyap(\prm_k) < \infty$ with probability one. This holds true by the stepsize conditions and the property $\sup_k \superlyap(\prm_k) < \infty$.
\end{proof}

\subsection{Proofs of \Cref{subsec:variance-reduction}}
\label{sec:RV:proofslemma}
\begin{proof}[Proof of \Cref{{lem:VR:bias}}]
Let $t \in \{1, \ldots, \kout\}$ and
$k \in \{0, \ldots, \kin-1\}$. Since
$\lbatchrv^{-1} \PE\left[ \sum_{i \in \batchrv_{t,k+1}} a_i \right]=
n^{-1} \sum_{i=1}^n a_i$ (see e.g. \cite[Lemma 7.1.]{fort:moulines:2022}) and $\Hgrv_{t,k} \in \mcf_{t,k}$, then
\[
\CPE{\Hgrv_{t,k+1}}{\mcf_{t,k}} = \Hgrv_{t,k} + n^{-1} \sum_{i=1}^n \left( \hg_i(\prm_{t,k}) - \hg_i(\prm_{t,k-1}) \right).
\]
This concludes the proof of the first claim. The second one follows by
induction, upon noting that $\CPE{{\bf U}}{\mcf_{t,0}}
= \CPE{ \CPE{\bf U}{\mcf_{t,\ell}} }{\mcf_{t,0}}$ for any $ \ell \geq
0$.
\end{proof}

\begin{proof}[Proof of \Cref{lem:VR:varL2error}]
Let $t \in \{1, \ldots, \kout\}$ and $k \in \{0, \ldots, \kin-1\}$. By \Cref{lem:VR:bias}, we have
\begin{align*}
& \Hgrv_{t,k+1} - \CPE{\Hgrv_{t,k+1}}{\mcf_{t,k}} \\ &
= \Hgrv_{t,k+1}- \hg(\prm_{t,k}) -\Hgrv_{t,k} + \hg(\prm_{t,k-1}) \\ &
= \frac{1}{\lbatchrv} \sum_{i \in \batchrv_{k+1}} \Delta_i(\prm_{t,k}, \prm_{t,k-1})
- \frac{1}{n} \sum_{i=1}^n  \Delta_i(\prm_{t,k}, \prm_{t,k-1})
\end{align*}
where $\Delta_i(\prm_{t,k}, \prm_{t,k-1}) \eqdef \hg_i(\prm_{t,k})
- \hg_i(\prm_{t,k-1})$. The conditional expectation of the
$L^2$-moment is upper bounded by $\lbatchrv^{-1}
n^{-1} \sum_{i=1}^n \|\Delta_i(\prm_{t,k}, \prm_{t,k-1}) \|^2$ (see
e.g. \cite[Lemma
7.1.]{fort:moulines:2022}). Under \Cref{assumVR:lipschitz:hi}, we have
$\|\Delta_i(\prm_{t,k}, \prm_{t,k-1}) \|^2 \leq \Liphrv_i^2 \| \prm_{t,k}
- \prm_{t,k-1} \|^2$. This concludes the proof of the first claim. For
the second one, we write $\Hgrv_{t,k+1} - \hg(\prm_{t,k}) = {\bf U_1}
+ {\bf U_2}$ where ${\bf U_1} \eqdef \Hgrv_{t,k+1}
- \CPE{\Hgrv_{t,k+1}}{\mcf_{t,k}}$. By definition of the conditional
expectation and since ${\bf U_2} \in \mcf_{t,k}$, we have $\CPE{\| {\bf
U_1} + {\bf U_2} \|^2}{\mcf_{t,k}} = \CPE{\| {\bf U_1} \|^2}{\mcf_{t,k}}
+ \CPE{\| {\bf U_2} \|^2}{\mcf_{t,k}} = \CPE{\|{\bf
U_1} \|^2}{\mcf_{t,k}} + \| {\bf U_2}\|^2$. The proof is concluded by
using the first statement and \Cref{lem:VR:bias}.
\end{proof}

\begin{proof}[Proof of \Cref{lem:VR:robbins-siegmund-ameliore}]
From \Cref{assum:lyapunov}-\ref{item:lyap}, we write for any $\prm, \bm{d}, \hg \in \rset^d$ and  $\step, \beta >0$,
\begin{multline*}
\lyap(\prm + \step \bm{d}) \leq \lyap(\prm) + \step \pscal{\nabla \lyap(\prm)}{\bm{d}} + \step^2 \frac{\Liplyap}{2} \| \bm{d} \|^2
\\
\leq \lyap(\prm) + \step \pscal{\nabla \lyap(\prm)}{\hg} + \step^2 \Liplyap \left( \| \bm{d} - \hg\|^2 + \| \hg \|^2 \right) \\
+\frac{\step}{2\beta} \| \nabla \lyap(\prm) \|^2 + \frac{\step \beta}{2} \| \bm{d} - \hg \|^2.
\end{multline*}
Applied with $\hg \leftarrow \hg(\prm)$,  and using \Cref{assum:field}-\ref{item:clyap-hg},
\Cref{assum:lyapunov}-\ref{item:rholyap} and \eqref{eq:definition:clyap-lyap}, this inequality implies
\begin{multline*}
\lyap(\prm + \step \bm{d})  \leq  \lyap(\prm) - \step \left( \rholyap - \frac{1}{2\beta} \clyap_\lyap^2 - \step \Liplyap \clyap_{\hg,1}\right) \superlyap(\prm) \\
+ \step\left( \frac{\beta}{2} + \step \Liplyap \right) \| \bm{d}
 - \hg(\prm)\|^2 + \step^2 \Liplyap \clyap_{\hg,0}.
\end{multline*}
We choose $\beta \eqdef \clyap_\lyap^2/\rholyap$ and obtain, for any $\mu > 0$,
\begin{multline*}
\lyap(\prm + \step \bm{d})  \leq  \lyap(\prm) - \step \left( \frac{\rholyap}{2} - \step \Liplyap \clyap_{\hg,1}\right) \superlyap(\prm) + \step^2 \Liplyap \clyap_{\hg,0}\\
- \step  \mu\| \bm{d}
 - \hg(\prm)\|^2  + \step\left( \frac{\clyap_\lyap^2}{2\rholyap} + \mu+ \step \Liplyap \right) \| \bm{d} - \hg(\prm)\|^2 .
\end{multline*}
Applying this inequality with $\step = \step_{t,k+1}$, $\prm
= \prm_{t,k}$, $\bm{d} =\Hgrv_{t,k+1}$, $\mu = \mu_{t,k+1}  \eqdef \rholyap/2
- \step_{t,k+1} \Liplyap$ and computing the conditional expectation
\wrt\ to $\mcf_{t,0}$, concludes the proof.
\end{proof}

\begin{proof}[Proof of \Cref{theo:VR:nonasymptotic}]
By \Cref{eq:VR:controlMSE} and \Cref{assum:field}-\ref{item:clyap-hg}
\begin{align*}
\sum_{k=0}^{\kin-1}& \step_{t,k+1} \CPE{\| \Hgrv_{t,k+1}  - \hg(\prm_{t,k})\|^2}{\mcf_{t,0}} \\
& \leq 2 \frac{\Liphrv^2 \kin}{\lbatchrv} \sum_{k=1}^{\kin}  \step^3_{t,k} \CPE{\| \Hgrv_{t,k} - \hg(\prm_{t,k-1})\|^2}{\mcf_{t,0}} \\
\textstyle & + 2 \frac{\Liphrv^2 \kin}{\lbatchrv} \sum_{k=1}^{\kin}  \step^3_{t,k} 
(\clyap_{\hg,0} + \clyap_{\hg,1} \CPE{\superlyap(\prm_{t,k-1})}{\mcf_{t,0}}).
\end{align*}
We use \Cref{lem:VR:robbins-siegmund-ameliore} and sum from $k=0$ to $k
=\kin-1$:
{\small
\begin{align*}
& \CPE{\lyap(\prm_{t,\kin})}{\mcf_{t,0}}  \leq  \CPE{\lyap(\prm_{t,0})}{\mcf_{t,0}}  + \Liplyap \clyap_{\hg,0} \sum_{k=0}^{\kin-1}\step^2_{t,k+1}  \\
& - \sum_{k=0}^{\kin-1} \step_{t,k+1} \nu_{t,k+1} \CPE{\superlyap(\prm_{t,k})}{\mcf_{t,0}} \\
&- \sum_{k=0}^{\kin-1} \step_{t,k+1} \,  \mu_{t,k+1}  \,  \CPE{\| \Hgrv_{t,k+1} - \hg(\prm_{t,k})\|^2}{\mcf_{t,0}} \\
& + 2 \frac{\Liphrv^2 \kin }{\lbatchrv} \operatorname{a} \sum_{k=0}^{\kin-1}   \step^3_{t,k+1} \CPE{\| \Hgrv_{t,k+1} - \hg(\prm_{t,k})\|^2}{\mcf_{t,0}} \\
& + 2 \frac{\Liphrv^2 \kin }{\lbatchrv}  \clyap_{\hg,0}  \operatorname{a}  \sum_{k=0}^{\kin-1}  \step^3_{t,k+1}  \\
 & + 2 \frac{\Liphrv^2 \kin}{\lbatchrv}  \clyap_{\hg,1}   \operatorname{a} \sum_{k=0}^{\kin-1}  \step^3_{t,k+1}  \CPE{\superlyap(\prm_{t,k})}{\mcf_{t,0}}.
\end{align*}
}
Hence, we get
{\small 
\begin{align*}
&\sum_{k=1}^{\kin}\step_{t,k} \left(\frac{\rholyap}{2} -    \clyap_{\hg,1} \, \step_{t,k} \, \lambda_{t,k} \right) \CPE{\superlyap(\prm_{t,k-1})}{\mcf_{t,0}}  \\
&\quad + \sum_{k=1}^{\kin} \step_{t,k} \left(\frac{\rholyap}{2} -  \step_{t,k}  \, \lambda_{t,k} \right) \CPE{\|\Hgrv_{t,k} - \hg(\prm_{t,k-1}) \|^2}{\mcf_{t,0}}  \\
&\leq   \CPE{\lyap(\prm_{t,0})}{\mcf_{t,0}}  -  \CPE{\lyap(\prm_{t,\kin})}{\mcf_{t,0}} \\
&\quad +  \clyap_{\hg,0}  \left\{  \Liplyap  \sum_{k=1}^{\kin}\step^2_{t,k} +
2 \frac{\Liphrv^2 \kin }{\lbatchrv}
\operatorname{a} \sum_{k=1}^{\kin} \step^3_{t,k} \right\}.
\end{align*}
}
We sum from $t= 1$ to $t = \kout$, and
use \Cref{assum:lyapunov}-\ref{item:vstar}, $\prm_{t-1,\kin}
= \prm_{t,0}$ (see \Cref{line:algo:SPIDERwithinSA:initinner} in \Cref{algo:SPIDERwithinSA}) to conclude the proof.
\end{proof}





\end{document}